\documentclass[12pt]{book}
\usepackage{amsmath,amssymb,amsbsy,amstext,amscd,amsfonts}
\usepackage{exscale,calrsfs,times,multind,theorem}
\usepackage{color}
\usepackage{article}
\usepackage{pictex,epsfig}
\usepackage{multind}
\makeindex{A} \makeindex{N} \makeindex{S}

\makeatletter
\renewcommand{\@oddhead}{ {\em\scriptsize 
		\hfill  Partial Sums and Marcinkiewicz and Fej\'er Means} \hfill}
\renewcommand{\@evenfoot}{\scriptsize\em Hardy spaces
	\hfill  \hfill,}
\renewcommand{\@oddfoot}{\scriptsize\em Hardy spaces
	\hfill  G. Tephnadze \hfill,}

\input article.def

\begin{document}
	
\thispagestyle{empty}

\thispagestyle{empty}

\begin{center}
	{\Large Ivane Javakhishvili Tbilisi State University}\\[4cm]
	
	{\Huge George Tephnadze}\\[20 mm]
	{\large Faculty of Exact and Natural Sciences\\[4 mm]
		Department of Mathematics}\\[20 mm]
	
	{\Large On the Partial Sums and Marcinkiewicz and Fej\'er Means \\ on the One- and Two-dimensional One-parameter Martingale Hardy Spaces}\\[2 cm]
	
	{Georgian PhD Thesis}\\[60 mm]
	
	{Tbilisi 2016}
\end{center}

\newpage

\thispagestyle{empty}

\newpage

\begin{center}
	{\Large \textbf{Abstract}}
\end{center}

\vspace{0.10cm}

Unlike the classical theory of Fourier series which deals with decomposition of a function into continuous waves, the Walsh functions are “rectangular waves”. Such waves have already been used frequently in the theory of signal transmission, codic theory, cryptography, filtering, image enhancement and digital signal processing. \vspace{0.2cm}

The problems we have studied in this PhD thesis are central to Mathematical Analysis. They involve  techniques  which  have  been  developed  a  great  deal  during the last three decades. \vspace{0.2cm}

In this  PhD thesis we are dealing with convergence and summability of partial sums, Fej\'er and Marcinkiewicz means with respect to one- and two-dimensional Walsh-Fourier series on the martingale Hardy spaces. \vspace{0.2cm}

This thesis is focus to achieve the following main results: \vspace{0.2cm}

• To find estimation of convergence and divergence of the subsequences of partial sums of the one-dimensional Walsh-Fourier series on the  martingale Hardy  spaces $H_p(G)$, when  $0<p\leq1$. \vspace{0.2cm}

• To find necessary and sufficient conditions in terms of modulus of continuity of martingale Hardy spaces, for which subsequences of partial sums of the one-dimensional Walsh-Fourier series convergence in $H_p(G)$ norm, when  $0<p\leq1$. \vspace{0.2cm}

• To find estimation of convergence and divergence of the subsequences of  Fej\'er means of the one-dimensional Walsh-Fourier series on the  martingale Hardy  spaces $H_p(G)$, when  $0<p\leq1/2$.  \vspace{0.2cm}

• To find necessary and sufficient conditions in terms  of modulus of continuity of martingale Hardy spaces, for which subsequences of Fej\'er means of the one-dimensional Walsh-Fourier series converge in $H_{p}(G)$ norm, when  $0<p\leq1/2$. \vspace{0.2cm}

• To prove strong convergence of one-dimensional Fej\'er means with respect to Walsh system  on the martingale Hardy  spaces $H_{p}(G)$, when  $0<p\leq 1/2$. \vspace{0.2cm}

• To prove strong convergence of diagonal partial sums with respect  to the two-dimensional Walsh-Fourier series  on the martingale Hardy  spaces $H_{p}(G^2)$, when  $0<p<1$.  \vspace{0.2cm}

• To prove strong convergence of Marcinkiewicz means with respect to the two-dimensional Walsh-Fourier series in $H_{2/3}(G^2)$ norm.  \vspace{0.2cm}

• To find necessary and sufficient conditions in terms of modulus of continuity  of Hardy spaces, for which Marcinkiewicz means of the two-dimensional Walsh-Fourier series converge in $H_{2/3}(G^2)$ norm.  \vspace{0.2cm}

\newpage

\thispagestyle{empty}

\tableofcontents

\thispagestyle{empty}

\newpage

\thispagestyle{empty}


\newpage

\begin{center}
{\Large \textbf{Preface}}
\end{center}

\vspace{0.9cm}

The classical theory of Fourier series deals with decomposition of a
function into sinusoidal waves. Unlike these continuous waves the Vilenkin
(Walsh) functions are rectangular waves. Such waves have already been used
frequently in the theory of signal transmission, multiplexing, filtering,
image enhancement, codic theory, digital signal processing and pattern
recognition. The development of the theory of Vilenkin-Fourier series has
been strongly influenced by the classical theory of trigonometric series.
Because of this it is inevitable to compare results of Vilenkin series to
those on trigonometric series. There are many similarities between these
theories, but there exist differences also. Much of these can be explained
by modern abstract harmonic analysis, which studies orthonormal systems from the point of view of the structure of a topological group. \vspace{0.1cm}

The problems studied in this PhD thesis are very important in Mathematical Analysis and its applications. In particular, we consider convergence and summability of partial sums, Fej\'er and Marcinkiewicz means with respect to the one- and two-dimensional Walsh-Fourier series in the martingale Hardy spaces. According to the problems considered in this PhD thesis, widely are used methods of real analysis combined with methods of abstract and non-linear harmonic analysis together with theory of approximation. Other research methods include theory of function spaces. They involve techniques which have been developed a great deal during the last three decades. \vspace{0.3cm}

This PhD thesis consists the following chapters: \vspace{0.2cm}

•	Preliminaries \vspace{0.2cm}

•   Partial sums with respect to the one-dimensional Walsh-Fourier series on the martingale Hardy spaces \vspace{0.2cm}

•   Fej\'er means means with respect to the one-dimensional Walsh-Fourier series on the martingale Hardy spaces \vspace{0.2cm}

•   Convergence and summability of partial sums with respect to the two-dimensional Walsh-Fourier series on the martingale Hardy spaces \vspace{0.2cm}

In Chapter 1 we first present some classical results and well-known facts, which are very important in the theory of Fourier analysis and is also very crucial for the further investigation of problems considered in this thesis. Moreover, there are presented results which are proved in the next chapters of this thesis and is emphasized the actuality of them, there is also shown some connections with known results. 

\vspace{0.1cm}

In Chapter 2 we first define Walsh group and system, which are important to develop theory of harmonic analysis on locally compact Abelian groups. We consider some expressions and estimations of Lebesgue constants and Dirichlet kernels, give basic definitions and notation of the theory of martingale Hardy spaces and fundamental theorems which are very important to prove main results of this thesis. We also construct martingales which we use to prove sharpness of our positive results. Next, we find rate of convergence and divergence  of the subsequences of partial sums with respect to the one-dimensional Walsh-Fourier series on the martingale Hardy spaces for $0<p\leq1$. Finally, we apply these results to find necessary and sufficient conditions for the modulus of continuity which provide norm convergence of subsequences of the partial sums in the martingale Hardy spaces, for $0<p\leq1$. \vspace{0.1cm}

In Chapter 3 we investigate Fej\'er means with respect to the one-dimensional Walsh-Fourier series. First, we consider some expressions and estimations of Fej\'er kernels and find rate of convergence and divergence  of the subsequences of Fej\'er means in the martingale Hardy spaces for $0<p\leq1/2$. After that, we apply these results to find necessary and sufficient conditions for the modulus of continuity, which provide convergence of subsequences of Fej\'er means in the martingale Hardy spaces. Finally, we prove some new strong convergence theorems of Fej\'er means, for $0<p\leq1/2$. We also prove sharpness of all our main results in this Chapter.\vspace{0.1cm}

In Chapter 4 we investigate basic definitions and notation of partial sums and Marcinkiewicz means with respect to the two-dimensional Walsh-Fourier series. First, we consider some expressions and estimations of Marcinkiewicz kernels, give basic definitions and notation of the theory of the two-dimensional martingale Hardy spaces and fundamental theorems which are very important to prove main results of this thesis. Next, we present and prove strong convergence results of diagonal partial sums with respect to the two-dimensional Walsh-Fourier series in the martingale Hardy spaces for $0<p<1$. Moreover, we consider strong convergence results of Marcinkiewicz means with respect to the two-dimensional Walsh-Fourier series for $p=2/3$ and find necessary and sufficient conditions for the modulus of continuity, which provide convergence in $H_{2/3}(G^2)$ norm of Marcinkiewicz means.\vspace{0.1cm}

This PhD thesis is written as a monograph based on the following
publications: \vspace{0.3cm}

{[1] G. Tephnadze,} Strong convergence of two-dimensional Walsh-Fourier series,  Ukr. Math. J., 65, (6), (2013), 822-834. \vspace{0.3cm}

{[2] G. Tephnadze,} Strong convergence theorems of Walsh-Fejér means, Acta Math. Hung., 142 (1) (2014), 244–259.  \vspace{0.3cm}

{[3] K. Nagy and G. Tephnadze,} Approximation by Walsh-Marcinkiewicz means on the Hardy space  $H_{2/3}$, Kyoto J. Math., 54 (3), (2014), 641-652. \vspace{0.3cm}

{[4] G. Tephnadze,} On the partial sums of Walsh-Fourier series, Colloquium Mathematicum, 141, 2 (2015), 227-242. \vspace{0.3cm}

{[5] G. Tephnadze,} On the convergence of Fejér means of Walsh-Fourier series in the space $ H_{p}$, J. Contemp. Math. Anal., 51, 2 (2016), 51-63. \vspace{0.3cm}

{[6] K. Nagy and G. Tephnadze,} Strong convergence theorem for Walsh-Marcin- kiewicz means, Math. Inequal. Appl., 19, 1 (2016), 185–195.
\vspace{0.3cm}

\newpage

\section{Preliminaries}

\text{ \  \ }\text{ \  \ } It is well-known that (for details see e.g. \cite{gol} and \cite{sws}) for every $p>1$ there exists an absolute constant $ c_{p} $, depending only on $p$, such that
\begin{equation*}
\left\Vert S_{n}f\right\Vert _{p}\leq c_{p}\left\Vert f\right\Vert _{p},
\text{ \ when \ }p>1 \text{ \ and \ } f\in H_{1}({G}).
\end{equation*}

Moreover,  Watari \cite{Wat1} (see also Gosselin \cite{goles} and Young \cite{Yo}) proved that there exists an absolute constant $ c $ such that, for $ n =1,2,..., $
\begin{eqnarray*}
	\lambda \mu
	\left( \vert S_n f\vert>\lambda \right)&\leq& c\left\Vert f\right\Vert_{1}, \ \ \ f\in L_1(G_m), \ \ \lambda>0.
\end{eqnarray*}

On the other hand, it is also well-known that (for details see e.g. \cite{AVD} and \cite{sws}) Walsh system is not Schauder basis in $L_{1}({G})$ space. Moreover, there exists function $f\in H_{1}({G})$, such that partial sums with respect to Walsh system are not uniforml5y bounded in $L_{1}({G})$.

By applying Lebesgue constants
\begin{equation*}
L(n):=\left\Vert D_{n}\right\Vert _{1}
\end{equation*}
we easily obtain that (for details see e.g. \cite{1} and \cite{sws}) subsequences of partial sums $S_{n_{k}}f$ with respect to Walsh system converge to $f$ in $L_{1}$ norm if and only if
\begin{equation} \label{1.0.1}
\sup_{k\in \mathbb{N}}L(n_{k})\leq c<\infty.
\end{equation}

Since $ n $-th Lebesgue constant with respect to Walsh system, where  $$n=\sum_{j=0}^{\infty }n_{j}2^{j}, (n_{j}\in Z_{2})$$  can be estimated by  variation of natural number 
\begin{equation*}
V\left( n\right) =n_{0}+\overset{\infty }{\underset{k=1}{\sum }}\left|
n_{k}-n_{k-1}\right|
\end{equation*}
and it is also well known that  (for details see e.g. \cite{BPT1} and \cite{sws}) the following two-sided estimate is true
\begin{equation*}
\frac{1}{8}V\left( n\right) \leq L(n) \leq
V\left( n\right)
\end{equation*}
to obtain convergence of subsequences of partial sums $S_{n_{k}}f$ with respect to Walsh system of $f\in L_{1}$ in  $f\in L_{1}$-norm.
Condition (\ref{1.0.1}) can be replaced by
\begin{equation*}
\sup_{k\in \mathbb{N}}V(n_{k})\leq c<\infty
\end{equation*}

It follows that (for details see e.g. \ \cite{sws} and \cite{We1}) subsequence of partial sums $S_{2^{n}}$ are bounded from $H_{p}(G)$ to $H_{p}(G)$ for every $p>0$, from which we obtain that
\begin{equation} \label{1.S2n000}
\left\Vert S_{2^{n}}f-f\right\Vert _{H_{p}(G)}\rightarrow 0,\text{ \ as \ }n\rightarrow \infty,
\end{equation}
On the other hand, (see e.g.  \cite{tep7}) there exist a martingale $f\in H_{p}(G)$ $\left( 0<p<1\right),$ such that
\begin{equation*}
\underset{n\in\mathbb{N}}{\sup }\left\Vert S_{2^{n}+1}f\right\Vert
_{weak-L_{p}(G)}=\infty.
\end{equation*}%

The main reason of divergence of subsequence  $S_{{2^n}+1}f$ of partial sums it that  (for details see \cite{tep9}) Fourier coefficients of $f\in H_{p}(G)$ are not uniformly bounded when $0<p<1$.

When $0<p<1$ in \cite{tep_thesis} was investigated boundedness of  subsequences of partial sums with respect to Walsh system from $H_{p}(G)$ to $H_{p}(G)$. In particular, the following result is true:

\textbf{Theorem T1.} Let $0<p<1$ and $f\in H_{p}(G)$. Then there exists a absolute constant $c_{p},$ depending only on $p$, such that
\begin{equation*}
\left\Vert S_{m_{k}}f\right\Vert _{H_{p}(G)}\leq c_{p}\left\Vert f\right\Vert _{H_{p}(G)}
\end{equation*}
if and only if the following condition holds
\begin{equation} \label{1.cond001}
\underset{k\in \mathbb{N}}{\sup } d\left( m_{k}\right)<c<\infty,
\end{equation}
where
\begin{equation*}
d\left( m_{k}\right) :=\left\vert m_{k}\right\vert -\left\langle m_{k}\right\rangle.
\end{equation*}

In particular, Theorem T1 immediately follows:

\textbf{Theorem T2.} Let $p>0$ and $f\in H_{p}(G)$. Then there exists a absolute constant $c_{p},$ depending only on $p$, such that
\begin{equation*}
\left\Vert S_{2^{n}}f\right\Vert_{H_p(G)}\leq c_{p}\left\Vert f\right\Vert _{H_p(G)}
\end{equation*}
and
\begin{equation*}
\left\Vert S_{2^{n}+2^{n-1}}f\right\Vert_{H_p(G)}\leq c_{p}\left\Vert f\right\Vert _{H_p(G)}.
\end{equation*}

On the other hand, we have the following result:

\textbf{Theorem T3.}  Let $p>0$. Then there exists a martingale $ f \in H_{p}(G) $, such that
$$  \sup_{n\in \mathbb{N}}\left\Vert S_{2^{n}+1}f\right\Vert_{_{H_{p}(G)}}=\infty. $$

Taking into account these results it is interesting to find behaviour of rate of divergence of subsequences of partials sums with respect to Walsh system of martingale $f\in H_{p}(G)$ in the martingale Hardy spaces $H_{p}(G)$.

In the second chapter of this thesis (see also \cite{tep12}) we investigate ebove mentioned problem. For $0<p<1$ we have the following result:

\textbf{Theorem \ref{th4.1.1}.} Let $f\in H_{p}(G)$. Then there exists a absolute constant $c_{p},$ depending only on $p$, such that the following inequality is true
\begin{equation} \label{1.vnk00002}
\text{ }\left\Vert S_{n}f\right\Vert _{H_{p}(G)}\leq c_{p}2^{d\left( n\right)\left( 1/p-1\right) }\left\Vert f\right\Vert _{H_{p}(G)}.
\end{equation}

On the other hand, if  $0<p<1,$ $\left\{ m_{k}:\text{ }k\geq 0\right\} $ be increasing subsequence of natural numbers, such that
\begin{equation}  \label{1.dmk}
\sup_{k\in \mathbb{N}}d\left( m_{k}\right) =\infty
\end{equation}%
and
$\Phi :\mathbb{N}_{+}\rightarrow \lbrack 1,\infty )$ be non-decreasing function satisfying the condition
\begin{equation*} \label{1010}
\overline{\underset{k\rightarrow \infty }{\lim }}\frac{2^{d\left(
m_{k}\right) \left( 1/p-1\right) }}{\Phi \left( m_{k}\right) }=\infty,
\end{equation*}%
then there exists a martingale $f\in H_{p}(G),$ such that
\begin{equation*}
\underset{k\in \mathbb{N}}{\sup }\left\Vert \frac{S_{m_{k}}f}{\Phi \left(m_{k}\right) }\right\Vert _{weak-L_{p}(G)}=\infty .
\end{equation*}

Theorem \ref{cor4.1.1} easily follows the following corollary:

\textbf{Corollary \ref{cor4.1.1}.} Let $0<p<1$ and $f\in H_{p}(G)$. Then there exists a absolute constant $c_{p},$ depending only on $p$, such that
\begin{equation*}
\text{ }\left\Vert S_{n}f\right\Vert _{H_{p}(G)}\leq c_{p}\left( n\mu \left\{\text{supp}\left( D_{n}\right) \right\} \right) ^{1/p-1}\left\Vert f\right\Vert _{H_{p}(G)}.
\end{equation*}

On the other hand, if $0<p<1$ and $\left\{ m_{k}:\text{ }k\geq 0\right\} $ be increasing sequence of natural numbers, such that
\begin{equation*}
\sup_{k\in \mathbb{N}}m_{k}\mu \left\{ \text{supp}\left( D_{m_{k}}\right)
\right\} =\infty
\end{equation*}%
and $\Phi :\mathbb{N}_{+}\rightarrow \lbrack 1,\infty )$
be non-decreasing function satisfying the condition
\begin{equation*}  \label{12e}
\overline{\underset{k\rightarrow \infty }{\lim }}\frac{\left( m_{k}\mu\left\{ \text{supp}\left( D_{m_{k}}\right) \right\} \right) ^{1/p-1}}{\Phi
\left( m_{k}\right) }=\infty,
\end{equation*}%
then there exists a martingale $f\in H_{p}(G),$  such that
\begin{equation*}
\underset{k\in \mathbb{N}}{\sup }\left\Vert \frac{S_{m_{k}}f}{\Phi \left(
m_{k}\right) }\right\Vert _{weak-L_{p}(G)}=\infty .
\end{equation*}

In particular, we also get the proofs of  Theorem T1 and Theorem T2.

In the second chapter of this thesis we also investigate case $p=1$. In this case the following result is true:

\textbf{Theorem \ref{th4.1.2}.} Let $n\in \mathbb{N}_{+}$ and $f\in H_{1}(G).$ Then there exists a absolute constant $c,$ such that
\begin{equation} \label{1.vnk00001}
\left\Vert S_{n}f\right\Vert _{H_{1}(G)}\leq cV\left( n\right) \left\Vert
f\right\Vert _{H_{1}(G)}.
\end{equation}

Moreover, if $\left\{ m_{k}:\text{ }k\geq 0\right\} $ be increasing sequence of natural numbers $\mathbb{N}_{+},$ such that
\begin{equation*}
\sup_{k\in \mathbb{N}}V\left( m_{k}\right) =\infty  \label{vnk}
\end{equation*}
and $\Phi :\mathbb{N}_{+}\rightarrow \lbrack 1,\infty )$
be non-decreasing function satisfying the condition
\begin{equation*}
\overline{\underset{k\rightarrow \infty }{\lim }}\frac{V\left( m_{k}\right)}{\Phi \left( m_{k}\right) }=\infty .
\end{equation*}
Then there exists a martingale $f\in H_{1}(G),$ such that
\begin{equation*}
\underset{{k\in \mathbb{N}}}{\sup }\left\Vert \frac{S_{m_{k}}f}{\Phi \left(m_{k}\right) }\right\Vert _{1}=\infty .
\end{equation*}

When $0<p<1$ in \cite{tep_thesis} was proved boundedness of maximal operators of subsequences of partial sums from $H_{p}(G)$ to $L_{p}(G)$. In particular, the following is true:

\textbf{Theorem T4.} Let $0<p<1$ and $f\in H_{p}(G)$. Then the maximal operator
\begin{equation*}
\underset{k\in \mathbb{N}}{\sup }|S_{m_{k}}f| \label{cond00}
\end{equation*}
is bounded from $H_{p}(G)$ to $L_{p}(G)$, if and only if condition (\ref{1.cond001}) is fulfilled.

In the special cases we obtain that the following is true:

\textbf{Theorem T5.} Let $p>0$ and $f\in H_{p}(G)$. Then there exists an absolute constant $c_{p},$ depending only on $p$, such that
\begin{equation}\label{1.S2n}
\left\Vert  \sup_{n\in \mathbb{N}}|S_{2^{n}}f|\right\Vert_{p}\leq c_{p}\left\Vert f\right\Vert _{H_p(G)}
\end{equation}
and
\begin{equation*}
\left\Vert  \sup_{n\in \mathbb{N}}|S_{2^{n}+2^{n-1}}f|\right\Vert_{p}\leq c_{p}\left\Vert f\right\Vert _{H_p(G)}.
\end{equation*}

On the other hand we have the following result:

\textbf{Theorem T6.}  Let $p>0$. Then there exists a martingale $ f \in H_{p}(G) $, such that
$$ \left\Vert  \sup_{n\in \mathbb{N}}|S_{2^{n}+1}f|\right\Vert_{p}=\infty. $$

Above mentioned condition (\ref{1.cond001}) is sufficient condition for the case $p=1$ also, but there exist subsequences which do not satisfy this condition, but maximal operators of these subsequences of partial sums with respect to Walsh system  are not bounded from  $H_{1}(G)$ to  $L_{1}(G)$.

Such necessary and sufficient conditions which provides boundedness of maximal operators of subsequences of partial sums with respect to Walsh system from $H_{1}(G)$ to  $L_{1}(G)$ is open problem.

In \cite{tep9} and \cite{tep_thesis} was investigated boundedness of weighted maximal operators from $H_{p}(G)$ to $L_{p}(G)$, when $ 0<p\leq 1$:

\textbf{Theorem T7.} Let $0<p\leq 1$. Then weighted maximal operator
\begin{equation*}
\overset{\sim }{S}_{p}^{\ast }f:=\underset{n\in \mathbb{N_+}}{%
\sup }\frac{\left\vert S_{n}f\right\vert }{\left( n+1\right) ^{1/p-1}\log^{[p]}\left( n+1\right)}
\end{equation*}
is bounded from  $H_{p}(G)$ to  $L_{p}(G)$, where $ [p] $ denotes integer part of $ p $.

Moreover, for any non-decreasing function  $\varphi :\mathbb{N}_{+}\rightarrow [1,$ $\infty )$ satisfying the condition
\begin{equation*}
\overline{\lim_{n\rightarrow \infty }}\frac{\left( n+1\right) ^{1/p-1}\log^{[p]}\left( n+1\right)}{%
\varphi \left( n+1\right)}=+\infty,
\end{equation*}%
there exists a martingale $f \in H_{p}(G)$ $ (0<p\leq 1), $ such that
\begin{equation*}
\sup_{n\in \mathbb{N}}\left\Vert \frac{S_{n}f}{\varphi \left( n\right)
}\right\Vert _{p}=\infty.
\end{equation*}

According to negative result for weighted maximal operator of partial sums of Walsh-Fourier series we immediately get the following result:

\textbf{Theorem S1.}  There exists a martingale $f\in H_{p}(G)$, $(0<p\leq1)$, such that
\begin{equation*}
\underset{n\in \mathbb{N}}{\sup }\left\Vert S_{n}f\right\Vert
_{p}=\infty .
\end{equation*}%

On the other hand, boundedness of weighted maximal operators immediately follows the following estimation:

\textbf{Theorem S2.} Let $0<p\leq1$. Then there exists a absolute constant $c_{p}$, depending only on $p$, such that 
\begin{equation*}
\left\Vert S_{n}f\right\Vert _{p}\leq c_{p} {(n+1)^{1/p-1}\log^{[p]} (n+1)}\left\Vert f\right\Vert _{H_p(G)},
\text{ \ for \ } 0<p\leq 1,
\end{equation*}
where $ [p] $ denotes integer part of $ p $.

By applying this inequality (see \cite{tep6}) we find necessary and sufficient conditions for martingale $f\in H_{p}(G)$ for which partial sums with respect to Walsh system of martingale $f\in H_{p}(G)$ converge in $H_{p}(G)$ norm.

\textbf{Theorem T8.} Let $0<p\leq 1$,  $[p]$ denotes integer part of $p$, $f\in H_{p}(G)$ and
\begin{equation*}
\omega _{H_{p}(G)}\left( \frac{1}{2^{N}},f\right)=o\left( \frac{1}{2^{N(1/p-1)}N^{[p]}}\right) ,\text{ \ as \ }N\rightarrow
\infty.
\end{equation*}%
Then
\begin{equation*}
\left\Vert S_{n}f-f\right\Vert _{p}\rightarrow 0,\text{ as }
n\rightarrow \infty.
\end{equation*}

Moreover, there exists a martingale $f\in H_{p}(G)$, where $0<p<1$, such that
\begin{equation*}
\omega _{H_{p}(G)} \left( \frac{1}{2^{N}},f\right)=O\left( \frac{1}{2^{N(1/p-1)}N^{[p]}}\right) ,\text{ \ as \ }N\rightarrow \infty
\end{equation*}%
and
\begin{equation*}
\left\Vert S_{n}f-f\right\Vert _{weak-L_{p}(G)}\nrightarrow 0,\,\,\,\text{ \ \ as \ \ }n\rightarrow \infty .
\end{equation*}%

By taking these results into account, it is interesting to find necessary and sufficient conditions for modulus of continuity, such that subsequences of partial sums with respect to Walsh system of martingale $f\in H_{p}(G)$ converge in $H_{p}(G)$ norm.

In the second chapter of this thesis (see also \cite{tep12}) we investigate this problem. By combining inequalities (\ref{1.vnk00002}) and (\ref{1.vnk00001}) we get the following theorem:

\textbf{Theorem \ref{theorem4.2.1}.}  Let $2^{k}<n\leq 2^{k+1}.$ Then there exists an absolute constant $c_{p},$ depending only on $p$, such that
\begin{equation} \label{1.sn3}
\left\Vert S_{n}f-f\right\Vert _{H_{p}(G)}\leq c_{p}2^{d\left( n\right) \left(1/p-1\right) }\omega _{H_{p}(G)}\left( \frac{1}{2^{k}},f\right) ,\text{ \ \ \ }\left( 0<p<1\right)
\end{equation}
and
\begin{equation}  \label{1.sn2}
\left\Vert S_{n}f-f\right\Vert _{H_{1}(G)}\leq c_{1}V\left( n\right) \omega_{H_{1}(G)}\left( \frac{1}{2^{k}},f\right)
\end{equation}

By applying inequality (\ref{1.sn3}) the following result is proved in the second chapter:

\textbf{Theorem \ref{th4.2.2}.} Let $0<p<1,$ $f\in H_{p}(G)$ and $\{m_{k}:k\geq 0\}$ be increasing sequence of natural number satisfying the condition
\begin{equation*}
\omega _{H_{p}(G)}\left( \frac{1}{2^{\left\vert m_{k}\right\vert }},f\right)
=o\left( \frac{1}{2^{d\left( m_{k}\right) \left( 1/p-1\right) }}\right)
\text{ as \ }k\rightarrow \infty.  \label{1.18a}
\end{equation*}%
Then
\begin{equation} \label{1.con1}
\left\Vert S_{m_{k}}f-f\right\Vert _{H_{p}(G)}\rightarrow 0\text{ \ as \ }k\rightarrow \infty .
\end{equation}

On the other hand, if $\{m_{k}:k\geq 0\}$ be increasing sequence of natural numbers satisfying the condition (\ref{1.dmk}), then there exists a martingale $f\in H_{p}(G)$ and subsequence $\{\alpha _{k}:k\geq 0\}\subset \{m_{k}:k\geq 0\},$ for which
\begin{equation*}
\omega _{H_{p}(G)}\left( \frac{1}{2^{\left\vert \alpha _{k}\right\vert }},f\right) =O\left( \frac{1}{2^{d\left( \alpha _{k}\right) \left(
1/p-1\right) }}\right) \text{ as \ }k\rightarrow \infty \text{\ }
\end{equation*}%
and
\begin{equation}
\limsup\limits_{k\rightarrow \infty }\left\Vert S_{\alpha_{k}}f-f\right\Vert _{weak-L_{p}(G)}>c_{p}>0\,\,\,\text{as\thinspace
\thinspace \thinspace }k\rightarrow \infty,  \label{1.con11}
\end{equation}
where $c_{p}$ is an absolute constant depending only on $p$.

According to this theorem we immediately get that the following result is true:

\textbf{Corollary \ref{cor4.2.1}.} Let $0<p<1,$ $f\in H_{p}(G)$ and $\{m_{k}:k\geq 0\}$  be increasing sequence of natural number, satisfying the condition
\begin{equation*} \label{1.cond2}
\omega _{H_{p}(G)}\left( \frac{1}{2^{\left\vert m_{k}\right\vert }},f\right)=o\left( \frac{1}{\left( m_{k}\mu \left( \text{supp}D_{m_{k}}\right) \right)^{1/p-1}}\right), \text{ \ as \ }k\rightarrow \infty.
\end{equation*}%
Then (\ref{1.con1}) holds.

On the other hand, if $\{m_{k}:k\geq 0\}$ be increasing sequence of natural number, satisfying the condition
\begin{equation*}
\overline{\underset{k\rightarrow \infty}{\lim }}\frac{\left( m_{k}\mu
\left\{\text{supp}\left(D_{m_{k}}\right)\right\}\right)^{1/p-1}}{\Phi
\left(m_{k}\right)}=\infty,
\end{equation*}%
then there exists a martingale $f\in H_{p}(G)$ and subsequence $\{\alpha _{k}:k\geq 0\}\subset \{m_{k}:k\geq 0\}$ such that
\begin{equation*}
\omega _{H_{p}(G)}\left( \frac{1}{2^{\left\vert \alpha _{k}\right\vert }},f\right) =O\left( \frac{1}{\left( \alpha _{k}\mu \left( \text{supp}%
D_{\alpha _{k}}\right) \right) ^{1/p-1}}\right), \text{\ \ as \ }k\rightarrow\infty
\end{equation*}%
and (\ref{1.con11}) holds.

By applying (\ref{1.sn2}) we prove that the following is true:

\textbf{Theorem \ref{th4.2.3}.} Let $f\in H_{1}(G)$ and $\{m_{k}:k\geq 0\}$ be increasing sequence of natural number, satisfying the condition
\begin{equation*}
\omega _{H_{1}(G)}\left( \frac{1}{2^{\left\vert m_{k}\right\vert }},f\right)=o\left( \frac{1}{V\left( m_{k}\right) }\right) \text{ \ as \ }k\rightarrow
\infty.
\end{equation*}%
Then
\begin{equation*}
\left\Vert S_{m_{k}}f-f\right\Vert _{H_{1}(G)}\rightarrow 0\text{ \ as \ }
k\rightarrow \infty.
\end{equation*}

Moreover, if $\{m_{k}:k\geq 0\}$ be increasing sequence of natural number, satisfying the condition (\ref{1.dmk}), then there exists a martingale $f\in H_{1}(G)$ and subsequence $\{\alpha _{k}:k\geq 0\}\subset \{m_{k}:k\geq 0\}$ for which
\begin{equation*}
\omega _{H_{1}(G)}\left( \frac{1}{2^{\left\vert \alpha _{k}\right\vert }}%
,f\right) =O\left( \frac{1}{V\left( \alpha _{k}\right) }\right) \text{ \ as
\ }k\rightarrow \infty
\end{equation*}%
and
\begin{equation*} \label{1.cond10}
\limsup\limits_{k\rightarrow \infty}\left\Vert S_{\alpha_k} f-f\right\Vert_1>c>0\text{ \ \ as \ \  }
k\rightarrow \infty, 
\end{equation*}%
where  $c$ is an absolute constant.

By applying Theorem \ref{th4.2.2} and Theorem \ref{th4.2.3} we immediately get proof of Theorem T8.

Weisz \cite{We3} consider convergence in norm of Fej\'er means of the one-dimensional Walsh-Fourier and proved the following:

\textbf{Theorem We1.} Let $p>1/2$ and $f\in H_{p}(G)$. Then there exists a absolute constant $c_{p}$, depending only on $p$, such that 
\begin{equation*}
\left\Vert \sigma _{k}f\right\Vert _{H_{p}(G)}\leq c_{p}\left\Vert f\right\Vert_{H_{p}(G)}.
\end{equation*}

Weisz (for details see e.g. \cite{We1}) also consider boundedness of subsequences of Fej\'er means $\sigma_{2^{n}}$ of the one-dimensional Walsh-Fourier series  from $H_{p}(G)$ to $H_{p}(G)$ when $p>0$:

\textbf{Theorem We2.} Let $p>0$ and $f\in H_{p}(G)$. Then
\begin{equation} \label{fe22222}
\left\Vert \sigma_{2^{k}}f-f\right\Vert _{H_{p}(G)}\rightarrow 0,\text{\ \ as \ \ }k\rightarrow\infty.
\end{equation}

On the other hand, in \cite{tep3} was proved the following result:

\textbf{Theorem T9.} There exists a martingale $f\in H_{p}(G)$ $\left( 0<p\leq 1/2\right)$ such that
\begin{equation*}
\underset{n\in\mathbb{N}}{\sup}\left\Vert\sigma_{2^{n}+1}f\right\Vert_{H_p(G)}=\infty.
\end{equation*}

Goginava \cite{gog16} (see also \cite{pt1}) proved that the following result is true:

\textbf{Theorem Gog1.} Let $0<p\leq 1.$ Then the sequence of operators $\left\vert \sigma _{2^{n}}f\right\vert $ are not bounded from $H_{p}(G)$ to $H_{p}(G)$.

When $ 0<p<1/2 $ then in \cite{pt2} was proved bondedness of subsequences of Fej\'er means of the one-dimensional Walsh-Fourier from $H_{p}(G)$ to $H_{p}(G)$. In particular, the following is true:

\textbf{Theorem T10.} Let $0<p<1/2$ and $f\in H_{p}(G)$. Then there exists a absolute constant $c_{p}$, depending only on $p$, such that 
\begin{equation*}
\left\Vert \sigma_{m_{k}}f\right\Vert _{H_{p}(G)}\leq c_{p}\left\Vert f\right\Vert _{H_{p}(G)}
\end{equation*}
estimation holds if and only if the condition (\ref{1.cond001}) is fulfilled.

Theorem T10 immediately follows theorem of Weisz (see Theorem We2) and and also interesting results:

\textbf{Theorem T11.} Let $p>0$ and $f\in H_{p}(G)$. Then there exists an absolute constant $c_{p},$ depending only on $p$, such that
\begin{equation*}
\left\Vert \sigma_{2^{n}}f\right\Vert_{H_p(G)}\leq c_{p}\left\Vert f\right\Vert _{H_p(G)}
\end{equation*}
and
\begin{equation*}
\left\Vert \sigma_{2^{n}+2^{n-1}}f\right\Vert_{H_p(G)}\leq c_{p}\left\Vert f\right\Vert _{H_p(G)}.
\end{equation*}

On the other hand, we have the following result:

\textbf{Theorem T12.}  Let $p>0$. Then there exists a martingale $ f \in H_{p}(G) $, such that
$$  \sup_{n\in \mathbb{N}}\left\Vert \sigma_{2^{n}+1}f\right\Vert_{_{H_{p}(G)}}=\infty. $$

According to above mentioned results it is interesting to find rate of divergence of subsequences  $\sigma_{n_{k}}f$ of Fej\'er means of the one-dimensional Walsh-Fourier series in the Hardy spaces $H_{p}(G)$.

In the third chapter of this thesis (see also \cite{tep13}) we find rate of divergence of subsequences of Fej\'er means of the one-dimensional Walsh-Fourier series on the martingale Hardy spaces $H_{p}(G)$, when
$ 0<p\leq1/2 $.

First, we consider case $p=1/2$:

\textbf{Theorem \ref{th5.1.1}.} Let $n\in \mathbb{N}_{+}$ and $f\in H_{1/2}(G).$ Then there exists an absolute constant $c,$ such that
\begin{equation} \label{vnk000f}
\left\Vert \sigma_{n}f\right\Vert _{H_{1/2}(G)}\leq c{{V}^2}\left( n\right) \left\Vert f\right\Vert _{H_{1/2}(G)}.
\end{equation}

Moreover, if $\left\{ m_{k}:\text{ }k\geq 0\right\} $ be increasing secuence of natural numbers, such that
\begin{equation*}
\sup_{k\in \mathbb{N}}V\left( m_{k}\right) =\infty  \label{vnkf}
\end{equation*}
and $\Phi :\mathbb{N}_{+}\rightarrow [1,\infty]$
be non-decreasing function satisfying the conditions
\begin{equation*}
\overline{\underset{k\rightarrow \infty }{\lim }}\frac{{V^2}\left( m_{k}\right)}{\Phi \left( m_{k}\right) }=\infty .  \label{17aaf}
\end{equation*}
then there exists a martingale $f\in H_{1/2}(G),$ such that
\begin{equation*}
\underset{{k\in \mathbb{N}}}{\sup }\left\Vert \frac{\sigma_{m_{k}}f}{\Phi \left(m_{k}\right) }\right\Vert _{1/2}=\infty.
\end{equation*}

There was also considered case $0<p<1/2$ and was proved that the following is true:

\textbf{Theorem \ref{th5.1.2}.} Let $0<p<1/2$ and $f\in H_{p}(G)$. Then there exists an absolute constant $c_{p},$ depending only on $p$ such that
\begin{equation} \label{vnk0000f}
\text{ }\left\Vert \sigma_{n}f\right\Vert _{H_{p}(G)}\leq c_{p}2^{d\left( n\right)\left( 1/p-2\right) }\left\Vert f\right\Vert _{H_{p}(G)}.
\end{equation}

On the other hand, if  $0<p<1/2,$ $\left\{ m_{k}:\text{ }k\geq 0\right\} $ be increasing sequence of natural numbers satisfying the condition (\ref{1.dmk})  and
$\Phi :\mathbb{N}_{+}\rightarrow \lbrack 1,\infty )$ be non-decreasing function such that
\begin{equation*}
\overline{\underset{k\rightarrow \infty }{\lim }}\frac{2^{d\left(m_{k}\right) \left( 1/p-2\right) }}{\Phi \left( m_{k}\right) }=\infty,
\end{equation*}%
then there exists a martingale $f\in H_{p}(G),$ such that
\begin{equation*}
\underset{k\in \mathbb{N}}{\sup }\left\Vert \frac{\sigma_{m_{k}}f}{\Phi \left(m_{k}\right) }\right\Vert _{weak-L_{p}(G)}=\infty .
\end{equation*}

From these results also follows proof of Theorem We2.

In 1975 Schipp \cite{Sc} (see also \cite{zy}) proved that the maximal operator of  Fej\'er means $\sigma ^{*}$ is of type weak-(1,1):
\begin{equation*}
\mu \left( \sigma ^{*}f>\lambda \right) \leq \frac{c}{\lambda }\left\|f\right\| _{1},\text{ \qquad }\left( \lambda >0\right).
\end{equation*}
By using  Marcinkiewicz interpolation theorem it follows that  $\sigma ^{*}$ is of strong type-$ (p,p) $, when $p>1:$
\begin{equation*}
\left\| \sigma ^{*}f \right\|_{p} \leq c\left\|
f\right\| _{p}, \text{ \qquad }\left( p >1\right).
\end{equation*}
The boundedness does not hold for  $ p=1,$ but Fujji \cite{Fu} (see also \cite{Yano}) proved that maximal operator of  Fej\'er means is bounded from  $H_{1}(G)$ to $L_{1}(G)$. Weisz in \cite{We2} generalized result of Fujii and proved that maximal operator of  Fej\'er means is bounded from $H_{p}(G)$ to $L_{p}(G)$, when $ p>1/2. $ Simon \cite{Si1} construct the counterexample, which shows that boundedness does not hold when  $0<p<1/2$. Goginava \cite{GoAMH} (see also \cite{BGG} and \cite{BGG2}) generalized this result for   $0<p\leq1/2$ and proved that the following is true:

\textbf{Theorem Gog2.} There exists a martingale $f\in H_{p}(G)$ $\left( 0<p\leq1/2\right)$ such that
\begin{equation*}
\underset{n\in \mathbb{N}}{\sup }\left\Vert \sigma_{n}f\right\Vert_{p}=\infty .
\end{equation*}%

Weisz \cite{we4} (see also Goginava \cite{GoPubl}) proved that the following is true:

\textbf{Theorem We3.} Let $f\in H_{1/2}(G)$. Then there exists an absolute constant $c,$ such that
\begin{equation*}
\left\|\sigma^{*}f\right\|_{weak-L_{1/2}(G)}\leq c\left\|f\right\|_{H_{1/2}(G)}.
\end{equation*}

In \cite{pt2} was considered boundedness of maximal operators of subsequences of Fej\'er means of the one-dimensional Walsh-Fourier series from $H_{p}(G)$ to $L_{p}(G)$ for  $0<p<1/2$. In particular, the following is true:

\textbf{Theorem T13.} Let $0<p<1/2$ and $f\in H_{p}(G)$. Then the maximal operator 
\begin{equation*}
\overset{\sim }{\sigma}^{*}f:=\underset{k\in \mathbb{N}}{\sup }|\sigma_{m_{k}}f|
\end{equation*}
is bounded from $H_{p}(G)$ to $L_{p}(G)$ if and only if when condition (\ref{1.cond001}) is fulfilled.

As consequences the following results are true:

\textbf{Theorem T14.} Let $p>0$ and $f\in H_{p}(G)$. Then there exists an absolute constant $c_p$ depending only on $ p, $ such that
\begin{equation}\label{sigmamax}
\left\Vert  \sup_{n\in \mathbb{N}}|\sigma_{2^{n}}f|\right\Vert_{p}\leq c_{p}\left\Vert f\right\Vert _{H_p(G)}
\end{equation}
and
\begin{equation*}
\left\Vert  \sup_{n\in \mathbb{N}}|\sigma_{2^{n}+2^{n-1}}f|\right\Vert_{p}\leq c_{p}\left\Vert f\right\Vert _{H_p(G)}.
\end{equation*}

On the other hand, we have the following negative result:

\textbf{Theorem T15.}  Let $0<p<1/2.$ Then there exists a martingale $ f \in H_{p}(G) $, such that
$$ \left\Vert  \sup_{n\in \mathbb{N}}|\sigma_{2^{n}+1}f|\right\Vert_{p}=\infty. $$

above mentioned condition is sufficient for the case $p=1/2$ also, but there exists subsequences, which do not satisfy condition (\ref{1.cond001}), but maximal operator of subsequences of Fej\'er means of the one-dimensional Walsh-Fourier series are bounded from $H_{1/2}(G)$ to $L_{1/2}(G).$

However, it is open problem to find necessary and sufficient conditions on the indexes, which provide boundedness of maximal operator of subsequences of Fej\'er means of the one-dimensional Walsh-Fourier series from $H_{1/2}(G)$ to $L_{1/2}(G)$.

In \cite{GoSzeged} and \cite{tep3} (see also \cite{ptw4}, \cite{tepkack4}, \cite{GNCz} and \cite{tep2}) is proved that the following is true:

\textbf{Theorem GT1.} Let $0<p\leq 1/2$ and $ f \in H_{p}(G) $. Then the maximal operator
\begin{equation*}
\overset{\sim }{\sigma }_{p}^{\ast }f:=\underset{n\in \mathbb{N}}{\sup }\frac{\left\vert \sigma _{n}f\right\vert }{\left( n+1\right) ^{1/p-2}\log ^{2[1/2+p]}\left( n+1\right)}
\end{equation*}
is bounded from  $H_{p}(G)$ to  $L_{p}(G)$.

Moreover, for any nondecreasing function  $\varphi :\mathbb{N}_{+}\rightarrow [1,$ $\infty )$ satisfying the condition
\begin{equation*}
\overline{\lim_{n\rightarrow \infty }}\frac{\left( n+1\right) ^{1/p-2}\log ^{2[1/2+p]}\left( n+1\right)}{\varphi \left( n\right) }=+\infty ,
\end{equation*}%
there exists a martingale $f \in H_{p}(G)$, $ (0<p<1/2) $  such that
\begin{equation*}
\sup_{n\in \mathbb{N}}\left\Vert \frac{\sigma _{n}f}{\varphi \left( n\right)}\right\Vert _{p}=\infty.
\end{equation*}

From the divergence of weighted maximal operators we immediately get that there exists a martingale $f\in H_{p}(G)$ $(0<p\leq 1/2)$, such that
\begin{equation*}
\underset{n\in \mathbb{N}}{\sup }\left\Vert \sigma_{n}f\right\Vert
_{p}=\infty.
\end{equation*}
and from the boundedness results of weighted maximal operators we immediately get that for any $f\in H_{p}(G)$ there exists an absolute constant $c_{p}$, such that the following inequality holds true:
\begin{equation} \label{1.cond33}
\left\Vert \sigma _{n}f\right\Vert _{p}\leq c_{p} {n^{1/p-2}\log ^{2[1/2+p]}\left( n+1\right)}\left\Vert f\right\Vert _{H_{p}(G)},\text{ \ as \ } 0<p\leq1/2.
\end{equation}

By applying inequality (\ref{1.cond33}) in \cite{tep6} was found necessary and sufficient conditions for modulus of continuity of martingale $f\in H_{p}(G)$, for which Fej\'er means of the one-dimensional Walsh-Fourier series converge in $H_{p}(G)$ norm.

\textbf{Theorem T16.} Let $ 0<p\leq 1/2 $, $f\in H_{p}(G)$ and
\begin{equation*}
\omega _{H_{p}(G)}\left( \frac{1}{2^{N}},f\right) =o\left(\frac{1}{2^{N(1/p-2)}N^{2[1/2+p]}}\right),\text{\  as \ } N\rightarrow\infty.
\end{equation*}%
Then
\begin{equation*}
\left\Vert\sigma_{n}f-f\right\Vert _{p}\rightarrow 0,\text{ as }n\rightarrow \infty .
\end{equation*}

Moreover, there exists a martingale $f\in H_{p}(G)$, for which
\begin{equation*}
\omega _{H_{1/2}(G)}\left( \frac{1}{2^{N}},f\right) =O\left( \frac{1}{2^{N(1/p-2)}N^{2[1/2+p]}}\right) ,\text{ \ as \ }N\rightarrow \infty
\end{equation*}%
and
\begin{equation*}
\left\Vert\sigma_{n}f-f\right\Vert_{p}\nrightarrow 0,\,\,\,\text{as\ \ }n\rightarrow \infty.
\end{equation*}

According above mentioned results, it is interesting to find necessary and sufficient conditions for the modulus of continuity, for which subsequences $\sigma_{n_{k}}f$ of Fej\'er means of the one-dimensional Walsh-Fourier series converge in $H_{p}(G)$ norm.

In the third chapter of this thesis we find necessary and sufficient conditions for the modulus of continuity, for which subsequences $\sigma_{n_{k}}f$ of Fej\'er means of the one-dimensional Walsh-Fourier series converge in $H_{p}(G)$ norm  (see also \cite{tep13}).

By applying inequality (\ref{vnk000f}) for the case $p=1/2$ the following necessary and sufficient conditions are found:

\textbf{Theorem \ref{theorem5.2.1}.} Let $f\in H_{1/2}(G)$ and $\{m_{k}:k\geq 0\}$ be increasing sequence of natural numbers, such that
\begin{equation*}
\omega _{H_{1/2}(G)}\left( \frac{1}{2^{\left\vert m_{k}\right\vert }},f\right)=o\left( \frac{1}{{V^2}\left( m_{k}\right) }\right) \text{ \ as \ }k\rightarrow\infty. \label{cond1f}
\end{equation*}%
Then
\begin{equation*}
\left\Vert \sigma_{m_{k}}f-f\right\Vert _{H_{1/2}(G)}\rightarrow 0\text{ \ as \ }k\rightarrow \infty .
\end{equation*}

Moreover, if $\{m_{k}:k\geq 0\}$ be increasing sequence of natural numbers, such that (\ref{1.dmk}) holds true, then there exists a martingale $f\in H_{1/2}(G)$ and subsequence $\{\alpha _{k}:k\geq 0\}\subset \{m_{k}:k\geq 0\}$ such that
\begin{equation*}
\omega _{H_{1/2}(G)}\left( \frac{1}{2^{\left\vert \alpha _{k}\right\vert }},f\right) =O\left( \frac{1}{{V^2}\left( \alpha _{k}\right) }\right) \text{ \ as\ }k\rightarrow \infty
\end{equation*}%
and
\begin{equation*}
\limsup\limits_{k\rightarrow \infty }\left\Vert \sigma_{\alpha_{k}}f-f\right\Vert _{1/2}>c>0\,\,\,\text{as\thinspace \thinspace \thinspace }%
k\rightarrow \infty ,  \label{cond10f}
\end{equation*}
wher  $c$ is an absolute constant.

By applying inequality (\ref{vnk0000f}) we also investigate case $0<p<1/2$. In  the third chapter of this thesis we prove that the following is true:

\textbf{Theorem \ref{theorem5.2.2}.} Let $0<p<1/2,$ $f\in H_{p}(G)$ and $\{m_{k}:k\geq 0\}$ be increasing sequence of natural numbers, such that
\begin{equation*}
\omega _{H_{p}(G)}\left( \frac{1}{2^{\left\vert m_{k}\right\vert }},f\right)=o\left( \frac{1}{2^{d\left( m_{k}\right) \left( 1/p-2\right) }}\right),
\text{ as \ }k\rightarrow \infty.
\end{equation*}%
Then
\begin{equation*} \label{con1f}
\left\Vert \sigma_{m_{k}}f-f\right\Vert _{H_{p}(G)}\rightarrow 0,\text{ \ as \ }k\rightarrow \infty .  
\end{equation*}

On the other hand, if $\{m_{k}:k\geq 0\}$ be increasing sequence of natural numbers satisfying the condition (\ref{1.dmk}), then there exists a martingale $f\in H_{p}(G)$ and subsequence $\{\alpha _{k}:k\geq 0\}\subset \{m_{k}:k\geq 0\},$ for which
\begin{equation*}
\omega _{H_{p}(G)}\left( \frac{1}{2^{\left\vert \alpha _{k}\right\vert }},f\right) =O\left( \frac{1}{2^{d\left( \alpha _{k}\right) \left(
1/p-2\right) }}\right), \text{ as \ }k\rightarrow \infty
\end{equation*}%
and
\begin{equation*}
\limsup\limits_{k\rightarrow \infty }\left\Vert \sigma_{\alpha_{k}}f-f\right\Vert _{weak-L_{p}(G)}>c_{p}>0,\text{\quad as\quad
}k\rightarrow \infty,  \label{con11f}
\end{equation*}%
where $c_{p}$ is constant depending only on $p$.

However, Simon in \cite{Si3} and \cite{si1} (see also \cite{cw, sw}) consider strong convergence theorems of the one-dimensional Walsh-Fouriere series and proved the following:

\textbf{Theorem Si1.}  Let $0<p\leq1$ and $f\in H_{1}\left( G\right)$. Then there exists an absolute constant $c_p$ depending only on $ p, $ such that the following inequality is true:
\begin{equation*}
\frac{1}{\log^{[p]} n}\overset{n}{\underset{k=1}{\sum }}\frac{\left\Vert
S_{k}f\right\Vert _{H_{p}(G)}}{k^{2-p}}\leq c_p\left\Vert f\right\Vert _{H_{p}(G)},
\end{equation*}

Analigical result for trigonometric system was proved in \cite{sm}, for unbounded Walsh systems in \cite{Ga1}.

In \cite{tep4} was proved that the following is true:

\textbf{Theorem T17.} for any $ 0<p<1 $ and non-decreasing function $\varphi :\mathbb{N}_{+}\rightarrow [1,$ $\infty )$ satisfying the condition
\begin{equation*}
\overline{\lim_{n\rightarrow \infty }}\frac{{n^{2-p}}}{\varphi \left( n\right) }=+\infty,
\end{equation*}
there exists a martingale $f\in H_{p}(G)$, such that
\begin{equation*}
\text{ }\underset{k=1}{\overset{\infty }{\sum }}\frac{\left\Vert S_{k}f\right\Vert _{weak-L_{p}(G)}^{p}}{\varphi (k)}=\infty ,\text{\qquad} (0<p<1).
\end{equation*}

Theorem Si1 follows that if $f\in H_{1}(G)$ then the following equalities are true:
\begin{equation*}
\underset{n\rightarrow \infty }{\lim }\frac{1}{\log n}\overset{n}{\underset{%
k=1}{\sum }}\frac{\left\Vert S_{k}f-f\right\Vert _{1}}{k}=0
\end{equation*}
and
\begin{equation*}
\underset{n\rightarrow \infty }{\lim }\frac{1}{\log n}\overset{n}{\underset{k=1}{\sum }}\frac{\left\Vert S_{k}f\right\Vert _{1}}{k}=\left\Vert
f\right\Vert _{H_{1}(G)}.
\end{equation*}

When $0<p<1$ and $f\in H_{p}(G),$ then Theorem Si1 follows that there exists an absolute constant $c_p$ depending only on $ p, $ such that 

\begin{equation*}
\frac{1}{n^{1/2-p/2}}\overset{n}{\underset{k=1}{\sum }}\frac{\left\Vert
S_{k}f\right\Vert _{H_{p}(G)}^{p}}{k^{3/2-p/2}}\leq c_{p}\left\Vert f\right\Vert
_{H_{p}(G)}^{p}.
\end{equation*}
Moreover,
\begin{equation*}
\frac{1}{n^{1/2-p/2}}\overset{n}{\underset{k=1}{\sum }}\frac{\left\Vert S_{k}f-f\right\Vert _{H_{p}(G)}^{p}}{k^{^{3/2-p}}}=0.
\end{equation*}
It follows the following equality
\begin{equation*}
\frac{1}{n^{1/2-p/2}}\overset{n}{\underset{k=1}{\sum }}\frac{\left\Vert
S_{k}f\right\Vert _{H_{p}(G)}^{p}}{k^{3/2-p/2}}=\left\Vert f\right\Vert
_{H_{p}(G)}^{p}.
\end{equation*}

In the third chapter of this thesis we consider strong convergence results of Fej\'er means of the one-dimensional Walsh-Fourier series. According to Theorem We1 and Theorem Gog2 we only have to consider case $0<p\leq1/2$ (for details see \cite{tep5}, see also \cite{BT1},  \cite{BTT1}, \cite{BT2}, \cite{BT3}, \cite{BPT1}):

\textbf{Theorem \ref{theorem5.3.1}.} Let $0<p\leq 1/2$ and $f\in H_{p}(G)$. Then there exists an absolute constant $c_p$ depending only on $p,$ such that

\begin{equation*}
\frac{1}{\log ^{\left[ 1/2+p\right] }n}\overset{n}{\underset{m=1}{\sum }}%
\frac{\left\Vert \sigma _{m}f\right\Vert _{H_{p}(G)}^{p}}{m^{2-2p}}\leq
c_{p}\left\Vert f\right\Vert _{H_{p}(G)}^{p}.
\end{equation*}%

Moreover, let $0<p<1/2$  and $\Phi :\mathbb{N}_{+}\rightarrow
\lbrack 1, $\textit{\ }$\infty )$ be non-decreasing, non-negative function, such that $\Phi \left( n\right) \uparrow \infty $ and
\begin{equation*}
\overline{\underset{k\rightarrow \infty }{\lim }}\frac{k^{
2-2p }}{\Phi \left({k}\right) }=\infty .
\end{equation*}
Then there exists a martingale $f\in H_{p}(G),$  such that
\begin{equation*}
\underset{m=1}{\overset{\infty }{\sum }}\frac{\left\Vert \sigma
_{m}f\right\Vert _{weak-L_{p}(G)}^{p}}{\Phi \left( m\right) }=\infty .
\end{equation*}

When $ p=1/2 $ is was also proved that the following is true:

\textbf{Theorem \ref{theorem5.3.2}.} Let $f\in H_{1/2}(G).$ Then
\begin{equation*}
\underset{n\in \mathbf{\mathbb{N}}_{+}}{\sup }\underset{\left\Vert f\right\Vert _{H_{p}(G)}\leq 1}{\sup } \frac{1}{n}\underset{m=1}{\overset{n}{\sum }}\left\Vert \sigma_{m}f\right\Vert_{1/2}^{1/2}=\infty .
\end{equation*}

Theorem \ref{theorem5.3.1} follows that if $f\in H_{1/2}(G)$ then the following equalities are true:
\begin{equation*}
\lim_{n\rightarrow \infty }\frac{1}{\log n}\overset{n}{\underset{k=1}{\sum }}\frac{\left\Vert \sigma _{k}f-f\right\Vert _{H_{1/2}(G)}^{1/2}}{k}=0
\end{equation*}
and
\begin{equation*}
\lim_{n\rightarrow \infty }\frac{1}{\log n}\overset{n}{\underset{k=1}{\sum }}%
\frac{\left\Vert \sigma _{k}f\right\Vert _{H_{1/2}(G)}^{1/2}}{k}=\left\Vert
f\right\Vert _{H_{1/2}(G)}^{1/2}.
\end{equation*}

When $0<p<1/2$ and $f\in H_{p}(G),$ then Theorem \ref{theorem5.3.1} follows that there exists an absolute constant $c_p$ depending only on $ p, $ such that
\begin{equation*}
\frac{1}{n^{1/2-p}}\overset{n}{\underset{k=1}{\sum }}\frac{\left\Vert \sigma
_{k}f\right\Vert _{H_{p}(G)}^{p}}{k^{3/2-p}}\leq c_{p}\left\Vert f\right\Vert
_{H_{p}(G)}^{p}.
\end{equation*}
Moreover,
\begin{equation*}
\frac{1}{n^{1/2-p}}\overset{n}{\underset{k=1}{\sum }}\frac{\left\Vert \sigma
_{k}f-f\right\Vert _{H_{p}(G)}^{p}}{k^{^{3/2-p}}}=0.
\end{equation*}%
It follows that 
\begin{equation*}
\frac{1}{n^{1/2-p}}\overset{n}{\underset{k=1}{\sum }}\frac{\left\Vert \sigma
_{k}f\right\Vert _{H_{p}(G)}^{p}}{k^{3/2-p}}=\left\Vert f\right\Vert
_{H_{p}(G)}^{p}.
\end{equation*}

For the two-dimensional case (for details see e.g. \cite{sws} and \cite{We1}) the following is true:

\textbf{Theorem S3.} Let $p>0$ and $ f\in H_{p}({G}^2) $. Then
\begin{equation} \label{1.Snn}
\left\Vert S_{2^{n},2^{n}}f-f\right\Vert _{p}\rightarrow 0, \text{ \ as \ }n\rightarrow \infty.
\end{equation}
Moreover,

\textbf{Theorem S4.} Let $p>0$  and $f\in H_{p}({G}^2)$. Then there exists an absolute constant $c_p$ depending only on $ p, $ such that
\begin{equation} \label{1.Snnmax}
\left\Vert \sup_{n\in \mathbb{N}}\left\vert S_{2^n,2^n} f\right\vert \right\Vert _{p}\leq c_{p}\left\Vert f\right\Vert _{H_{p}({G}^2)},
\end{equation}

By applying Theorem S4 we can conclude that the following holds true (for details see e.g. \cite{sws} and \cite{We1}):

\textbf{Theorem S5.} Let $p>0$ and $ f\in H_{p}({G}^2) $. Then there exists an absolute constant $ c_p$ depending only on $ p, $ such that
\begin{equation} \label{1.Snn0}
\left\Vert S_{2^{n},2^{n}}f\right\Vert _{H_{p}({G}^2)}\leq c_{p}\left\Vert f\right\Vert _{H_{p}({G}^2)}.
\end{equation}

On the other hand, (see \cite{tep7}) the following is true:

\textbf{Theorem T18.} Let $0<p\leq1.$ Then there exists a martingale $f\in H_{p}({G}^2)$  such that
\begin{equation*}
\underset{n\in \mathbb{N}}{\sup }\left\Vert S_{n,n}f\right\Vert
_{p}=\infty.
\end{equation*}

However, for the two-dimensional case Weisz \cite{We} proved the following:

\textbf{Theorem We4.} Let $\alpha \geq 0$ and $f\in H_{p}({G}^2).$ Then there exists an absolute constant $c_p$ depending only on $ p, $ such that
\begin{equation*}
\underset{n,m\geq 2}{\sup}\left(\frac{1}{\log n\log m}\right)^{\left[ p\right]}\underset{2^{-\alpha}\leq k/l\leq 2^{\alpha },\text{ }\left(
k,l\right) \leq \left(n,m\right)}{\sum}\frac{\left\Vert S_{k,l}f\right\Vert _{p}^{p}}{\left(kl\right)^{2-p}}\leq c_{p}\left\Vert f\right\Vert _{H_{p}({G}^2)}^{p},
\end{equation*}%
where $0<p\leq1$ and $\left[ p\right] $ denotes integer part of real number $p$.

Moreover, sharpness of of the rates of weights are proved  \cite{tep22} was proved that rate of the  ${\left(kl\right)^{2-p}}$ is sharp. 

Goginava and Gogoladze in \cite{gg} generalized this result in the case when $\alpha=0$:

\textbf{Theorem GG1.} Let $f\in H_{1}({G}^2) $. Then there exists an absolute constant $c$, such that
\begin{equation*}
\sum\limits_{n=1}^{\infty }\frac{\left\Vert S_{n,n}f\right\Vert _{1}}{n\log^{2}n}\leq c\left\Vert f\right\Vert _{H_{1}{({G}^2)}}.
\end{equation*}

In \cite{tep14} was proved that rate of the weights $\left\{ n\log^{2}n \right\} _{n=1}^{\infty }$ is sharp. The following is true:

\textbf{Theorem T19.} Let
 \ $\Phi :\mathbb{N}\rightarrow \lbrack 1,$ $\infty )$ be non-decreasing function satisfying the condition $\lim_{n\rightarrow \infty}\Phi \left( n\right) =+\infty.$ Then
\begin{equation*}
\underset{\left\Vert f\right\Vert _{H_{1}({G}^2)}\leq 1}{\sup }\underset{n=1}{\overset{\infty }{\sum }}\frac{\left\Vert S_{n,n}f\right\Vert _{1}\Phi
\left( n\right) }{n\log ^{2}\left( n+1\right) }=\infty .
\end{equation*}

Theorem GG1 follows that if $f\in H_{1}({G}^2),$ then
\begin{equation*}
\frac{1}{\log ^{1/2} n}\overset{n}{\underset{k=1}{\sum }}\frac{\left\Vert S_{k,k}f\right\Vert _{H_{1}({G}^2)}}{k\log ^{3/2} k}\leq c\left\Vert f\right\Vert
_{H_{1}({G}^2)},\text{ }
\end{equation*}
Moreover,
\begin{equation*}
\lim_{n\rightarrow \infty }\frac{1}{\log^{1/2} n}\overset{n}{\underset{k=1}{\sum }}%
\frac{\left\Vert S_{k,k}f-f\right\Vert _{H_{1}({G}^2)}^{2/3}}{k\log^{3/2} k}=0
\end{equation*}%
It follows the following equality
\begin{equation*}
\lim_{n\rightarrow \infty }\frac{1}{\log^{1/2} n}\overset{n}{\underset{k=1}{\sum }}%
\frac{\left\Vert S_{k,k}f\right\Vert _{H_{2/3}({G}^2)}^{2/3}}{k{\log^{3/2} k}}=\left\Vert
f\right\Vert _{H_{2/3}({G}^2)}^{2/3}.
\end{equation*}

In the fourth chapter (see also \cite{tep15}) of this thesis we consider strong convergence of Marcinkiewicz means with respect to the two-dimensional partial sums of Walsh-Fourier series when $ 0<p<1 $:

\textbf{Theorem  \ref{theorem6.1.1}} Let $ 0<p<1 $ and $f\in H_{p}({G}^2).$ Then there exists an absolute constant $ c_p$ depending only on $ p, $ such that
\begin{equation}  \label{1.th}
\sum\limits_{n=1}^{\infty }\frac{\left\Vert S_{n,n}f\right\Vert _{p}^{p}}{n^{3-2p}}\leq c_{p}\left\Vert f\right\Vert _{H_{p}({G}^2)}^{p}.
\end{equation}%

Moreover, if $0<p<1$ and $\Phi :\mathbb{N}\rightarrow [1,$ $\infty )$ be non-decreasing function satifying condition  $\underset{n\rightarrow
\infty }{\lim }\Phi \left( n\right) =+\infty$, then there exists martingale $f\in H_{p}({G}^2) $ such that

\begin{equation*}
\underset{n=1}{\overset{\infty}{\sum}} \frac{\left\|S_{n,n}f\right\|_{weak-L_{p}(G^2)}^{p}\Phi\left(n\right)} {n^{3-2p}}=\infty.
\end{equation*}

Theorem  \ref{theorem6.1.1} follows that, if $0<p<1$ and $f\in H_{p}({G}^2),$ then there exists an absolute constant $ c_p$ depending only on $ p, $ such that

\begin{equation*}
\frac{1}{n^{1/2-p}}\overset{n}{\underset{k=1}{\sum }}\frac{\left\Vert S_{k,k}f\right\Vert _{H_{p}({G}^2)}^{p}}{k^{3/2-p}}\leq c_{p}\left\Vert f\right\Vert
_{H_{p}({G}^2)}^{p},
\end{equation*}
Moreover,
\begin{equation*}
\frac{1}{n^{1/2-p}}\overset{n}{\underset{k=1}{\sum }}\frac{\left\Vert S_{k,k}f-f\right\Vert _{H_{p}({G}^2)}^{p}}{k^{3/2-p}}=0.
\end{equation*}
It follows the following equality
\begin{equation*}
\frac{1}{n^{1/2-p}}\overset{n}{\underset{k=1}{\sum }}\frac{\left\Vert S_{k,k}f\right\Vert _{H_{p}(G)}^{p}}{k^{3/2-p}}=\left\Vert f\right\Vert
_{H_{p}(G)}^{p}.
\end{equation*}

Weisz (for details see e.g. \cite{We1}) consider Marcinkiewicz means with respect to the two-dimensional partial sums of Walsh-Fourier series and proved the following:

\textbf{Theorem We5.} Let $p>2/3$ and $f\in H_{p}({G}^2)$. Then there exists an absolute constant $ c_p$ depending only on $ p, $ such that
\begin{equation*}
\left\Vert \mathcal{M}_{n}f-f\right\Vert _{H_{p}({G}^2)}\rightarrow 0, \text{ \  \ }n\rightarrow \infty.
\end{equation*}

Goginava \cite{gog3} proved that the following is true:

\textbf{Theorem Gog2.} Let $0<p\leq2/3$. Then there exists a martingale $f\in H_{p}({G}^2)$, such that
\begin{equation*}
\underset{n\in \mathbb{N}}{\sup }\left\Vert \mathcal{M}_{n}f\right\Vert_{H_{p}({G}^2)}=\infty.
\end{equation*}

Goginava in \cite{Go} consider subsequence $\mathcal{M}_{2^{n}}$  of Marcinkiewicz means with respect to the two-dimensional partial sums of Walsh-Fourier series and proved that the following is true:

\textbf{Theorem Gog3.} Let $p>1/2$ and $f\in H_{p}({G}^2)$. Then
\begin{equation} \label{1.Mnn}
\left\Vert \mathcal{M}_{2^{k}}f-f\right\Vert _{H_{p}({G}^2)}\rightarrow 0,  \text{ \ as \ }k\rightarrow \infty.
\end{equation}

Moreover, there exists a martingale $f\in H_{p}({G}^2),$ $\text{ } (0<p\leq1/2)$ such that
\begin{equation*}
\underset{n\in \mathbb{N}}{\sup }\left\Vert \mathcal{M}_{2^{n}}f\right\Vert_{H_{p}({G}^2)}=\infty .
\end{equation*}

In \cite{tep17} was investigated strong convergence theorems of Marcinkiewicz means with respect to the two-dimensional partial sums of Walsh-Fourier series when $0<p<2/3$:

\textbf{Theorem NT1.} Let $0<p<2/3$ and $f\in H_{p}({G}^2)$. Then there exists an absolute constant $ c_p$ depending only on $ p, $ such that
\begin{equation*}
\sum_{m=1}^\infty\frac{\left\Vert \mathcal{M}_{m}f\right\Vert _{H_{p}({G}^2)}^{p}}{m^{3-3p}}\leq c_{p}\left\Vert f\right\Vert_{H_{p}({G}^2)}^{p}.
\end{equation*}

Moreover, if $0<p<2/3$ and $\Phi :\mathbb{N}_{+}\rightarrow \lbrack 1,$ $\infty )$ be non-decreasing function satisfying the condition $\Phi \left(n\right) \uparrow \infty $ and
\begin{equation*}
\overline{\lim_{k\rightarrow\infty}}\frac{k^{3-3p}}{\Phi\left(k\right)}=\infty,
\end{equation*}
then there exists a martingale $f\in H_{p}\left( G^{2}\right) ,$ such that
\begin{equation*}
{\sum_{m=1}^\infty }\frac{\left\Vert \mathcal{M}_{m}f\right\Vert _{weak-L_{p}({G}^2)}^{p}}{\Phi \left( m\right) }=\infty.
\end{equation*}

Theorem NT1 follows that if $0<p<2/3$ and $f\in H_{p}({G}^2),$ then there exists an absolute constant $ c_p$ depending only on $ p,$ such that

\begin{equation*}
\frac{1}{n^{1-3p/2}}\overset{n}{\underset{k=1}{\sum }}\frac{\left\Vert \mathcal{M}
_{k}f\right\Vert _{H_{p}({G}^2)}^{p}}{k^{2-3p/2}}\leq c_{p}\left\Vert f\right\Vert
_{H_{p}({G}^2)}^{p},
\end{equation*}
Moreover,
\begin{equation*}
\frac{1}{n^{1-3p/2}}\overset{n}{\underset{k=1}{\sum }}\frac{\left\Vert \mathcal{M}
_{k}f-f\right\Vert _{H_{p}({G}^2)}^{p}}{k^{2-3p/2}}=0,
\end{equation*}
It follows the following equality:
\begin{equation*}
\frac{1}{n^{1-3p/2}}\overset{n}{\underset{k=1}{\sum }}\frac{\left\Vert \mathcal{M}
_{k}f\right\Vert _{H_{p}(G)}^{p}}{k^{2-3p/2}}=\left\Vert f\right\Vert
_{H_{p}(G)}^{p}.
\end{equation*}

In the fourth chapter (see \cite{tep16}) we consider strong convergence results of Marcinkiewicz means with respect to the two-dimensional partial sums of Walsh-Fourier series, when $p=2/3$:

\textbf{Theorem \ref{th7.1.1}} Let $f\in H_{2/3}({G}^2).$ Then there exists an absolute constant $ c, $ such that
\begin{equation*}  \label{ineq2}
\frac{1}{\log n}\sum_{m=1}^{n}\frac{\left\Vert\mathcal{M}_{m}f \right\Vert_{H_{2/3}({G}^2)}^{2/3}}{m} \leq c \left\Vert f \right\Vert_{H_{2/3}({G}^2)}^{2/3}.
\end{equation*}

From these results we obtain that, if $f\in H_{2/3}({G}^2)$ then
\begin{equation*}
\lim_{n\rightarrow \infty }\frac{1}{\log n}\overset{n}{\underset{k=1}{\sum }}%
\frac{\left\Vert \mathcal{M} _{k}f-f\right\Vert _{H_{2/3}({G}^2)}^{2/3}}{k}=0
\end{equation*}%
and
\begin{equation*}
\lim_{n\rightarrow \infty }\frac{1}{\log n}\overset{n}{\underset{k=1}{\sum }}%
\frac{\left\Vert \mathcal{M} _{k}f\right\Vert _{H_{2/3}({G}^2)}^{2/3}}{k}=\left\Vert
f\right\Vert _{H_{2/3}({G}^2)}^{2/3}.
\end{equation*}

For the two-dimensional case Weisz \cite{We6} proved that the following is true:

\textbf{Theorem We6.} Let $p>2/3$ and $f\in H_{p}\left(G^{2}\right)$. Then the maximal operator of Marcinkiewicz means with respect to the two-dimensional partial sums of Walsh-Fourier series $\mathcal{M}^{\ast}f $ is bounded from $H_{p}\left(G^{2}\right)$ to $L_{p}\left(G^{2}\right)$:

\begin{equation*}
\left\|\mathcal{M}^{\ast}f\right\|_{p}\leq c_p\left\|f\right\|_{H_p(G^2)},
\end{equation*}
where $c_p$ is an absolute constant, depending only on $p$.

Goginava \cite{gog1} also proved that the following is true:

\textbf{Theorem Gog4.} Let $f\in H_{2/3}\left( G^{2}\right) $. Then there exists an absolute constant $ c, $ such that

\begin{equation*}
\left\| \mathcal{M}^{\ast}f \right\|_{weak-L_{2/3}(G^2)} \leq c\left\|
f\right\| _{H_{2/3}}.
\end{equation*}

Moreover, there exists a martingale $f\in H_{p}({G}^2),$ $\text{ } (0<p\leq2/3)$ such that
\begin{equation*}
\underset{n\in \mathbb{N}}{\sup }\left\Vert \mathcal{M}_{n}f\right\Vert_{p}=\infty .
\end{equation*}

Goginava \cite{gog3} also consider restricted maximal operator of Marcinkiewicz means with respect to the two-dimensional partial sums of Walsh-Fourier series $\sup_{n\in\mathbb{N}}\left\vert \mathcal{M}_{2^{n}}\right\vert$ and show that the following is true:

\textbf{Theorem Gog5.} Let $p>1/2$ and $f\in H_{p}({G}^2).$ Then there exists an absolute constant $ c_p$ depending only on $ p, $ such that
\begin{equation} \label{1.Mnnmax}
\left\Vert\sup_{n\in\mathbb{N}}\left\vert \mathcal{M}_{2^{n}}f\right\vert\right\Vert_{p}\leq c_{p} \left\Vert f\right\Vert _{H_{p}({G}^2)}.
\end{equation}

Moreover, there exists a martingale $f\in H_{p}({G}^2),$ $\text{ } (0<p\leq1/2)$ such that
\begin{equation*}
\underset{n\in \mathbb{N}}{\sup }\left\Vert \mathcal{M}_{2^{n}}f\right\Vert_{p}=\infty .
\end{equation*}

In \cite{nagy} and \cite{tep17} we investigate boundedness of weighted maximal operators when $0<p\leq2/3$:

\textbf{Theorem NT2.} Let $0<p\leq2/3.$ Then the maximal operator 
\begin{equation*}
\widetilde{\mathcal{M}}^{\ast }:=\sup_{n\in \mathbb{N}}\frac{\left\vert
\mathcal{M}_{n}\right\vert }{{(n+1)}^{2/p-3}\log ^{3[1/3+p]/2}\left( n+1\right)}
\end{equation*}%
is bounded from $H_{p}(G^{2})$ to $L_{p}(G^{2})$.

Moreover, if $\varphi :\mathbb{N}\rightarrow \lbrack 1,\infty )$ be non-decreasing function satisfying the condition
$$
{\lim_{n\rightarrow \infty }}\frac{ n^{2/p-3}\log ^{3[1/3+p]/2}\left( n\right)}
{\varphi \left( n\right) }=+\infty,
$$
then
\begin{equation*}
\sup_{n\in \mathbb{N}}\left\Vert \frac{{\mathcal M}_{n}f}{\varphi \left( n\right)}\right\Vert_{p}=\infty .
\end{equation*}

From Theorem NT2 we get that for $ 0<p\leq1/2 $ and $ f\in H_{p}(G^2)$, there exists a absolute constant $c_{p},$ depending only on $p$, such that:
\begin{equation}  \label{1.c22}
\left\Vert \mathcal{M}_{n}f\right\Vert _{p}\leq c_{p} {{(n+1)}^{2/p-3}\log ^{3[1/3+p]/2}\left( n+1\right)} \left\Vert f\right\Vert _{H_{p}({G}^2)}.
\end{equation}

By applying inequality (\ref{1.c22}) in \cite{tep17} (see also \cite{NT4}) we obtain necessary and sufficient conditions for modulus of continuity of martingale $f\in H_{p}(G^2)$, for which Marcinkiewicz means with respect to the two-dimensional partial sums of Walsh-Fourier series of $f \in H_{p}(G^2)$ converge in $H_{p}(G^2)$ norm.

\textbf{Theorem NT3.} Let $1/2<p<2/3,$ $f\in H_{p}\left( G^{2}\right) \ $ and
\begin{equation*}
\omega_{H_{p}({G}^2)} \left( \frac{1}{2^{k}},f\right)=o\left( \frac{1}{2^{k\left(
2/p-3\right) }}\right) ,\text{ as \ }k\rightarrow \infty .  \label{10A}
\end{equation*}
Then
\begin{equation*}
\left\Vert \mathcal{M}_{n}f-f\right\Vert _{H_{p}({G}^2)}\rightarrow
0,\text{ as }n\rightarrow \infty .
\end{equation*}

Moreover, if $0<p<2/3,$ then there exists a martingale $f\in H_{p}(G^{2}),$ such that
\begin{equation*}
\omega_{H_{p}({G}^2)}\left( \frac{1}{2^{k}},f\right)=O\left( \frac{1}{2^{k\left(
2/p-3\right) }}\right) ,\text{ as \ }k\rightarrow \infty
\end{equation*}
and
\begin{equation*}
\left\Vert \mathcal{M}_{n}f-f\right\Vert _{weak-L_p({G}^2)}\nrightarrow
0,\text{ \ as \  }n\rightarrow \infty.
\end{equation*}

If we apply again (\ref{1.c22}) and improve method which was investigated in \cite{tep17} in the fourth chapter (see also \cite{tep18}) we obtain that the following is true:

\textbf{Theorem \ref{theorem7.2.1}.} Let $f\in H_{2/3}(G^{2})$  and
\begin{equation*}
\omega_{H_{2/3}({G}^2)} \left( \frac{1}{2^{k}},f\right)=o\left( \frac{1}{{k}^{3/2}}\right) ,\text{ as \ }k\rightarrow \infty.
\end{equation*}
Then
\begin{equation*}
\left\Vert \mathcal{M}_{n}f-f\right\Vert _{H_{2/3}({G}^2)}\rightarrow
0,\text{ as }n\rightarrow \infty .
\end{equation*}

On the other hand, there exists a martingale  $f\in H_{2/3}(G^{2}),$ such that
\begin{equation*}
\omega_{H_{2/3}({G}^2)}\left(\frac{1}{2^{k}},f\right)=O\left( \frac{1}{{k}^{3/2}}\right),\text{ as \ }k\rightarrow \infty
\end{equation*}
and
\begin{equation*}
\left\Vert \mathcal{M}_{n}f-f\right\Vert _{2/3}\nrightarrow
0,\text{ \ as \ }n\rightarrow \infty.
\end{equation*}

\newpage

\section{Partial sums with respect to the one-dimensional Walsh-Fourier series on the martingale Hardy spaces}

\subsection{Basic notations}

\ \ Denote by $\mathbb{N}_{+}$ the set of the positive integers and by  $\mathbb{N}:=\mathbb{N}_{+}\cup \{0\}$ the set of non-negative integers. Denote by $Z_{2}$ the additive group of integers modulo-$2,$ which contains only two elements  $Z_{2}:=\{0,1\},$ group operation is modulo-$2$ sum and all sets are open. 

Define the group $G$ as the complete direct product of the groups $Z_{2}$ with the product of the discrete topologies  $Z_{2}$. The direct product $\mu $ of the measures $\mu _{n}\left( \{j\}\right)
:=1/2,\ (j\in Z_{2})$ is the Haar measure on $G$ with $\mu
\left( G\right) =1.$

The elements of $G$ are represented by sequences
\begin{equation*}
x:=(x_{0},x_{1},...,x_{j},...)\qquad \left( x_{k}=0,1\right) .
\end{equation*}

It is easy to give a base for the neighbourhood of $G$
\begin{equation*}
I_{0}\left( x\right) :=G,
\end{equation*}%
\begin{equation*}
I_{n}(x):=\{y\in G\mid y_{0}=x_{0},...,y_{n-1}=x_{n-1}\}\text{ }(x\in G,%
\text{ }n\in \mathbb{N}).
\end{equation*}

Set  $I_{n}:=I_{n}\left( 0\right) $ for any $n\in \mathbb{N}$ and $\overline{I_{n}}:=G$ $\backslash $ $I_{n}$.

It is evident that
\begin{equation} \label{2.1.2}
\overline{I_{M}}=\left( \overset{M-2}{\underset{k=0}{\bigcup }}\overset{M-1}{%
\underset{l=k+1}{\bigcup }}I_{l+1}\left( e_{k}+e_{l}\right) \right) \bigcup \left( \underset{k=0}{\bigcup\limits^{M-1}}I_{M}\left( e_{k}\right) \right)= \underset{k=0}{\bigcup\limits^{M-1}}I_{k} \backslash I_{k+1}.
\end{equation}

If $n\in \mathbb{N},$ then it can be uniquely expressed as

\begin{equation*}
n=\sum_{k=0}^{\infty }n_{j}2^{j}
\end{equation*}
where  $n_{j}\in Z_{2}$  $~(j\in \mathbb{N})$ and only a finite number of $n_{j}$s differ from zero.

Set
\begin{equation*}
\left\langle n\right\rangle :=\min \{j\in \mathbb{N},n_{j}\neq 0\}\text{ \ \
and \ \ \ }\left\vert n\right\vert :=\max \{j\in \mathbb{N},n_{j}\neq 0\},
\end{equation*}%
It is evident that $2^{\left\vert n\right\vert }\leq n\leq 2^{\left\vert n\right\vert
+1}.$

Let
\begin{equation*}
d\left( n\right) :=\left\vert n\right\vert -\left\langle n\right\rangle,
\text{ \ for any \ }n\in \mathbb{N}.
\end{equation*}

Denote by $V(n)$ variation of natural number $n\in \mathbb{N}$ 

\begin{equation*}
V\left( n\right) =n_{0}+\overset{\infty }{\underset{k=1}{\sum }}\left|
n_{k}-n_{k-1}\right| .
\end{equation*}

Define $k$-th Rademacher functions by
\begin{equation*}
r_{k}\left( x\right) :=\left( -1\right) ^{x_{k}}\text{\qquad }\left( \text{ }%
x\in G,\text{ }k\in \mathbb{N}\right) .
\end{equation*}

By using Rademacher functions we define Walsh system $w:=(w_{n}:n\in \mathbb{N})$ $G$ as:
\begin{equation*}
w_{n}(x):=\underset{k=0}{\overset{\infty }{\prod }}r_{k}^{n_{k}}\left(x\right) =r_{\left\vert n\right\vert }\left( x\right) \left( -1\right) ^{%
\underset{k=0}{\overset{\left\vert n\right\vert -1}{\sum }}n_{k}x_{k}}\text{%
\qquad }\left( n\in \mathbb{N}\right).
\end{equation*}

The norm (quasi-norm) of space $L_{p}(G), \ \left( 0<p<\infty \right)$ is defined as

\begin{equation*}
\left\Vert f\right\Vert _{p}^p:=\left(\int_{G}\left\vert f(x)\right\vert
^{p}d\mu (x)\right).
\end{equation*}
and norm (quasi-norm) of space $weak-L_{p}(G) $ is defined by

\begin{equation*}
\left\Vert f\right\Vert _{weak-L_{p}(G)}^p:=\underset{\lambda >0}{\sup}
\lambda^p \mu \left (x\in{G}:  |f|>\lambda \right)<+\infty.
\end{equation*}

Walsh system is orthonormal and complete in $L_{2}\left(G \right) $ (see  \cite{sws}).

For any $f\in L_{1}\left(G\right) $ the numbers
\begin{eqnarray*}
\widehat{f}\left(n\right) :=\int\limits_{G}f(x)w_{n}(x)d\mu({x})
\end{eqnarray*}
are called $n$-th Walsh-Fourier coefficient of $f$.

$n$-th partial sum is denoted by
\begin{equation*}
S_{n}(f;{x}):=\sum\limits_{i=0}^{n-1}
\widehat{f}\left(i\right) w_{i}\left( x\right) .
\end{equation*}

Dirichlet kernels are defined by
\begin{equation*}
D_{n}\left( x\right) :=\sum\limits_{i=0}^{n-1} w_{i}\left( x\right) .
\end{equation*}

We also define the following maximal operators

\begin{eqnarray*}
{S}^{\ast}f&=&\sup_{n\in \mathbb{N}}\left\vert S_{{n}}f\right\vert  \\
\widetilde{S}_{\#}^{\ast}f&=&\sup_{n\in \mathbb{N}}\left\vert S_{2^{n}}f\right\vert
\end{eqnarray*}

The $\sigma $-algebra generated by the intervals $ I_{n}(x) $
with measure $2^{-n}$ is denoted by $\digamma_{n}\left(n\in \mathbb{N}\right).$ Conditional exponential operator with respect to
$\digamma _{n}\left( n\in \mathbb{N}\right) $ is denoted by $E_{n}$ and it is given by 
\begin{eqnarray*}
 E_{n}f(x)
&=& S_{2^{n}}f\left( x\right) \\
&=&\sum_{k=0}^{2^{n}-1}
\widehat{f}\left(k\right) w_{k}(x) \\
&=& \frac{1}{\left| I_{n}\left(x\right) \right| }\int_{I_{n}\left(x\right)
}f(x)d\mu (x),
\end{eqnarray*}
where $\left| I_{n}\left(x\right) \right| =2^{-n}$  denotes length of set
 $I_{n}\left(x\right)$.

Sequence $f=\left(f_{n},\text{ }n\in \mathbb{N}\right) $ of functions $f_{n}\in L_{1}\left( G\right) $ is called dyadic martingale  (for details see \cite{Nev}, \cite{sws}) if

$\left( i\right) $ $f_{n}$ is measurable with respect to  $\sigma-$ algebras $\digamma _{n}$ for any $n\in
\mathbb{N}$,

$\left( ii\right) $ $E_{n}f_{m}=f_{n}$ for any $n\leq m$.

The maximal function of a martingale $f$ is defined by
\begin{equation*}
f^{\ast }=\sup_{n\in \mathbb{N}}\left\vert f_{n}\right\vert.
\end{equation*}

In the case $f\in L_{1}\left(G\right) ,$ the maximal functions are also be given by:
\begin{equation*}
f^{\ast }\left(x\right) =\sup\limits_{n\in \mathbb{N}}\frac{1}{\mu \left(
I_{n}\left( x\right) \right) }\left\vert \int\limits_{I_{n}\left( x\right)
}f\left( u\right) d\mu \left( u\right) \right\vert.
\end{equation*}

For $0<p<\infty $ the Hardy martingale space $H_{p}\left( G\right)$
consists of all martingales, for which
\begin{equation*}
\left\Vert f\right\Vert _{H_{p}\left( G\right)}:=\left\Vert f^{\ast }\right\Vert
_{p}<\infty .
\end{equation*}

A bounded measurable function $a$ is said to be a $p$-atom if there exists an
dyadic interval $I$, such that%
\begin{equation*}
\int_{I}ad\mu =0,\text{ \ \ }\left\Vert a\right\Vert _{\infty }\leq \mu
\left( I\right) ^{-1/p},\text{ \ \ supp}\left( a\right) \subset I.\qquad
\end{equation*}%

It is easy to show that for martingale $f=(f_{n},n\in \mathbb{N}) $ and for any ${k}\in {\mathbb{N}}$ there exists a limit
\begin{equation*}
\widehat{f}\left(k\right):=\lim_{{n}\rightarrow \infty }\int_{{G}}f_{{{n}}}\left(
{x}\right) w_{{k}}\left( {x}\right) d\mu \left({x}\right)
\end{equation*}
and it is called  $k$-th Walsh-Fourier coefficients of $f$.

If $f_{0}\in L_{1}\left( G\right)$ and $f:=(E_{n}f_{0}:n\in \mathbb{N})$ is regular martingale then
\begin{eqnarray*}
\widehat{f}\left(k\right)&=&\int_{{G}}f\left( {x}\right) w_{k}(x)
d\mu \left({x}\right)=\widehat{f_{0}}\left(k\right),\ \ k\in \mathbb{N}.
\end{eqnarray*}

The modulus of continuity in  $H_{p}(G)$ space is defined by
\begin{equation*}
\omega _{H_{p}\left( G\right)}\left( \frac{1}{2^{n}},f\right) :=\left\Vert
f-S_{2^{n}}f\right\Vert _{H_{p}\left( G\right)}.
\end{equation*}

It is important to describe how can be understood difference $f-S_{2^{n}}f$, where $f$ be martingale $S_{2^{n}}f$ is a function:

\begin{remark} \label{lemma2.3.6666}
Let $0<p\leq 1.$ Since
\begin{equation*}
S_{2^{n}}f=f^{\left( n\right) }\in L_1(G),\text{ \ where \ }f=\left( f^{\left( n\right)}:n\in \mathbb{N}\right) \in H_{p}(G)
\end{equation*}%
and
\begin{eqnarray*}
\left( S_{2^{k}}f^{\left( n\right) }:k \in \mathbb{N}\right)
&=&\left( S_{2^{k}}S_{2^{n}},k \in \mathbb{N}\right) \\
&=&\left( S_{2^{0}}f,\ldots ,S_{2^{n-1}}f,S_{2^{n}}f,S_{2^{n}}f,\ldots \right) \\
&=&\left( f^{\left( 0\right) },\ldots ,f^{\left( n-1\right) },f^{\left(
n\right) },f^{\left( n\right) },\ldots \right).
\end{eqnarray*}
Under the difference $  f-S_{2^{n}}f $ we mean the following martingale:
$$ f:=(\left( f-S_{2^{n}}f\right) ^{\left( k\right) },  \ \ k\in \mathbb{N} )$$
where
\begin{equation*}
\left( f-S_{2^{n}}f\right) ^{\left( k\right) }=\left\{
\begin{array}{ll}
0, & k=0,\ldots ,n, \\
f^{\left( k\right) }-f^{\left( n\right) }, & k\geq n+1,\end{array}
\right.
\end{equation*}
\end{remark}

Consequently, the norm $ \left\Vert f-S_{2^{n}}f\right\Vert _{H_{p}\left( G\right)}  $ is understood as  $ H_p $-norm of 
$$  f-S_{2^{n}}f=(\left( f-S_{2^{n}}f\right) ^{\left( k\right) }, k\in \mathbb{N}) $$.

Watari \cite{wat} showed that there are strong connections between
\begin{equation*}
\omega _{p}\left( \frac{1}{2^{n}},f\right) ,\text{ \ }E_{2^{n}}\left(
L_{p},f\right) \text{ \ \ and\ \ \ }\left\Vert f-S_{2^{n}}f\right\Vert
_{p},\ \ p\geq 1,\text{ \ }n\in \mathbb{N}.
\end{equation*}

In particular,%
\begin{equation*}
\frac{1}{2}\omega _{p}\left( \frac{1}{2^{n}},f\right) \leq \left\Vert
f-S_{2^{n}}f\right\Vert _{p}\leq \omega _{p}\left( \frac{1}{2^{n}},f\right)
\label{eqvi}
\end{equation*}%
and%
\begin{equation*}
\frac{1}{2}\left\Vert f-S_{2^{n}}f\right\Vert _{p}\leq E_{2^{n}}\left(
L_{p},f\right) \leq \left\Vert f-S_{2^{n}}f\right\Vert _{p}.
\end{equation*}

\subsection{Auxiliary lemmas}

 \ First we present and prove equalities and estimations of Dirichlet kernel and Lebesgue constants with respect to the one-dimensional Walsh-Fourier systems (see Lemmas \ref{lemma0}-\ref{lemma7}).

First equality of the following Lemma is proved in \cite{sws} and second identity is proved in  \cite{Ga2}:

\begin{lemma} \label{lemma0}
Let $j, n\in \mathbb{N}$. Then
\begin{equation*}
D_{j+2^n}=D_{2^n}+w_{2^n}D_j,\text{ \ when \ } j\leq 2^n,
\end{equation*}
and
\begin{equation*}
D_{2^n-j}=D_{2^n}-\psi_{2^n-1}{D_j},\text{ \ when \ }j<2^n.
\end{equation*}
\end{lemma}

The following estimation of Dirichlet kernel with respect to the one-dimensional Walsh-Fourier systems is proved in \cite{sws}:

\begin{lemma} \label{lemma1}

Let $n\in \mathbb{N}$. Then

\begin{equation*}
D_{2^{n}}\left( x\right) =\left\{
\begin{array}{l}
2^{n},\text{ \ if  \  }x\in I_{n},
\\ 0,\text{ \ if \ }x\notin I_{n},
\end{array}
\right.
\end{equation*}
and
\begin{equation*}
D_{n}=w_{n}\overset{\infty }{\underset{k=0}{\sum }}n_{k}r_{k}D_{2^{k}}=w_{n}%
\overset{\infty }{\underset{k=0}{\sum }}n_{k}\left(
D_{2^{k+1}}-D_{2^{k}}\right) ,\text{ \ for \ }n=\overset{\infty }{\underset{i=0
}{\sum }}n_{i}2^{i}.
\end{equation*}%
\end{lemma}

The following two-sided estimations of Lebesgue constants with respect to the one-dimensional Walsh-Fourier systems is proved in \cite{sws} and second equality is proved in \cite{fi}:

\begin{lemma} \label{lemma2}

Let $n\in \mathbb{N}.$ Then

\begin{equation*}
\frac{1}{8}V\left( n\right) \leq \left\Vert D_{n}\right\Vert _{1}\leq
V\left( n\right)
\end{equation*}
and
\begin{equation*}
\frac{1}{n\log n}\underset{k=1}{\overset{n}{\sum }}V\left( k\right) =\frac{1%
}{4\log 2}+o\left( 1\right).
\end{equation*}

\end{lemma}

Hardy martingale space $H_{p}\left( G\right) $ for any $0<p\leq 1$ can be characterize by simple functions which are called $ p $-atoms. The following is true (for details see \cite{S}, \cite{We1} and \cite{We5}):

\begin{lemma} \label{lemma3.2.4}
A martingale $f=\left(f_{n},\text{ }n\in \mathbb{N}%
\right) $ belongs to $H_{p}(G)\left( 0<p\leq 1\right) $ if and only if there exists a sequence of $ p $-atoms of
 $\left( a_{k},k\in \mathbb{N}\right) $ and sequence of real numbers 
 $\left( \mu _{k},\text{ }k\in \mathbb{N}\right) $ such that for all $n\in \mathbb{N}$,

\begin{equation} \label{3.2.1.lemma3.2.4}
\qquad \sum_{k=0}^{\infty }\mu _{k}S_{2^{n}}a_{k}=f_{n}
\end{equation}%
and
\begin{equation*}
\qquad \sum_{k=0}^{\infty }\left\vert \mu _{k}\right\vert ^{p}<\infty.
\end{equation*}%
Moreover,
\begin{equation*}
\left\Vert f\right\Vert _{H_{p}(G)}\backsim \inf \left(
\sum_{k=0}^{\infty }\left\vert \mu _{k}\right\vert ^{p}\right) ^{1/p},
\end{equation*}
where the infimum is taken over all decomposition of  $f$ of the form (\ref{3.2.1.lemma3.2.4}).
\end{lemma}

\text{ \qquad }  The next five Examples of martingales will be used many times to prove sharpness of our main results. Such counterexamples first appear in the papers  of Goginava  \cite{Go} (see also \cite{GoPubl}). Such constructions of martingales are also used in the papers \cite{BPTW}, \cite{BT1}, \cite{BNPT}, \cite{MPT}, \cite{MST}, \cite{tep18}, \cite{PTW1}, \cite{PTW2}, \cite{tep1}, \cite{tep5}, \cite{tep12}, \cite{tep13}, \cite{tep15},  \cite{tep20}, \cite{tep19}, \cite{tep22}, \cite{tep23}, \cite{tep21}, \cite{tep23}.
For the one-dimensional case we use martingales which were used in \cite{tep_thesis}. So, we leave out the details of proof.

\begin{example}  \label{example2.2.1}
Let $0<p\leq 1,$  $\left\{ \lambda _{k}:k\in \mathbb{N}\right\}$
be sequence of real numbers
\begin{equation} \label{3.3.2aa}
\sum_{k=0}^{\infty }\left\vert \lambda _{k}\right\vert ^{p}\leq c_{p}<\infty
\end{equation}%
and $\left\{ a_{k}:k\in \mathbb{N}\right\} $ be sequence of $p$-atoms, given by 
\begin{equation*}
a_{k}(x):={2^{\left\vert \alpha _{k}\right\vert(1/p-1)}}\left(
D_{2^{\left\vert \alpha _{k}\right\vert +1}}(x)-D_{2^{{\left\vert \alpha_{k}\right\vert }}}(x)\right),
\end{equation*}%
where $\left\vert \alpha _{k}\right\vert :=\max $ $\{j\in \mathbb{N}:$ $%
\left(\alpha_{k}\right)_{j}\neq 0\}$ and $\left(\alpha_{k}\right)_{j}$
denotes $j$-th binary coefficients of real number of $\alpha _{k}\in\mathbb{N}_{+}$. Then
$\,f=\left( f_n:\text{ }n\in \mathbb{N}\right),$ where
\begin{equation*}
f_n(x):=\sum_{\left\{ k:\text{ }\left\vert \alpha
_{k}\right\vert <n\right\} }\lambda _{k}a_{k}(x)
\end{equation*}
is martingale, which belongs to $H_{p}(G)$ for any $0<p\leq 1$.

It is easy to show that
\begin{equation}  \label{3.3.10AA}
\widehat{f}(j)
\end{equation}%
\begin{equation*}
=\left\{
\begin{array}{ll}
{\lambda _{k}2^{(1/p-1)\left\vert \alpha _{k}\right\vert }},
& j\in \left\{ 2^{\left\vert \alpha _{k}\right\vert },...,
2^{\left\vert \alpha _{k}\right\vert +1}-1\right\} ,\text{ }k\in \mathbb{N}%
_{+}, \\
0, & \text{\thinspace }j\notin \bigcup\limits_{k=1}^{\infty }\left\{
2^{\left\vert \alpha _{k}\right\vert },...,2^{\left\vert \alpha_{k}\right\vert +1}-1\right\}.
\end{array}%
\right.
\end{equation*}
Let
$2^{\left\vert \alpha _{l-1}\right\vert +1}\leq j\leq 2^{\left\vert\alpha_{l}\right\vert },$\ $l\in \mathbb{N}_{+}.$
Then
\begin{eqnarray}  \label{3.3.12AA}
S_{j}f&=& S_{2^{\left\vert \alpha _{l-1}\right\vert +1}} \\ \notag
&=& \sum_{\eta =0}^{l-1}{\lambda _{\eta }2^{_{\left\vert \alpha _{\eta }\right\vert }(1/p-1)}}\left( D_{2^{\left\vert \alpha _{\eta }\right\vert
+1}}-D_{2^{\left\vert \alpha _{\eta }\right\vert }}\right) .
\end{eqnarray}

Let $2^{\left\vert \alpha _{l}\right\vert }\leq j<2^{\left\vert \alpha_{l}\right\vert +1},$ $l\in \mathbb{N}_{+}.$ Then
\begin{eqnarray}  \label{3.3.11AA}
&& S_{j}f \\ \notag
&=& S_{2^{\left\vert \alpha _{l}\right\vert }}+{\lambda
_{l}2^{(1/p-1)\left\vert \alpha _{l}\right\vert}w_{2^{\left\vert
\alpha _{l}\right\vert }}D_{j-2^{_{\left\vert \alpha _{l}\right\vert}}}} \\  \notag
&=&\sum_{\eta =0}^{l-1}{\lambda _{\eta }2^{(1/p-1)\left\vert \alpha _{\eta}\right\vert}}\left( D_{2^{{\left\vert \alpha _{\eta
}\right\vert +1}}}-D_{2^{{\left\vert \alpha _{\eta }\right\vert }}}\right) \\ \notag
&+&{\lambda _{l}2^{(1/p-1)\left\vert \alpha _{l}\right\vert }w_{2^{\left\vert \alpha _{l}\right\vert }}D_{j-2^{\left\vert \alpha_{l}\right\vert }}}.
\end{eqnarray}

Moreover, \ for \ the \ modulus \ of \ continuity \ for \ $0<p\leq 1$ \ we \  have \ the \ following \ estimation:
\begin{equation} \label{3.3.2aa0}
\omega _{H_{p}}\left( \frac{1}{2^{n}},f\right) =O\left( \sum_{\left\{ k:%
\text{ }\left\vert \alpha _{k}\right\vert \geq n\right\} }^{\infty
}\left\vert \lambda _{k}\right\vert ^{p}\right) ^{1/p},\text{ \ \ as \ \ }n\rightarrow \infty.
\end{equation}%
\end{example}

By applying Lemma \ref{lemma3.2.4} we easily obtain that the following is true (see \cite{We5}):

\begin{lemma} \label{lemma3.2.5}
\label{lemma2.3} Let $0<p\leq 1$ and $T$ be $ \sigma  $-sub-linear operator, such that, for any $p$-atom  $a$,
\begin{equation*}
\int\limits_{G}\left\vert Ta({x})\right\vert ^{p}d\mu({x}) \leq
c_{p}<\infty.
\end{equation*}
Then
\begin{equation} \label{3.2.5000}
\left\Vert Tf\right\Vert _{p}\leq c_{p}\left\Vert f\right\Vert
_{H_{p}(G)}.
\end{equation}

In addition, if $T$  is bounded from  $L_{\infty}(G)$ to $L_{\infty}(G)$ then to prove (\ref{3.2.5000}) it is suffices to show that  
\begin{equation*}
\int\limits_{\overset{-}{I}}\left\vert Ta({x})\right\vert ^{p}d\mu({x}) \leq
c_{p}<\infty,
\end{equation*}
for every $p$-atom  $a$, where $I$ denotes support of the atom $a$.
\end{lemma}

In the concrete cases the norm of Hardy martingale spaces can be calculated  by simpler formulas (for details see \cite{S}, \cite{We1} and \cite{We3}):

\begin{lemma}  \label{lemma3.2.3}
If $g\in L_{1}\left( G\right)$ and $f:=(E_{n}g:n\in \mathbb{N})$ be regular martingale, then  $H_{p}\left( G\right)\text{ \ for \ }0<p\leq1$ norm can be calculated by
\begin{equation*}
\left\Vert f\right\Vert _{H_{p}(G)}=\left\Vert \sup\limits_{n\in \mathbb{N}}|S_{2^n}g|\right\Vert _{p}.
\end{equation*}
\end{lemma}

The following lemmas are proved in \cite{tep5}, \cite{tep12}, \cite{tep13}.
\begin{lemma} \label{lemma3.2.8}
\label{example002} Let $0<p\leq 1$, $2^{k}\leq n<2^{k+1}$ and $S_{n}f$ be $n$-th partial sum with respect to the one-dimensional Walsh-Fourier series, where $f\in H_{p}(G)$. Then for any fixed $n\in \mathbb{N}$,
\begin{eqnarray*}
\left\Vert S_{n}f\right\Vert _{H_{p}(G)}^p&\leq& \left\Vert \sup_{0\leq l\leq
k}\left\vert S_{2^{l}}f\right\vert \right\Vert _{p}^p+\left\Vert
S_{n}f\right\Vert _{p}^p \\ \notag
&\leq &\left\Vert \widetilde{S}_{\#}^{\ast }f\right\Vert _{p}^p+\left\Vert
S_{n}f\right\Vert _{p}^p.
\end{eqnarray*}
\end{lemma}
{\bf Proof}:
Let consider the following martingales
\begin{equation*}
f_{\#}:=\left( S_{2^{k}}S_{n}f,\text{ }k\in\mathbb{N}_+\right)
\end{equation*}
\begin{equation*}
=\left(
S_{2^{0}},S_{2^{k}}f,\text{ }S_{n}f,...,S_{n}f,...\right),
\end{equation*}%
Hence, Lemma \ref{lemma3.2.3} immediately follows that
\begin{eqnarray*}
&&\left\Vert S_{n}f\right\Vert _{H_{p}(G)}^p \\
&\leq& \left\Vert \sup_{0\leq l\leq k}\left\vert S_{2^{l}}f\right\vert \right\Vert _{p}^p+\left\Vert
S_{n}f\right\Vert _{p}^p \\
&\leq &\left\Vert \widetilde{S}_{\#}^{\ast }f\right\Vert
_{p}^p+\left\Vert S_{n}f\right\Vert _{p}^p.
\end{eqnarray*}

Lemma is proved.
\QED

\subsection{Boundedness of subsequences of partial sums with respect to the one-dimensional Walsh-Fourier series on the martingale Hardy spaces}

\text{ \qquad } In this section we consider boundedness of subsequences of partial sums with respect to the one-dimensional Walsh-Fourier series in the martingale Hardy spaces (for details see \cite{tep12}).

\begin{theorem} \label{th4.1.1}
a) Let $0<p<1$ and $f\in H_{p}(G)$.  Then there exists an absolute constant $c_{p}$ depending only on $p,$ such that
\begin{equation*}
\text{ }\left\Vert S_{n}f\right\Vert _{H_{p}(G)}\leq c_{p}2^{d\left( n\right)
\left( 1/p-1\right) }\left\Vert f\right\Vert _{H_{p}(G)}.
\end{equation*}

b) Let $0<p<1,$ $\left\{ m_{k}:\text{ }k \in\mathbb{N}_{+}\right\}$ be non-negative, increasing  sequence of natural numbers such that
\begin{equation}
\sup_{k\in \mathbb{N}}d\left( m_{k}\right) =\infty  \label{4.1.dnk}
\end{equation}%
and $\Phi :\mathbb{N}_{+}\rightarrow \lbrack 1,\infty )$ be non-decreasing function satisfying the condition
\begin{equation}
\overline{\underset{k\rightarrow \infty }{\lim }}\frac{2^{d\left(
m_{k}\right) \left( 1/p-1\right) }}{\Phi \left( m_{k}\right) }=\infty.
\label{4.1.1010}
\end{equation}
Then there exists a martingale $f\in H_{p}(G)$ such that
\begin{equation*}
\underset{k\in \mathbb{N}}{\sup }\left\Vert \frac{S_{m_{k}}f}{\Phi \left(
m_{k}\right) }\right\Vert _{weak-L_{p}(G)}=\infty .
\end{equation*}
\end{theorem}

{\bf Proof}:
Suppose that
\begin{equation}
\left\Vert 2^{\left( 1-1/p\right) d\left( n\right) }S_{n}f\right\Vert
_{p}\leq c_{p}\left\Vert f\right\Vert _{H_{p}(G)}.  \label{4.1.11.1}
\end{equation}

By combining Lemma \ref{lemma3.2.8} and inequalities (\ref{1.S2n}) and (\ref{4.1.11.1}), since $  2^{\left( 1-1/p\right) d\left( n\right) }\leq c_{p} $ we obtain that
\begin{eqnarray} \label{4.1.11.2}
&&\left\Vert 2^{\left( 1-1/p\right) d\left( n\right) }S_{n}f\right\Vert
_{H_{p}(G)}^p \\ \notag
&\leq & \left\Vert  2^{\left( 1-1/p\right) d\left( n\right) }S_nf\right\Vert _{p}^p+2^{\left( 1-1/p\right) d\left( n\right) }\left\Vert\widetilde{S}_{\#}^{\ast}f\right\Vert_{p}^p \\ \notag
&\leq & c_{p}\left\Vert f\right\Vert _{H_{p}(G)}^p+c_{p}\left\Vert\widetilde{S}_{\#}^{\ast}f\right\Vert
_{p}^p \\ \notag
&\leq & c_{p}\left\Vert f\right\Vert _{H_{p}(G)}^p.
\end{eqnarray}
By combining Lemma \ref{lemma3.2.5} and (\ref{4.1.11.2}) it is suficies to show that 
\begin{equation}
\int\limits_{G}\left\vert 2^{\left( 1-1/p\right) d\left( n\right)
}S_{n}a\right\vert ^{p}d\mu \leq c_{p}<\infty ,  \label{4.1.25a}
\end{equation}%
for every $ p $-atom $a$, with support $I$, such that
$\mu \left(I\right) =2^{-M}$.

Without loss the generality we may assume that $p$-atom $a$ has support $I=I_{M}.$ Then it is easy to see that  $S_{n}a =0,$ where $2^{M}$ $\geq n$.
So, we may assume that $2^{M}<n$. Since $\left\Vert a\right\Vert_{\infty}\leq 2^{M/p}$ we can conclude that
\begin{eqnarray} \label{4.1.11a}
&&\left\vert 2^{\left( 1-1/p\right) d\left( n\right) }S_{n}a\left( x\right)
\right\vert \\ \notag
&\leq& 2^{\left( 1-1/p\right) d\left( n\right) }\left\Vert
a\right\Vert _{\infty }\int_{I_{M}}\left\vert D_{n}\left( x+t\right)
\right\vert d\mu \left( t\right)  \\ \notag
&\leq & 2^{M/p}2^{\left( 1-1/p\right) d\left( n\right) }\int_{I_{M}}\left\vert
D_{n}\left( x+t\right) \right\vert d\mu \left(t\right).
\end{eqnarray}

Let $x\in I_{M}$. Since $V\left( n\right) \leq 2d\left( n\right) ,$  by using first estimations of Lemma \ref{lemma2} we can conclude that
\begin{eqnarray*}
&&\left\vert 2^{\left( 1-1/p\right) d\left( n\right)}S_{n}a\right\vert \\
&\leq& 2^{M/p}2^{\left( 1-1/p\right) d\left( n\right) }V\left( n\right) \\
&\leq & 2^{M/p}d\left( n\right) 2^{\left( 1-1/p\right) d\left( n\right) }
\end{eqnarray*}%
and
\begin{eqnarray} \label{4.1.11b}
&& \int_{I_{M}}\left\vert 2^{\left( 1-1/p\right) d\left( n\right)
}S_{n}a\right\vert ^{p}d\mu \\ \notag
&\leq & d\left( n\right) 2^{\left( 1-1/p\right)
d\left( n\right) }<c_{p}<\infty .
\end{eqnarray}

Let $t\in I_{M}$ and $x\in I_{s}\backslash I_{s+1},\,$ where $0\leq s\leq
M-1<\left\langle n\right\rangle $ or $0\leq s<\left\langle n\right\rangle
\leq M-1.$ Then $x+t$ $\in I_{s}\backslash I_{s+1}$ and if we use both equality of Lemma \ref{lemma1} we get that $D_{n}\left( x+t\right) =0$ and it follows that
\begin{equation}
\left\vert 2^{\left( 1-1/p\right) d\left( n\right) }S_{n}a\left( x\right)
\right\vert =0. \label{4.1.11bbb}
\end{equation}

Let $x\in I_{s}\backslash I_{s+1},$ $\,\left\langle n\right\rangle \leq
s\leq M-1.$ Then $x+t\in I_{s}\backslash I_{s+1},$ where $t\in I_{M}$. Then by using again both equality of Lemma \ref{lemma1} we have that
\begin{equation*}
\left\vert D_{n}\left( x+t\right) \right\vert \leq
\sum_{j=0}^{s}n_{j}2^{j}\leq c2^{s}.
\end{equation*}%
If we apply again (\ref{4.1.11a}) we can conclude that
\begin{eqnarray} \label{4.1.12}
&& \left\vert 2^{\left( 1-1/p\right) d\left( n\right) }S_{n}a\left( x\right)
\right\vert \\ \notag
&\leq & 2^{\left( 1-1/p\right) d\left( n\right) }2^{M/p}\frac{2^{s}%
}{2^{M}} \\ \notag
&\leq & 2^{\left< n\right> \left(1/p-1\right) }2^{M(1/p-1)}\frac{2^{s}%
}{2^{\left| n\right| \left(1/p-1\right) }} \\ \notag
&\leq& 2^{\left\langle n\right\rangle \left( 1/p-1\right) }2^{s}.
\end{eqnarray}

By identity (\ref{2.1.2})  and inequalities  (\ref{4.1.11bbb}) and (\ref{4.1.12}) we find that
\begin{eqnarray*}
&&\int_{\overline{I_{M}}}\left\vert 2^{\left( 1-1/p\right) d\left( n\right)
}S_{n}a\left( x\right) \right\vert ^{p}d\mu \left( x\right) \\
&=&\overset{M-1}{\underset{s=\left\langle n\right\rangle }{\sum }}%
\int_{I_{s}\backslash I_{s+1}}\left\vert 2^{\left\langle n\right\rangle
\left( 1/p-1\right) }2^{s}\right\vert ^{p}d\mu \left( x\right) \\
&\leq & c\overset%
{M-1}{\underset{s=\left\langle n\right\rangle }{\sum }}\frac{2^{\left\langle
n\right\rangle \left( 1-p\right) }}{2^{s\left( 1-p\right) }}\leq
c_{p}<\infty .
\end{eqnarray*}

Now, we prove part b) of Theorem \ref{th4.1.1}. By using condition (\ref{4.1.1010}) there exists sequence of natural numbers $\left\{ \alpha _{k}:\text{ }k\in\mathbb{N}_+\right\} \subset \left\{
m_{k}:\text{ }k\in\mathbb{N}_+\right\},$ such that
\begin{equation}
\sum_{\eta =0}^{\infty }\frac{\Phi ^{p/2}\left( \alpha _{\eta }\right) }{%
2^{d\left( \alpha _{\eta }\right) \left( 1-p\right) /2}}<\infty ,
\label{4.1.12f}
\end{equation}

Let $f=\left(f_{n},n\in\mathbb{N}_+\right) \in H_{p}(G)$ be a martingale from the Example \ref{example2.2.1}, where
\begin{equation} \label{4.1.100}
\lambda _{k}=\frac{\Phi ^{1/2}\left( \alpha _{k}\right) }{2^{d\left( \alpha
_{k}\right) \left( 1/p-1\right) /2}}.
\end{equation}

Then, if we use (\ref{4.1.12f}) we obtain that condition (\ref{3.3.2aa})  is fulfilled and it follows that $f=\left(f_{n},n\in\mathbb{N}_+\right) \in H_{p}(G).$

If we apply (\ref{3.3.10AA}) when $\lambda _{k}$ are given by the formula (\ref{4.1.100}) then we get that
\begin{equation} \label{4.1.6}
\widehat{f}(j)
\end{equation}
\begin{equation*}
=\left\{
\begin{array}{c}
\Phi ^{1/2}\left( \alpha _{k}\right) 2^{\left( \left\vert \alpha
_{k}\right\vert +\left\langle \alpha _{k}\right\rangle \right) \left(
1/p-1\right) /2},\text{ \ if \thinspace \thinspace }j\in \left\{
2^{_{\left\vert \alpha _{k}\right\vert }},...,2^{_{\left\vert \alpha
_{k}\right\vert +1}}-1\right\} ,\text{ }k\in\mathbb{N}_+ \\
0,\text{ \ if \ }j\notin \bigcup\limits_{k=0}^{\infty }\left\{ 2^{_{\left\vert \alpha _{k}\right\vert}},...,2^{_{\left\vert \alpha _{k}\right\vert +1}}-1\right\} .\text{ }%
\end{array}\right.
\end{equation*}

In the view of (\ref{3.3.11AA}) when $\lambda _{k}$ are given by (\ref{4.1.100}) we get that
\begin{eqnarray} \label{4.1.6aaa}
&&\frac{S_{\alpha _{k}}f}{\Phi \left( \alpha _{k}\right) }  \\ \notag
&=&\frac{1}{\Phi \left( \alpha _{k}\right) }\sum_{\eta =0}^{k-1}\Phi
^{1/2}\left(\alpha _{\eta }\right) 2^{\left( \left\vert \alpha _{\eta
}\right\vert +\left\langle \alpha _{\eta }\right\rangle \right) \left(
1/p-1\right) /2}\left( D_{2^{\left\vert \alpha _{\eta }\right\vert
+1}}-D_{2^{\left\vert \alpha _{\eta }\right\vert}}\right) \\ \notag
&+&\frac{2^{\left( \left\vert \alpha _{k}\right\vert +\left\langle \alpha
_{k}\right\rangle \right) \left( 1/p-1\right) /2}w_{2^{\left\vert \alpha
_{k}\right\vert }}D_{\alpha _{k}-2^{\left\vert \alpha _{k}\right\vert }}}{%
\Phi ^{1/2}\left( \alpha _{k}\right) }:=I+II.
\end{eqnarray}

by using (\ref{4.1.12f}) for $I$ we have that
\begin{eqnarray} \label{4.1.8a}
&&\left\Vert I\right\Vert _{weak-L_{p}(G)}^{p}  \\ \notag
&\leq &\frac{1}{\Phi ^{p}\left( \alpha _{k}\right) }\sum_{\eta =0}^{k-1}\frac{%
\Phi ^{p/2}\left( \alpha _{\eta }\right) }{2^{d\left( \alpha _{\eta }\right)
\left( 1-p\right) /2}}\left\Vert 2^{\left\vert \alpha _{\eta }\right\vert
\left( 1/p-1\right) }\left( D_{2^{\left\vert \alpha _{\eta }\right\vert
+1}}-D_{2^{\left\vert \alpha _{\eta }\right\vert}}\right) \right\Vert
_{weak-L_{p}(G)}^{p} \\ \notag
&\leq & \frac{1}{\Phi ^{p}\left( \alpha _{k}\right) }\sum_{\eta =0}^{\infty }%
\frac{\Phi ^{p/2}\left( \alpha _{\eta }\right) }{2^{d\left( \alpha _{\eta
}\right) \left( 1-p\right) /2}}\leq c<\infty.
\end{eqnarray}

Let $x\in I_{\left\langle \alpha _{k}\right\rangle }\backslash
I_{\left\langle \alpha _{k}\right\rangle +1}.$  Since
$\left\vert \alpha _{k}\right\vert \neq\left\langle \alpha _{k}\right\rangle$ and
$\left\langle
\alpha _{k}-2^{\left\vert \alpha _{k}\right\vert }\right\rangle
=\left\langle \alpha _{k}\right\rangle.$

By using both inequalities of Lemma  \ref{lemma1} we get that
\begin{eqnarray} \label{4.1.77}
&&\left\vert D_{\alpha _{k}-2^{\left\vert \alpha _{k}\right\vert }}(x)\right\vert \\ \notag
&=&\left\vert \left( D_{2^{\left\langle \alpha _{k}\right\rangle
+1}}(x)-D_{2^{\left\langle \alpha _{k}\right\rangle }}(x)\right) +\overset{\left\vert \alpha _{k}\right\vert -1}{\underset{j=\left\langle \alpha_{k}\right\rangle +1}{\sum }}\left( \alpha _{k}\right)_{j}\left(
D_{2^{i+1}}(x)-D_{2^{i}}(x)\right) \right\vert \\ \notag
&=&\left\vert -D_{2^{\left\langle
\alpha _{k}\right\rangle }}(x)\right\vert =2^{\left\langle \alpha
_{k}\right\rangle }
\end{eqnarray}%
and
\begin{eqnarray} \label{4.1.12aa}
\left\vert II\right\vert &=&\frac{2^{\left( \left\vert \alpha _{k}\right\vert
+\left\langle \alpha _{k}\right\rangle \right) \left( 1/p-1\right) /2}}{\Phi
^{1/2}\left( \alpha _{k}\right) }\left\vert D_{\alpha _{k}-2^{\left\vert
\alpha _{k}\right\vert }}\left( x\right) \right\vert  \\ \notag
&=&\frac{2^{\left\vert \alpha _{k}\right\vert \left( 1/p-1\right)
/2}2^{\left\langle \alpha _{k}\right\rangle \left( 1/p+1\right) /2}}{\Phi
^{1/2}\left( \alpha _{k}\right) }.
\end{eqnarray}

By combining (\ref{4.1.8a}) and (\ref{4.1.12aa}) we get that
\begin{eqnarray*}
&&\left\Vert \frac{S_{\alpha _{k}}f}{\Phi \left( \alpha _{k}\right) }%
\right\Vert _{weak-L_{p}(G)}^{p} \\ \notag
&\geq & \left\Vert II\right\Vert _{weak-L_{p}(G)}^{p}-\left\Vert I\right\Vert _{weak-L_{p}(G)}^{p} \\
&\geq &\frac{2^{\left( \left\vert \alpha _{k}\right\vert \right) \left(
1/p-1\right) /2}2^{\left\langle \alpha _{k}\right\rangle \left( 1/p+1\right)
/2}}{\Phi ^{1/2}\left( \alpha _{k}\right) }\mu \left\{ x\in G:\text{ }%
\left\vert II\right\vert \geq \frac{2^{\left( \left\vert \alpha
_{k}\right\vert \right) \left( 1/p-1\right) /2}2^{\left\langle \alpha
_{k}\right\rangle \left( 1/p+1\right) /2}}{\Phi ^{1/2}\left( \alpha
_{k}\right) }\right\} ^{1/p} \\
&\geq & \frac{2^{\left( \left\vert \alpha _{k}\right\vert \right) \left(
1/p-1\right) /2}2^{\left\langle \alpha _{k}\right\rangle \left( 1/p+1\right)
/2}}{\Phi ^{1/2}\left( \alpha _{k}\right) }\left( \mu \left\{
I_{\left\langle \alpha _{k}\right\rangle }\backslash I_{\left\langle \alpha
_{k}\right\rangle +1}\right\} \right) ^{1/p} \\
&\geq & c\frac{2^{d\left( \alpha _{k}\right) \left( 1/p-1\right) /2}}{\Phi
^{1/2}\left( \alpha _{k}\right) }\rightarrow \infty, \ \text{ \ as \  }
k\rightarrow \infty.
\end{eqnarray*}

The proof of Theorem \ref{th4.1.1} is complete.
\QED

\begin{corollary} \label{cor4.1.1}
a) Let $n\in \mathbb{N}_{+}$, $0<p<1$ and $f\in H_{p}(G)$. Then there exists an absolute constant  $c_{p},$ depending only on $p$ such that
\begin{equation*}
\text{ }\left\Vert S_{n}f\right\Vert _{H_{p}(G)}\leq c_{p}\left( n\mu \left\{
\text{supp}\left( D_{n}\right) \right\} \right) ^{1/p-1}\left\Vert
f\right\Vert _{H_{p}(G)}.
\end{equation*}

\textit{b) Let} $0<p<1,$ $\left\{ m_{k}:\text{ }k\in \mathbb{N}_+\right\} $  be increasing sequence of natural numbers, such that
\begin{equation}\label{suppdnk}
\sup_{k\in \mathbb{N}}m_{k}\mu \left\{ \text{supp}\left( D_{m_{k}}\right)
\right\} =\infty
\end{equation}%
and
$\Phi :\mathbb{N}_{+}\rightarrow \lbrack 1,\infty )$ be non-decreasing function satisfying the condition
\begin{equation}
\overline{\underset{k\rightarrow \infty }{\lim }}\frac{\left( m_{k}\mu
\left\{ \text{supp}\left( D_{m_{k}}\right) \right\} \right) ^{1/p-1}}{\Phi
\left(m_{k}\right) }=\infty .  \label{4.1.12e}
\end{equation}%
\textit{Then there exists a martingale} $f\in H_{p}(G)$ \textit{such that}
\begin{equation*}
\underset{k\in \mathbb{N}}{\sup }\left\Vert \frac{S_{m_{k}}f}{\Phi \left(
m_{k}\right) }\right\Vert _{weak-L_{p}(G)}=\infty.
\end{equation*}
\end{corollary}

{\bf Proof}:
By applying both inequalities of Lemma \ref{lemma1} we get that
\begin{equation*}
I_{\left\langle n\right\rangle }\backslash I_{\left\langle n\right\rangle
+1}\subset \text{supp}\left\{ D_{n}\right\} \subset I_{\left\langle
n\right\rangle }\text{ \ and \ }2^{-\left\langle n\right\rangle -1}\leq \mu
\left\{ \text{supp}\left( D_{n}\right) \right\} \leq 2^{-\left\langle
n\right\rangle }.
\end{equation*}

Hence,
\begin{equation*}
\frac{2^{d\left( n\right) \left( 1/p-1\right) }}{4}\leq \left( n\mu \left\{
\text{supp}\left( D_{n}\right) \right\} \right) ^{1/p-1}\leq 2^{d\left(
n\right) \left( 1/p-1\right) }.
\end{equation*}%

Corollary \ref{cor4.1.1} is proved.
\QED

\begin{theorem} \label{th4.1.2}
a) Let $n\in \mathbb{N}_{+}$ and $f\in H_{1}(G).$ Then there exists an absolute constant $c,$ such that
\begin{equation*}
\left\Vert S_{n}f\right\Vert _{H_{1}(G)}\leq cV\left( n\right) \left\Vert f\right\Vert _{H_{1}(G)}.
\end{equation*}

\textit{b) Let} $\left\{ m_{k}:\text{ }k\in \mathbb{N}_{+}\right\} $ be non-negative increasing  sequence of natural numbers such that
\begin{equation} \label{4.1.vnk}
\sup_{k\in \mathbb{N}}V\left( m_{k}\right) =\infty
\end{equation}%
and $\Phi :\mathbb{N}_{+}\rightarrow \lbrack 1,\infty )$ be non-decreasing function satisfying the condition
\begin{equation}  \label{4.1.17aa}
\overline{\underset{k\rightarrow \infty }{\lim }}\frac{V\left( m_{k}\right)
}{\Phi \left( m_{k}\right) }=\infty .
\end{equation}

Then there exists a martingale $f\in H_{1}(G),$ such that
\begin{equation*}
\underset{{k\in \mathbb{N}}}{\sup }\left\Vert \frac{S_{m_{k}}f}{\Phi \left(
m_{k}\right) }\right\Vert _{1}=\infty.
\end{equation*}
\end{theorem}

{\bf Proof}:
Since
\begin{equation} \label{4.1.12k}
\left\Vert \frac{S_{n}f}{V\left( n\right)}\right\Vert _{1}\leq \left\Vert
f\right\Vert _{1}\leq \left\Vert f\right\Vert _{H_{1}(G)}
\end{equation}
by combining Lemmas \ref{lemma3.2.8} and (\ref{4.1.12k}) we can conclude that
\begin{eqnarray} \label{4.1.50}
&&\left\Vert \frac{S_{n}f}{V\left( n\right) }\right\Vert _{H_{1}(G)} \\ \notag
&\leq & \left\Vert \frac{S_{n}f}{V\left( n\right) }\right\Vert_{1}+\frac{1}{V\left( n\right) }\left\Vert\widetilde{S}_{\#}^{\ast}f\right\Vert_{1} \\ \notag
&\leq & c\left\Vert f\right\Vert _{H_{1}(G)}+c\left\Vert\widetilde{S}_{\#}^{\ast}f\right\Vert
_{1}\leq c\left\Vert f\right\Vert _{H_{1}(G)}.
\end{eqnarray}

Now prove second part of Theorem \ref{th4.1.2}. Let $\left\{ m_{k}:\text{ }k\in \mathbb{N}_{+}\right\} $ be increasing sequence of natural numbers and function $\ \Phi :\mathbb{N}_{+}\rightarrow \lbrack 1,\infty )$ satisfies conditions (\ref{4.1.vnk}) and (\ref{4.1.17aa}). Then there exists  non-negative, increasing sequence  $\left\{ \alpha _{k}:k\in \mathbb{N}_{+} \right\} \subset\left\{ m_{k}:k \in\mathbb{N}_{+}\right\} $ such that
\begin{equation} \label{4.1.2aaa}
\sum_{k=1}^{\infty}\frac{\Phi ^{1/2}\left(\alpha_{k}\right)}{
V^{1/2}\left(\alpha_{k}\right)}\leq \beta <\infty.
\end{equation}

Let  $f=\left(f_{n},n\in \mathbb{N}_+\right)$ be martingale from Example \ref{example2.2.1}, where
\begin{equation} \label{4.1.101}
\text{\ }\lambda _{k}=\frac{\Phi ^{1/2}\left( \alpha _{k}\right) }{%
V^{1/2}\left( \alpha _{k}\right) }.
\end{equation}

By applying condition (\ref{4.1.2aaa}) we can conclude that condition (\ref{3.3.2aa}) is fulfilled and it follows that $f=\left(f_{n},n\in \mathbb{N}_+\right) \in H_{1}(G).$

In the view of (\ref{3.3.10AA}) when $\lambda _{k}$ are given by (\ref{4.1.101}) we get that
\begin{equation} \label{4.1.5aa}
\widehat{f}(j)=\left\{
\begin{array}{ll}
\frac{\Phi ^{1/2}\left( \alpha _{k}\right) }{V^{1/2}\left( \alpha
_{k}\right) }, & \text{ \  if \ }j\in \left\{
2^{_{\left\vert \alpha _{k}\right\vert }},...,2^{_{\left\vert \alpha
_{k}\right\vert +1}}-1\right\} ,\text{ }k=0,1,... \\
0, & \text{ \ if \ }j\notin
\bigcup\limits_{k=0}^{\infty }\left\{ 2^{_{\left\vert \alpha _{k}\right\vert
}},...,2^{_{\left\vert \alpha _{k}\right\vert +1}}-1\right\} .\text{ }%
\end{array}%
\right.
\end{equation}

Analogously to (\ref{4.1.6aaa}) if we apply (\ref{3.3.11AA}) when $\lambda _{k}$ are given by (\ref{4.1.101}) we get that
\begin{eqnarray*}
 S_{\alpha _{k}}f &=& \sum_{\eta =0}^{k-1}\frac{\Phi ^{1/2}\left( \alpha _{\eta
}\right) }{V^{1/2}\left( \alpha _{\eta }\right) }\left( D_{2^{\left\vert
\alpha _{\eta }\right\vert +1}}-D_{2^{\left\vert \alpha _{\eta }\right\vert
}}\right) \\
&+&\frac{\Phi ^{1/2}\left( \alpha _{k}\right) }{V^{1/2}\left( \alpha
_{k}\right) }w_{2^{\left\vert \alpha _{k}\right\vert }}D_{\alpha
_{k}-2^{\left\vert \alpha _{k}\right\vert }}.
\end{eqnarray*}

By applying first estimation of Lemma \ref{lemma2} and (\ref{4.1.2aaa}) we can conclude that
\begin{eqnarray*}
&&\left\Vert \frac{S_{\alpha _{k}}f}{\Phi \left( \alpha _{k}\right) }%
\right\Vert _{1} \\
&\geq &\frac{\Phi ^{1/2}\left( \alpha _{k}\right) }{\Phi
\left( \alpha _{k}\right) V^{1/2}\left( \alpha _{k}\right) }\left\Vert
D_{\alpha _{k}-2^{\left\vert \alpha _{k}\right\vert }}\right\Vert _{1} \\
&-&\frac{1}{\Phi \left( \alpha _{k}\right) }\sum_{\eta =0}^{k-1}\frac{\Phi
^{1/2}\left( \alpha _{\eta }\right) }{V^{1/2}\left( \alpha _{\eta }\right) }\left\Vert  D_{2^{\left\vert
\alpha _{\eta }\right\vert +1}}-D_{2^{\left\vert \alpha _{\eta }\right\vert}}\right\Vert _{1} \\
&\geq &\frac{V\left( \alpha _{k}-2^{\left\vert \alpha _{k}\right\vert }\right)
\Phi ^{1/2}\left( \alpha _{k}\right) }{8\Phi \left( \alpha _{k}\right)
V^{1/2}\left( \alpha _{k}\right) } \\
&-&\frac{1}{\Phi \left( \alpha _{k}\right) }%
\sum_{\eta =0}^{\infty }\frac{\Phi ^{1/2}\left( \alpha _{\eta }\right) }{%
V^{1/2}\left( \alpha _{\eta }\right) } \\
&\geq & \frac{cV^{1/2}\left( \alpha _{k}\right) }{\Phi ^{1/2}\left( \alpha
_{k}\right) }\rightarrow \infty ,\text{ \ as \ }k\rightarrow \infty.
\end{eqnarray*}

Theorem \ref{th4.1.2} is proved.
\QED

\begin{corollary} \label{cor4.1.2}
Let $n\in \mathbb{N}$, $0<p\leq 1$ and $f\in H_{p}(G)$. Then there exists an absolute constant $ c_{p}, $ depending only on $ p, $ such that
\begin{equation} \label{4.1.sn2n20}
\left\Vert S_{2^{n}}f\right\Vert _{H_{p}(G)}\leq c_{p}\left\Vert
f\right\Vert _{H_{p}(G)}.
\end{equation}
\end{corollary}

{\bf Proof}:
To prove Theorem \ref{cor4.1.2} we only have to show that
\begin{equation*}
\left\vert 2^{n}\right\vert =n,\text{ \  \ } \left\langle
2^{n}\right\rangle =n-1\text{ \ და \ }d\left(
2^{n}\right)=0.
\end{equation*}
By applying first part of Theorem \ref{th4.1.1} we immediately get that (\ref{4.1.sn2n20}) for any $0<p\leq 1$ and proof of Corollary \ref{cor4.1.2} is proved.
\QED

\begin{corollary} \label{cor4.1.3}
Let $n\in \mathbb{N}$, $0<p\leq 1$ and $f\in H_{p}(G)$. Then there exists an absolute constant $ c_p, $  depending only on $p$ such that 
\begin{equation} \label{4.1.sn2n2}
\left\Vert S_{2^{n}+2^{n-1}}f\right\Vert _{H_{p}(G)}\leq c_{p}\left\Vert
f\right\Vert _{H_{p}(G)}.
\end{equation}
\end{corollary}

{\bf Proof}:
Since
\begin{equation*}
\left\vert 2^{n}+2^{n-1}\right\vert =n,\left\langle
2^{n}+2^{n-1}\right\rangle =n-1\text{ \ and \ }d\left(
2^{n}+2^{n-1}\right)=1
\end{equation*}
by first part of Theorem \ref{th4.1.1} we get that  (\ref{4.1.sn2n2}) holds, for any $0<p\leq 1$ and proof of Corollary \ref{cor4.1.3} is complete.
\QED

\begin{corollary} \label{cor4.1.4}
Let $n\in \mathbb{N}$ and $0<p<1.$ Then there exists a martingale $f\in
H_{p}(G),$ such that
\begin{equation} \label{4.1.sn2n1}
\underset{n\in \mathbb{N}}{\sup }\left\Vert S_{2^{n}+1}f\right\Vert
_{weak-L_{p}(G)}=\infty.
\end{equation}

On the other hand, there exists an absolute constant $ c $, such that

\begin{equation} \label{4.1.sn2n22}
\left\Vert S_{2^{n}+1}f\right\Vert _{H_{1}(G)}\leq c\left\Vert
f\right\Vert _{H_{1}(G)}.
\end{equation}

\end{corollary}

{\bf Proof}:
Since
\begin{equation}
\left\vert 2^{n}+1\right\vert =n,\left\langle 2^{n}+1\right\rangle =0\text{
\ და \ }d\left( 2^{n}+1\right) =n.  \label{4.1.cor1}
\end{equation}%
by applying second part of Theorem \ref{th4.1.1} we get that there exists a martingale $f=\left(f_{n},n\in \mathbb{N}_+\right) \in H_{p}(G),$ for $0<p<1,$ such that (\ref{4.1.sn2n1}) holds.

On the other hand, proof of  (\ref{4.1.sn2n22}) follows simple observation that $$ V( 2^{n}+1)=4<\infty. $$

Corollary \ref{cor4.1.4} is proved.
\QED

\bigskip

\subsection{Modulus of continuity and convergence in norm of subsequences of partial sums with respect to the one-dimensional Walsh-Fourier series on the martingale Hardy spaces}

\text{ \qquad } In this section we apply Theorem \ref{th4.1.1} and  Theorem \ref{th4.1.2} to find necessary and sufficient conditions for modulus of continuity, for which subsequences of partial sums with respect to the one-dimensional Walsh-Fourier series are bounded in the martingale Hardy spaces.

First, we prove the following estimation:

\begin{theorem} \label{theorem4.2.1}
Let $ n\in \mathbb{N}_+ $ and $2^{k}<n\leq 2^{k+1}.$ Then there exists an absolute constant $c_{p},$ depending only on $p$ such that
\begin{equation} \label{4.2.sn1}
\left\Vert S_{n}f-f\right\Vert _{H_{p}(G)}\leq c_{p}2^{d\left( n\right) \left(1/p-1\right)}\omega _{H_{p}(G)}\left( \frac{1}{2^{k}},f\right),
\text{ \ \ \ } (f\in H_p(G))\text{ \ \ \ }\left( 0<p<1\right)
\end{equation}%
and
\begin{equation} \label{4.2.sn2}
\left\Vert S_{n}f-f\right\Vert _{H_{1}(G)}\leq c_{1}V\left( n\right) \omega
_{H_{1}(G)}\left( \frac{1}{2^{k}},f\right),\text{ \ \ \ } (f\in H_1(G)) .
\end{equation}
\end{theorem}

{\bf Proof}:
Let $0<p<1$ and $2^{k}<n\leq 2^{k+1}.$ By applying first part of Theorem \ref{th4.1.1} we get that
\begin{eqnarray} \label{4.2.1000}
&&\left\Vert S_{n}f-f\right\Vert _{H_{p}(G)}^p \\ \notag
&\leq& c_{p}\left\Vert
S_{n}f-S_{2^{k}}f\right\Vert _{H_{p}(G)}^p+c_{p}\left\Vert
S_{2^{k}}f-f\right\Vert _{H_{p}(G)}^p  \\ \notag
&=&c_{p}\left\Vert S_{n}\left( S_{2^{k}}f-f\right) \right\Vert
_{H_{p}(G)}^p+c_{p}\left\Vert S_{2^{k}}f-f\right\Vert _{H_{p}(G)}^p \\ \notag
&\leq & c_{p}\left( 1+2^{d\left( n\right) \left( 1-p\right) }\right) \omega
_{H_{p}(G)}^p\left( \frac{1}{2^{k}},f\right) \\ \notag
&\leq & c_{p}2^{d\left( n\right) \left(
1-p\right) }\omega _{H_{p}(G)}^p\left( \frac{1}{2^{k}},f\right).
\end{eqnarray}

The proof of (\ref{4.2.sn2}) is analogical to (\ref{4.2.sn1}). Analogously to (\ref{4.2.sn1}) we can also prove estimation (\ref{4.2.sn2}). So, we leave out the details.

Theorem \ref{theorem4.2.1} is proved.
\QED

\begin{theorem} \label{th4.2.2}
a) Let $k\in\mathbb{N}_+$, $0<p<1,$ $f\in H_{p}(G)$ and $\{m_{k}:k\in\mathbb{N}_+\}$ be increasing sequence of natural numbers, such that
\begin{equation}  \label{4.2.18a}
\omega _{H_{p}(G)}\left( \frac{1}{2^{\left\vert m_{k}\right\vert}},f\right)
=o\left(\frac{1}{2^{d\left(m_{k}\right)\left( 1/p-1\right)}}\right)
\text{ \ as \ }k\rightarrow \infty.
\end{equation}%
Then
\begin{equation} \label{4.2.con1}
\left\Vert S_{m_{k}}f-f\right\Vert _{H_{p}(G)}\rightarrow 0\text{ \ as \ }%
k\rightarrow \infty .
\end{equation}

b) Let $\{m_{k}:k\in\mathbb{N}_+\}$ be increasing sequence of natural numbers,such that condition (\ref{4.1.dnk}) is fulfilled. Then there exists a martingale  $f\in H_{p}(G)$ and increasing sequence of natural numbers  $\{\alpha _{k}:k\in\mathbb{N}_+\}\subset \{m_{k}:k\in\mathbb{N}_+\},$ such that
\begin{equation*}
\omega _{H_{p}(G)}\left( \frac{1}{2^{\left\vert \alpha _{k}\right\vert }}%
,f\right) =O\left( \frac{1}{2^{d\left( \alpha _{k}\right) \left(
1/p-1\right) }}\right) \text{ \ as \ }k\rightarrow \infty \text{\ }
\end{equation*}%
and
\begin{equation} \label{4.2.con11}
\limsup\limits_{k\rightarrow \infty }\left\Vert S_{\alpha
_{k}}f-f\right\Vert _{weak-L_{p}(G)}>c_{p}>0,\text{ \ as  \ }k\rightarrow \infty,
\end{equation}%
where $c_{p}$ is an absolute constant depending only on  $p$.
\end{theorem}

{\bf Proof}:
Let $0<p<1$, $f\in H_{p}(G)$ and $\{m_{k}:k\in\mathbb{N}_+\}$ be increasing sequence of natural numbers, such that condition  (\ref{4.2.18a}) is fulfilled. By combining Theorem \ref{theorem4.2.1} and estimation (\ref{4.2.sn1}) we get that (\ref{4.2.con1}) holds true.

Now, prove second part of Theorem \ref{th4.2.2}. In the view of  (\ref{4.1.dnk}) we simply get that there exists sequence $\{\alpha _{k}:k\in\mathbb{N}_+\}\subset \{m_{k}:k\in\mathbb{N}_+\},$ such that
\begin{equation} \label{4.2.4.18}
\ 2^{d\left( \alpha _{k}\right) }\uparrow \infty ,\,\,\,\,\text{ \ as\ }\ \
k\rightarrow \infty \text{, \ \ \ }2^{2\left( 1/p-1\right) d\left( \alpha
_{k}\right) }\leq 2^{\left( 1/p-1\right) d\left( \alpha _{k+1}\right) }.
\end{equation}%

Let  $f=\left(f_{n},n\in \mathbb{N}\right)$ be a martingale from Example \ref{example2.2.1}, such that
\begin{equation} \label{4.2.4.188}
\lambda_{i}={2^{-(1/p-1)d\left(\alpha_{i}\right)}}.
\end{equation}

By applying (\ref{4.2.4.18}) we obtain that condition (\ref{3.3.2aa}) is fulfilled and it follows that  $\ f\in H_{p}(G).$

By applying (\ref{3.3.10AA}), whene $\lambda _{k}$ are given by (\ref{4.2.4.188}), then
\begin{equation} \label{4.2.4.22}
\widehat{f}(j)=\left\{
\begin{array}{ll}
2^{\left( 1/p-1\right) \left\langle \alpha _{k}\right\rangle }, & \text{ \ if
\ }j\in \left\{ 2^{_{\left\vert \alpha _{k}\right\vert
}},...,2^{_{\left\vert \alpha _{k}\right\vert +1}}-1\right\} ,\text{ }%
k\in\mathbb{N}_+, \\
0, & \text{\ if \ }j\notin
\bigcup\limits_{k=0}^{\infty }\left\{ 2^{_{\left\vert \alpha _{k}\right\vert
}},...,2^{_{\left\vert \alpha _{k}\right\vert +1}}-1\right\} .\text{ }%
\end{array}
\right.
\end{equation}

By combining (\ref{4.2.4.18}) and (\ref{3.3.2aa0}) we have that
\begin{eqnarray} \label{4.2.4.21}
&&\omega _{H_{p}(G)}(\frac{1}{2^{\left\vert \alpha _{k}\right\vert }},f) \\ \notag
&\leq&\sum\limits_{i=k}^{\infty }\frac{1}{2^{\left( 1/p-1\right) d\left( \alpha_{i}\right) }} \\ \notag
&=& O\left( \frac{1}{2^{\left( 1/p-1\right) d\left( \alpha
_{k}\right) }}\right) ,\text{ \ as \ }k\rightarrow \infty.
\end{eqnarray}

By using (\ref{4.1.77}) we get that
\begin{equation*}
\left\vert D_{\alpha _{k}-2^{\left\langle \alpha
_{k}\right\rangle }}\right\vert \geq 2^{\left\langle \alpha
_{k}\right\rangle },\text{ \ \ \ \ where \ \ \ }I_{\left\langle \alpha_{k}\right\rangle }\backslash I_{\left\langle \alpha _{k}\right\rangle +1}.
\end{equation*}

In the view of (\ref{3.3.11AA}) we can conclude that
\begin{equation*}
S_{\alpha _{k}}f=S_{2^{|\alpha _{k}|}}f+2^{\left( 1/p-1\right) \left\langle \alpha _{k}\right\rangle }{w_{2^{|\alpha _{k}|}}D_{{\alpha _{k}}-2^{|\alpha _{k}|}}},
\end{equation*}

Since
\begin{eqnarray*}
&&\Vert D_{\alpha _{k}}\Vert _{weak-L_{p}(G)} \\
&\geq & 2^{\left\langle \alpha
_{k}\right\rangle }\mu \left\{ x\in I_{\left\langle \alpha _{k}\right\rangle} \backslash I_{\left\langle \alpha _{k}\right\rangle +1}:\text{ }\left \vert D_{\alpha _{k}}\right\vert \geq 2^{\left\langle \alpha _{k}\right\rangle}\right\} ^{1/p} \\
&\geq & 2^{\left\langle \alpha _{k}\right\rangle }\left( \mu \left\{
I_{\left\langle \alpha _{k}\right\rangle }\backslash I_{\left\langle \alpha_{k}\right\rangle +1}\right\} \right) ^{1/p}\geq 2^{\left\langle \alpha_{k}\right\rangle \left(1-1/p\right) },
\end{eqnarray*}
if we apply (\ref{1.S2n000}) (see also Theorem T2) we obtain that
\begin{eqnarray*}
&&\Vert f-S_{\alpha _{k}}f\Vert_{weak-L_{p}(G)}^p \\
&\geq& 2^{\left( 1-p\right) \left\langle \alpha _{k}\right\rangle }\Vert {w_{2^{|\alpha _{k}|}}D_{{\alpha _{k}}-2^{|\alpha _{k}|}}}\Vert _{weak-L_{p}(G)}^p \\
&-&\Vert f-S_{2^{|\alpha _{k}|}}f\Vert_{weak-L_{p}(G)}^p \\
&\geq & c-o(1)>c>0, \text{\ \ \ as  \ \ }k\rightarrow \infty .
\end{eqnarray*}
Proof of Theorem \ref{th4.2.2} is complete.
\QED

\begin{corollary} \label{cor4.2.1}
a) Let $0<p<1,$ $f\in H_{p}(G)$ and $\{m_{k}:k\in \mathbb{N}_+\}$ be increasing sequence of natural numbers, such that
\begin{equation} \label{4.2.cond2}
\omega_{H_{p}(G)}\left(\frac{1}{2^{\left\vert m_{k}\right\vert }},f\right)
=o\left( \frac{1}{\left(m_{k}\mu \left( \text{supp}D_{m_{k}}\right) \right)
^{1/p-1}}\right) \text{ \ as\ }k\rightarrow \infty .
\end{equation}%
Then (\ref{4.2.con1}) holds.

b) Let $\{m_{k}:k\in \mathbb{N}_+\}$ be increasing sequence of natural numbers, such that
\begin{equation} \label{4.2.con1ab}
\sup_{k\in \mathbb{N}_+}m_{k}\mu \left\{ \text{supp}\left( D_{m_{k}}\right)
\right\} =\infty.
\end{equation}
Then there exist a martingale $f\in H_{p}(G)$ and sequence $\{\alpha _{k}:k\in \mathbb{N}_+\}\subset \{m_{k}:k\in \mathbb{N}_+\},$ such that
\begin{equation*}
\omega _{H_{p}(G)}\left( \frac{1}{2^{\left\vert \alpha _{k}\right\vert }},f\right) =O\left( \frac{1}{\left( \alpha _{k}\mu \left( \text{supp}D_{\alpha _{k}}\right) \right)^{1/p-1}}\right) \text{\ \ as \ }k\rightarrow\infty
\end{equation*}%
and (\ref{4.2.con11}) holds.
\end{corollary}

\begin{theorem} \label{th4.2.3}
a) Let $f\in H_{1}(G)$ and  $\{m_{k}:k\in \mathbb{N}_+\}$ be increasing sequence of natural numbers, such that
\begin{equation} \label{4.2.cond1}
\omega _{H_{1}(G)}\left( \frac{1}{2^{\left\vert m_{k}\right\vert }},f\right)=o\left( \frac{1}{V\left( m_{k}\right) }\right) \text{ \ as \ }k\rightarrow \infty.
\end{equation}%
Then
\begin{equation} \label{4.2.cond1a}
\left\Vert S_{m_{k}}f-f\right\Vert _{H_{1}(G)}\rightarrow 0\text{ \ as \ }%
k\rightarrow \infty .
\end{equation}

b) Let $\{m_{k}:k\in \mathbb{N}_+\}$ be increasing sequence of natural numbers, such that condition (\ref{4.1.vnk}) is fulfilled. Then there exists a martingale $f\in H_{1}(G)$ and increasing sequence of natural numbers
$\{\alpha _{k}:k\in \mathbb{N}_+\}\subset \{m_{k}:k\in \mathbb{N}_+\}$ such that
\begin{equation*}
\omega _{H_{1}(G)}\left( \frac{1}{2^{\left\vert \alpha _{k}\right\vert }},f\right) =O\left( \frac{1}{V\left( \alpha _{k}\right) }\right) \text{ \ as
\ }k\rightarrow \infty
\end{equation*}%
and
\begin{equation} \label{cond10}
\limsup\limits_{k\rightarrow \infty }\left\Vert S_{\alpha
_{k}}f-f\right\Vert _{1}>c>0\,\,\,\text{\ as \ }
k\rightarrow \infty,
\end{equation}
where $c$ is an absolute constant.
\end{theorem}

{\bf Proof}:
Let $f\in H_{1}(G)$ and $\{m_{k}:k\in \mathbb{N}_+\}$  be increasing sequence of natural numbers, such that (\ref{4.2.cond1}). By applying Theorem \ref{theorem4.2.1} we get that condition (\ref{4.2.cond1a}) is fulfilled.

Now, we prove second part of Theorem \ref{th4.2.3}. By applying (\ref{4.1.vnk}) we conclude that there exists sequence  $\{\alpha _{k}:k\in \mathbb{N}_+\}\subset
\{m_{k}:k\in \mathbb{N}_+\}$, such that
\begin{equation}
V(\alpha _{k})\uparrow \infty , \text{ \ \ as \ \ }k\rightarrow \infty \text{ \ \ and \ \ }V^{2}(\alpha _{k})\leq V(\alpha _{k+1}) \text{ \ \ }k\in \mathbb{N}_+.  \label{4.2.4.7}
\end{equation}

Let  $f=\left(f_{n},n\in \mathbb{N}_+\right)$ be a martingale from the Example \ref{example2.2.1}, where
\begin{equation*}
\lambda_{k}=\frac{1}{V(\alpha _{k})}.
\end{equation*}%

By applying (\ref{4.2.4.7}) we conclude that (\ref{3.3.2aa}) is fulfilled and we conclude that $f=\left( f_{n},n\in
\mathbb{N}_+\right) \in H_{1}(G).$

In the view of (\ref{3.3.10AA}) we have that
\begin{equation}
\widehat{f}(j)=\left\{
\begin{array}{ll}
\frac{1}{V(\alpha _{k})}, & \text{ თუ \thinspace \thinspace }j\in \left\{
2^{_{\left\vert \alpha _{k}\right\vert }},...,2^{_{\left\vert \alpha
_{k}\right\vert +1}}-1\right\} ,\text{ }k=0,1,... \\
0, & \text{\ თუ \thinspace \thinspace \thinspace }j\notin
\bigcup\limits_{k=0}^{\infty }\left\{ 2^{_{\left\vert \alpha _{k}\right\vert}},...,2^{_{\left\vert \alpha _{k}\right\vert +1}}-1\right\} .\text{ }%
\end{array}%
\right.  \label{4.13}
\end{equation}

According to (\ref{3.3.2aa0}) we get that
\begin{eqnarray} \label{4.2.4.12}
&& w_{H_{1}(G)}(1/2^n,f) \\ \notag
&=&\left\Vert f-S_{2^{n}}f\right\Vert _{H_{1}(G)}
\leq \sum\limits_{i=n+1}^{\infty }%
\frac{1}{V(\alpha _{i})} \\ \notag
&=& O\left( \frac{1}{V(\alpha _{n})}\right), \text{\ \ as \ \ }n\rightarrow \infty.   \\ \notag
\end{eqnarray}

By applying (\ref{3.3.11AA}) we can conclude that
\begin{equation*} \label{4.2.4.1222222}
S_{\alpha _{k}}f=S_{2^{|\alpha _{k}|}}f+\frac{w_{2^{|\alpha _{k}|}}D_{{\alpha _{k}}-2^{|\alpha _{k}|}}}{V(\alpha _{k})},
\end{equation*}

If we use (\ref{1.S2n000}) and Theorem T2  we get that
\begin{eqnarray*}
&&\Vert f-S_{\alpha _{k}}f\Vert _{1} \\
&\geq & \Vert \frac{w_{2^{|\alpha _{k}|}}D_{{\alpha _{k}}-2^{|\alpha _{k}|}}}{V(\alpha _{k})}\Vert _{1}-\Vert f-S_{2^{|\alpha _{k}|}}f\Vert _{1} \\
&\geq & \frac{V(\alpha _{k}-2^{|\alpha _{k}|})}{8V(\alpha _{k})}-o(1)>c>0, \text{\ \ as \ \ }k\rightarrow \infty .
\end{eqnarray*}

The proof of Theorem \ref{th4.2.3} is proved.
\QED

Theorem \ref{theorem5.2.2} follows the following corollaries which are \cite{tep7}:

\begin{corollary}\label{corollary5.2.2snsn}
a) Let $0<p<1,$ $f\in H_{p}(G)$ and
\begin{equation*}
\omega _{H_{p}(G)}\left( 1/2^{k},f\right) =o\left(
1/2^{k( 1/p-1)}\right) ,\text{ as  }
k\rightarrow \infty .
\end{equation*}%
Then
\begin{equation*}
\left\Vert S_{k}f-f\right\Vert _{H_{p}(G)}\rightarrow
0, \text{ \ as \ }k\rightarrow \infty .
\end{equation*}

b) There exists a martingale $f\in H_{p}(G)$ $\left( 0<p<1\right),$ such that
\begin{equation*}
\omega _{H_{p}(G)}\left( 1/2^{k},f\right) =O\left(
1/2^{k\left( 1/p-1\right) }\right),\text{ \ as \ }%
k\rightarrow \infty
\end{equation*}%
and
\begin{equation*}
\left\Vert S_{k}f-f\right\Vert _{weak-L_{p}(G)}\nrightarrow
0, \text{ \ as \ }k\rightarrow\infty .
\end{equation*}
\end{corollary}

\begin{corollary}\label{corollary5.2.2snsn2}
a) Let  $f\in H_{1}(G)$ and
\begin{equation*}
\omega _{H_{1}(G)}\left( 1/2^{k},f\right) =o\left(\frac{1}{ k}\right) ,\text{ as \ }
k\rightarrow \infty .
\end{equation*}%
Then
\begin{equation*}
\left\Vert S_{k}f-f\right\Vert _{H_{1}(G)}\rightarrow
0,\text{\ as \ }k\rightarrow\infty .
\end{equation*}

b) There exists a martingale $f\in H_{1}(G),$ such that
\begin{equation*}
\omega _{H_{1}(G)}\left( 1/2^{k},f\right)=O\left(\frac{1}{k}\right),\text{ \ as \ }k\rightarrow \infty
\end{equation*}
and
\begin{equation*}
\left\Vert S_kf-f\right\Vert_1\nrightarrow 0,\text{ \ as \ }k\rightarrow \infty.
\end{equation*}
\end{corollary}
\newpage

\section{Fej\'er means with respect to the one-dimensional Walsh-Fourier series on the martingale Hardy spaces}
\subsection{Basic notations}

\text{ \qquad } Fot the one-dimensional case Fej\'er means with respect to the one-dimensional Walsh-Fourier series $\sigma _{n}$ is defined by:

\begin{eqnarray*}
\qquad \sigma _{n}f(x) &:&=\frac{1}{n}\sum_{k=1}^{n}S_{k}f(x)\text{\qquad\ \ \ }\left( \text{ }n\in \mathbb{N}_{+}\right).
\end{eqnarray*}

The following equality is true (for details see \cite{1} and \cite{sws}):

\begin{equation*}
\sigma _{n}f\left( x\right) =\frac{1}{n}\overset{n-1}{\underset{k=0}{\sum }}%
\left( D_{k}\ast f\right) \left( x\right)
\end{equation*}%
\begin{equation*}
=\left( f\ast K_{n}\right) \left( x\right) =\int_{G_{m}}f\left( t\right)
K_{n}\left( x-t\right) d\mu \left( t\right) .
\end{equation*}
where
\begin{eqnarray*}
K_{n}(x) &:&=\frac{1}{n}\overset{n}{\underset{k=1}{\sum }}D_{k}(x)\text{ \qquad
\thinspace }\left( \text{ }n\in \mathbb{N}_{+}\text{ }\right)
\end{eqnarray*}

In the literature $ K_{n} $ is called $ n $-th Fej\'er kernel. 

We also define the following maximal operators
\begin{eqnarray*}
{\sigma }^{\ast}f&=&\sup_{n\in \mathbb{N}}\left\vert \sigma_{{n}}f\right\vert  \\
\widetilde{\sigma }_{\#}^{\ast}f&=&\sup_{n\in \mathbb{N}}\left\vert \sigma_{2^{n}}f\right\vert.
\end{eqnarray*}.

For any natural number $n\in \mathbb{N}$ we also need the following expression
\begin{equation*}
n=\sum_{i=1}^{s}2^{n_{i}}, \text{ \qquad
\thinspace } n_{1}<n_{2}<...<n_{s}.
\end{equation*}

Set
\begin{equation*}
n^{\left( i\right) }:=2^{n_{1}}+...+2^{n_{i-1}},\text{ \ }i=2,...,s
\end{equation*}%
and
\begin{equation*}
\mathbb{A}_{0,2}:=\left\{ n\in \mathbb{N}:\text{ }n=2^{0}+2^{2}+%
\sum_{i=3}^{s_{n}}2^{n_{i}}\right\} .
\end{equation*}

Then, for any natural number $n\in \mathbb{N}$ there exists numbers
\begin{equation*}
0 \leq l_{1}\leq m_{1}\leq l_{2}-2<l_{2}\leq m_{2}\leq ...\leq l_{s}-2<l_{s}\leq m_{s}
\end{equation*}
such that it can be written as
\begin{equation*}
n=\sum_{i=1}^{s}\sum_{k=l_{i}}^{m_{i}}2^{k},
\end{equation*}
where $s$ is depending on $n$.

It is evident that
\begin{equation*}
s\leq V\left(n\right) \leq 2s+1.
\end{equation*}

\subsection{Auxiliary lemmas}

\text{ \qquad } The following equality and estimation of Fej\'er kernels with respect to the one-dimensional Walsh-Fourier series is proved in \cite{sws}:

\begin{lemma} \label{lemma3}

Let $n\in \mathbb{N}$ and $n=\sum_{i=1}^{s}2^{n_{i}},$ $n_{1}<n_{2}<...<n_{s}$. Then
\begin{equation*}
nK_{n}=\sum_{r=1}^{s}\left( \underset{j=r+1}{\overset{s}{\prod }}%
w_{2^{n_{j}}}\right) 2^{n_{r}}K_{2^{n_{r}}}+\sum_{t=2}^{s}\left( \underset{%
j=t+1}{\overset{s}{\prod }}w_{2^{n_{j}}}\right) n^{\left( t\right)
}D_{2^{n_{t}}},
\end{equation*}
and
\begin{equation*}
\sup_{n\in \mathbb{N}}\int_{G}\left\vert K_{n}\left( x\right) \right\vert d\mu \left(
x\right) \leq c<\infty ,
\end{equation*}%
where $ c $ is an absolute constant.
\end{lemma}

The following equality is proved in \cite{sws} (see also \cite{gat}):

\begin{lemma} \label{lemma4}
Let $n>t$ and $t,n\in \mathbb{N}$. Then we have the following expression for $2^{n}$-th Fej\'er kernels with respect to the one-dimensional Walsh-Fourier series:
\begin{equation*}
K_{2^{n}}\left( x\right) =\left\{
\begin{array}{c}
\text{ }2^{t-1},\text{\  if \ }x\in I_{n}\left( e_{t}\right) , \\
\frac{2^{n}+1}{2},\text{ \  if \ }x\in I_{n},\text{\ } \\
0,\text{ \  otherwise.\  }
\end{array}%
\right.
\end{equation*}%

\end{lemma}

The following estimation is proved by Goginava \cite{GoSzeged}:
\begin{lemma} \label{lemma5}
Let  $x\in I_{l+1}\left(e_{k}+e_{l}\right), \ \ k=0,...,M-2, \ \ l=0,...,M-1.$
Then
\begin{equation*}
\int_{I_{M}}\left\vert K_{n}\left( x+t\right) \right\vert d\mu \left(
t\right) \leq \frac{c2^{l+k}}{n2^{M}},\text{ \  where \ }n>2^{M}.
\end{equation*}

Let $x\in I_{M}\left(e_{k}\right),m=0,...,M-1.$ Then
\begin{equation*}
\int_{I_{M}}\left\vert K_{n}\left( x+t\right) \right\vert d\mu \left(t\right) \leq \frac{c2^{k}}{2^{M}},\text{ \ \ for \ \ }n>2^{M},
\end{equation*}
where $ c $ is an absolute constant.
\end{lemma}

The following estimations of Fej\'er kernels with respect to the one-dimensional Walsh-Fourier series is proved in \cite{tep13}:

\begin{lemma} \label{lemma6}
Let
 $$n=\sum_{i=1}^{r}\sum_{k=l_{i}}^{m_{i}}2^{k},$$
where
$$m_{1}\geq l_{1}>l_{1}-2\geq m_{2}\geq l_{2}>l_{2}-2>...>m_{s}\geq l_{s}\geq0.$$
Then
\begin{equation*}
\left\vert nK_{n}\right\vert \leq c\sum_{A=1}^{r}\left( 2^{l_{A}}\left\vert
K_{2^{l_{A}}}\right\vert +2^{m_{A}}\left\vert K_{2^{m_{A}}}\right\vert
+2^{l_{A}}\sum_{k=l_{A}}^{m_{A}}D_{2^{k}}\right) +cV\left( n\right),
\end{equation*}
where $ c $ is an absolute constant.
\end{lemma}

{\bf Proof}:
Let $$n=\sum_{i=1}^{r}2^{n_{i}}, n_{1}>n_{2}>...>n_{r}\geq 0.$$ By using Lemma \ref{lemma3} for $n$-th Fej\'er kernels we can conclude that
\begin{eqnarray*}
nK_{n}&=&\sum_{A=1}^{r}\left( \underset{j=1}{\overset{A-1}{\prod }}%
w_{2^{n_{j}}}\right) \left( \left( 2^{n_{A}}K_{2^{n_{A}}}+\left(
2^{n_{A}}-1\right) D_{2^{n_{A}}}\right) \right) \\
&-&\sum_{A=1}^{r}\left( \underset{j=1}{\overset{A-1}{\prod }}%
w_{2^{n_{j}}}\right) \left( 2^{n_{A}}-1-n^{\left( A\right) }\right)
D_{2^{n_{A}}}=I_{1}-I_{2}.
\end{eqnarray*}

For $I_{1}$ we have the following equality
\begin{eqnarray*}
I_{1}&=&\sum_{v=1}^{r}\left( \underset{j=1}{\overset{v-1}{\prod }}\underset{%
i=l_{j}}{\overset{m_{j}}{\prod }}w_{2^{i}}\right) \left(
\sum_{k=l_{v}}^{m_{v}}\left( \underset{j=k+1}{\overset{m_{v}}{\prod }}%
w_{2^{j}}\right) \left( 2^{k}K_{2^{k}}-\left( 2^{k}-1\right)
D_{2^{k}}\right) \right) \\
&=&\sum_{v=1}^{r}\left( \underset{j=1}{\overset{v-1}{\prod }}\underset{i=l_{j}}%
{\overset{m_{j}}{\prod }}w_{2^{i}}\right) \left(
\sum_{k=0}^{m_{v}}-\sum_{k=0}^{l_{v}-1}\right) \left( \underset{j=k+1}{%
\overset{m_{v}}{\prod }}w_{2^{j}}\right) \left( 2^{k}K_{2^{k}}-\left(
2^{k}-1\right) D_{2^{k}}\right) \\
&=&\sum_{v=1}^{r}\left( \underset{j=1}{\overset{v-1}{\prod }}\underset{i=l_{j}}%
{\overset{m_{j}}{\prod }}w_{2^{i}}\right) \left( \sum_{k=0}^{m_{v}}\left(
\underset{j=k+1}{\overset{m_{v}}{\prod }}w_{2^{j}}\right) \left(
2^{k}K_{2^{k}}-\left( 2^{k}-1\right) D_{2^{k}}\right) \right) \\
&-&\sum_{v=1}^{r}\left( \underset{j=1}{\overset{v}{\prod }}\underset{i=l_{j}}{%
\overset{m_{j}}{\prod }}w_{2^{i}}\right) \left( \sum_{k=0}^{l_{v}-1}\left(
\underset{j=k+1}{\overset{l_{v}-1}{\prod }}w_{2^{j}}\right) \left(
2^{k}K_{2^{k}}-\left( 2^{k}-1\right) D_{2^{k}}\right) \right) .
\end{eqnarray*}

Since
$$2^{n}-1=\sum_{k=0}^{n-1}2^{k}$$
and
\begin{equation*}
\left( 2^{n}-1\right) K_{2^{n}-1}=\sum_{k=0}^{n-1}\left(
\prod_{j=k+1}^{n-1}w_{2^{j}}\right) \left( 2^{k}K_{2^{k}}-\left(
2^{k}-1\right) D_{2^{k}}\right),
\end{equation*}%
we obtain that
\begin{eqnarray*}
I_{1}&=&\sum_{v=1}^{r}\left( \underset{j=1}{\overset{v-1}{\prod }}\underset{%
i=l_{j}}{\overset{m_{j}}{\prod }}w_{2^{i}}\right) \left(
2^{m_{v}+1}-1\right) K_{2^{m_{v}+1}-1} \\
&-&\sum_{v=1}^{r}\left( \underset{j=1}{\overset{v}{\prod }}\underset{i=l_{j}}{%
\overset{m_{j}}{\prod }}w_{2^{i}}\right) \left( 2^{l_{v}}-1\right)
K_{2^{l_{v}}-1}.
\end{eqnarray*}

If we apply estimations
$$\left\vert K_{2^{n}}\right\vert \leq c\left\vert
K_{2^{n-1}}\right\vert $$
and
$$ \left\vert K_{2^{n}-1}\right\vert \leq c\left\vert
K_{2^{n}}\right\vert +c\  $$
we get that
\begin{equation}
\left\vert I_{1}\right\vert \leq c\sum_{v=1}^{r}\left( 2^{l_{v}}\left\vert
K_{2^{l_{v}}}\right\vert +2^{m_{v}}\left\vert K_{2^{m_{v}}}\right\vert
+cr\right) .  \label{1.12c}
\end{equation}

Let $l_{j}<n_{A}\leq m_{j},$ where $j=1,...,s.$ Then
$$n^{\left( A\right)
}\geq \sum_{v=l_{j}}^{n_{A}-1}2^{v}\geq 2^{n_{A}}-2^{l_{j}}$$ 
and
$$
2^{n_{A}}-1-n^{\left( A\right) }\leq 2^{l_{j}}.$$

If $l_{j}=n_{A}$ where $j=1,...,s,$ then
$$n^{\left( A\right) }\leq 2^{m_{j-1}+1}<2^{l_{j}}.$$

By using these estimations we can conclude that
\begin{equation}
\left\vert I_{2}\right\vert \leq
c\sum_{v=1}^{r}2^{l_{v}}\sum_{k=l_{v}}^{m_{v}}D_{2^{k}}.  \label{1.12d}
\end{equation}

By combining (\ref{1.12c})-(\ref{1.12d}) we get the proof of Lemma \ref{lemma6}.
\QED

The following estimations of Fej\'er kernels with respect to the one-dimensional Walsh-Fourier series is proved in \cite{tep13}:

\begin{lemma} \label{lemma7}
\label{lemma(Tephnadze)} Let $$n=\sum_{i=1}^{s}\sum_{k=l_{i}}^{m_{i}}2^{k},$$
where
$$0\leq l_{1}\leq m_{1}\leq l_{2}-2<l_{2}\leq m_{2}\leq ...\leq
l_{s}-2<l_{s}\leq m_{s}.$$

Then
\begin{equation*}
n\left\vert K_{n}\left( x\right) \right\vert \geq \frac{2^{2l_{i}}}{16},%
\text{ \ \ for \ \ }x\in I_{l_{i}+1}\left( e_{l_{i}-1}+e_{l_{i}}\right).
\end{equation*}
\end{lemma}

{\bf Proof}: 
If we apply Lemma \ref{lemma3} for $n=\sum_{i=1}^{s}\sum_{k=l_{i}}^{m_{i}}2^{k}$ we can write that
\begin{eqnarray*}
nK_{n} &=&\sum_{r=1}^{s}\sum_{k=l_{r}}^{m_{r}}\left( \underset{j=r+1}{%
\overset{s}{\prod }}\underset{q=l_{j}}{\overset{m_{j}}{\prod }}w_{2^{q}}%
\underset{j=k+1}{\overset{m_{r}}{\prod }}w_{2^{j}}\right) 2^{k}K_{2^{k}} \\
&+&\sum_{r=1}^{s}\sum_{k=l_{r}}^{m_{r}}\left( \underset{j=r+1}{\overset{s}{%
\prod }}\underset{q=l_{j}}{\overset{m_{j}}{\prod }}w_{2^{q}}\underset{j=k+1}{%
\overset{m_{r}}{\prod }}w_{2^{j}}\right) \left(
\sum_{t=1}^{r-1}\sum_{q=l_{t}}^{m_{t}}2^{q}+\sum_{q=l_{r}}^{k-1}2^{q}\right)
D_{2^{k}}.
\end{eqnarray*}

Let $x\in I_{l_{i}+1}\left( e_{l_{i}-1}+e_{l_{i}}\right) .$ Then

\begin{eqnarray*}
n\left\vert K_{n}\right\vert \geq\left\vert
2^{l_{i}}K_{2^{l_{i}}}\right\vert
-\sum_{r=1}^{i-1}\sum_{k=l_{r}}^{m_{r}}\left\vert 2^{k}K_{2^{k}}\right\vert
-\sum_{r=1}^{i-1}\sum_{k=l_{r}}^{m_{r}}\left\vert 2^{k}D_{2^{k}}\right\vert= I-II-III.
\end{eqnarray*}

Lemma \ref{lemma4} follows that

\begin{equation}
I=\left\vert 2^{l_{i}}K_{2^{l_{i}}}\left( x\right) \right\vert =\frac{%
2^{2l_{i}}}{4}.  \label{2.2.10.0}
\end{equation}

Since $m_{i-1}\leq l_{i}-2$, we easily obtain that the following estimation is true:
\begin{eqnarray} \label{2.2.10.1}
II&\leq &\sum_{n=0}^{l_{i}-2}\left\vert 2^{n}K_{2^{n}}\left( x\right)
\right\vert \\ \notag
&\leq & \sum_{n=0}^{l_{i}-2}2^{n}\frac{\left( 2^{n}+1\right) }{2} \\ \notag
&\leq & \frac{2^{2l_{i}}}{24}+\frac{2^{l_{i}}}{4}-\frac{2}{3}.
\end{eqnarray}

For $III$ we get that
\begin{equation} \label{2.2.10.2}
III\leq \sum_{k=0}^{l_{i}-2}\left\vert 2^{k}D_{2^{k}}\left( x\right)
\right\vert \leq \sum_{k=0}^{l_{i}-2}4^{k}=\frac{2^{2l_{i}}}{12}-\frac{1}{3}.
\end{equation}

By combining (\ref{2.2.10.0}-\ref{2.2.10.2}) we can conclude that
\begin{equation}
n\left\vert K_{n}\left( x\right) \right\vert \geq I-II-III\geq \frac{%
2^{2l_{i}}}{8}-\frac{2^{l_{i}}}{4}+1.  \label{10.3}
\end{equation}

Suppose that $l_{i}\geq 2$. Then
\begin{equation*}
n\left\vert K_{n}\left( x\right) \right\vert \geq \frac{2^{2l_{i}}}{8}-\frac{%
2^{2l_{i}}}{16}\geq \frac{2^{2l_{i}}}{16}.
\end{equation*}

If $l_{i}=0$ or $l_{i}=1$, then by applying (\ref{10.3}) we get that
\begin{equation*}
n\left\vert K_{n}\left( x\right) \right\vert \geq \frac{7}{8}\geq \frac{%
2^{2l_{i}}}{16},
\end{equation*}
Lemma is proved.
\QED

The following estimations of Fej\'er kernels with respect to the one-dimensional Walsh-Fourier series is proved in \cite{tep13} (see also \cite{tep_thesis}):

\begin{lemma} \label{lemma3.2.10}
Let $0<p\leq1$, $2^{k}\leq n<2^{k+1}$ and $\sigma_{n}f$ be Fej\'er means with respect to the one-dimensional Walsh-Fourier series, where $f\in H_{p}(G)$. Then, for any fixed $n\in \mathbb{N}$,
\begin{eqnarray*}
&&\left\Vert \sigma_{n}f\right\Vert _{H_{p}(G)} \\
&&\leq \left\Vert \sup_{0\leq l \leq k}\left\vert \sigma_{2^{l}}f\right\vert \right\Vert _{p}+\left\Vert \sup_{0\leq l\leq k}\left\vert S_{2^{l}}f\right\vert \right\Vert _{p}+\left\Vert \sigma_{n}f\right\Vert _{p} \\
&&\leq \left\Vert \widetilde{\sigma}_{\#}^{\ast
}f\right\Vert _{p}+\left\Vert \widetilde{S}_{\#}^{\ast
}f\right\Vert _{p}+\left\Vert \sigma_{n}f\right\Vert _{p}.
\end{eqnarray*}
\end{lemma}

{\bf Proof}:
Let consider the following martingale
\begin{eqnarray*}
f_{\#}&=&\left( S_{2^{k}}\sigma_{n}f,\text{ }k \in \mathbb{N}\right) \\
&=&\left( \frac{2^{0}\sigma_{2^{0}}}{n}+\frac{(n-2^{0})S_{2^{0}}f}{n},...,\frac{2^{k}\sigma_{2^{k}}f}{n}+\frac{(n-2^{k})S_{2^{k}}f}{n},
\sigma_{n}f,...,\sigma_{n}f,...\right).
\end{eqnarray*}
By using Lemma \ref{lemma3.2.3} we immediately get
\begin{eqnarray*}
&&\left\Vert \sigma_{n}f\right\Vert _{H_{p}(G^2)}^{p} \\
&&\leq \left\Vert \sup_{0\leq l\leq k}\left\vert \sigma_{2^{l}}f\right\vert \right\Vert _{p}^{p}+\left\Vert \sup_{0\leq l\leq k}\left\vert S_{2^{l}}f\right\vert \right\Vert _{p}^{p}+\left\Vert
S_{n}f\right\Vert _{p}^{p} \\
&&\leq \left\Vert \widetilde{\sigma}_{\#}^{\ast }f\right\Vert
_{p}^{p}+\left\Vert \widetilde{S}_{\#}^{\ast }f\right\Vert
_{p}^{p}+\left\Vert \sigma_{n}f\right\Vert _{p}^{p}.
\end{eqnarray*}

Lemma is proved.
\QED

\subsection{Boundedness of subsequences of Fej\'er means with respect to the one-dimensional Walsh-Fourier series on the martingale Hardy spaces}

\text{ \qquad } In this section we study boundedness of subsequences of Fej\'er means with respect to the one-dimensional Walsh-Fourier series in the martingale Hardy spaces (For details see \cite{tep13}).

First, we consider case $ p=1/2 $. The following estimation is true:
\begin{theorem} \label{th5.1.1}
\label{theorem1}a) Let $f\in H_{1/2}(G).$ Then there exists an absolute constant $c,$ such that
\begin{equation*}
\left\Vert \sigma _{n}f\right\Vert _{H_{1/2}(G)}\leq cV^{2}\left( n\right)
\left\Vert f\right\Vert _{H_{1/2}(G)}.
\end{equation*}

b) Let $\left\{ n_{k}:k\in\mathbb{N}_{+}\right\} $  be increasing sequence of natural numbers, such that $\sup_{k \in \mathbb{N}_{+}}V\left(
n_{k}\right) =\infty $ and $\Phi :\mathbb{N}_{+}\rightarrow \lbrack
1,\infty )$ be non-decreasing function satisfying the conditions
$\Phi \left( n\right) \uparrow \infty $
and
\begin{equation}
\overline{\underset{k\rightarrow \infty }{\lim }}\frac{V^{2}\left(
n_{k}\right) }{\Phi \left( n_{k}\right) }=\infty .  \label{5.1.30}
\end{equation}%
Then there exists a martingale $f\in H_{1/2}(G),$ such that
\begin{equation*}
\underset{k\in \mathbb{N}}{\sup }\left\Vert \frac{\sigma _{n_{k}}f}{\Phi
\left(n_{k}\right) }\right\Vert _{1/2}=\infty .
\end{equation*}
\end{theorem}

{\bf Proof}:
Suppose that
\begin{equation} \label{5.1.12k}
\left\Vert \frac{\sigma _{n}f}{V^{2}\left( n\right) }\right\Vert _{1/2}\leq
c\left\Vert f\right\Vert_{H_{1/2}(G)}.
\end{equation}%

By combining estimations (\ref{1.S2n}), (\ref{sigmamax}) and Lemma \ref{lemma3.2.10} we can conclude that
\begin{eqnarray} \label{5.1.12l}
&&\left\Vert \frac{\sigma _{n}f}{V^{2}\left( n\right) }\right\Vert
_{H_{1/2}(G)}^{1/2} \\ \notag
&\leq &\left\Vert \frac{\sigma _{n}f}{V^{2}\left( n\right) }\right\Vert_{1/2}^{1/2}+\frac{1}{V^{2}\left( n\right) }\left\Vert \sigma_{\#}^{\ast }f \right\Vert_{1/2}^{1/2} +\frac{1}{V^{2}\left( n\right) }\left\Vert\widetilde{S}_{\#}^{\ast }\right\Vert_{1/2}^{1/2}\\ \notag
&\leq & \left\Vert \frac{\widetilde{\sigma} _{n}f}{V^{2}\left( n\right) }\right\Vert_{1/2}^{1/2}+\left\Vert \widetilde{\sigma}_{\#}^{\ast }f \right\Vert_{1/2}^{1/2}+\left\Vert \widetilde{S}_{\#}^{\ast }f \right\Vert_{1/2}^{1/2} \leq  c\left\Vert f\right\Vert _{H_{1/2}(G)}^{1/2}.
\end{eqnarray}

By combining Lemma \ref{lemma3.2.5} and (\ref{5.1.12l}), Theorem \ref{theorem1} will be proved if we show that
\begin{equation*}
\int_{\overline{I_{M}}}\left( \frac{\left\vert \sigma _{n}a\right\vert }{V^{2}\left( n\right) }\right) ^{1/2}d\mu  \leq c<\infty ,
\end{equation*}%
for any $ 1/2 $-atom $a$.

Without loss the generality we may assume that $a$ is $1/2$-atom, with support $I,$ for which $\mu \left( I\right) =2^{-M},$  $I=I_{M}.$ It is easy to check that $\sigma _{n}\left( a\right) =0,$ when $n\leq 2^{M}.$ Therefore, we may assume that $n>2^{M}.$ Set
\begin{eqnarray*}
&& II_{\alpha _{A}}^{1}\left( x\right):=2^{M}\int_{I_{M}}2^{\alpha
_{A}}\left\vert K_{2^{\alpha _{A}}}\left( x+t\right) \right\vert d\mu \left(t\right) ,\text{\ } \\
&& II_{l_{A}}^{2}\left( x\right)
=2^{M}\int_{I_{M}}2^{l_{A}}\sum_{k=l_{A}}^{m_{A}}D_{2^{k}}\left(
x+t\right) d\mu \left( t\right).
\end{eqnarray*}

Let $x\in I_{M}.$ Since $\sigma _{n}$ is bounded from $L_{\infty }(G)$ to $L_{\infty }(G)$, for $n>2^{M}$ and $\left\Vert a\right\Vert _{\infty }\leq 2^{2M},$ by using Lemma \ref{lemma5} we can conclude that
\begin{eqnarray*}
&&\frac{\left\vert \sigma _{n}a\left( x\right) \right\vert }{V^{2}\left(n\right)} \\
&& \leq \frac{c}{V^{2}\left( n\right) }\int_{I_{M}}\left\vert a\left(
x\right) \right\vert \left\vert K_{n}\left( x+t\right) \right\vert d\mu\left( t\right) \\
&& \leq \frac{c\left\Vert a\right\Vert _{\infty }}{V^{2}\left(n\right) }\int_{I_{M}}\left\vert K_{n}\left( x+t\right) \right\vert d\mu \left(t\right) \\
&& \leq \frac{c2^{2M}}{V^{2}\left( n\right) }\int_{I_{M}}\left\vert
K_{n}\left( x+t\right) \right\vert d\mu \left( t\right)\\
&&\leq \frac{c2^{M}}{V^{2}\left( n\right) }\left\{
\sum_{A=1}^{s}\int_{I_{M}}2^{l_{A}}\left\vert K_{2^{l_{A}}}\left( x+t\right)
\right\vert d\mu \left( t\right) +\int_{I_{M}}2^{m_{A}}\left\vert
K_{2^{m_{A}}}\left( x+t\right) \right\vert d\mu \left( t\right) \right\}\\
&&+\frac{c2^{M}}{V^{2}\left( n\right) }\sum_{A=1}^{s}\int_{I_{M}}2^{l_{A}}%
\sum_{k=l_{A}}^{m_{A}}D_{2^{k}}\left( x+t\right) d\mu \left( t\right) +\frac{c2^{M}}{V^{2}\left( n\right) }\int_{I_{M}}V\left( n\right) d\mu \left(t\right)\\
&&=\frac{c}{V^{2}\left( n\right) }\sum_{A=1}^{s}\left( II_{^{l_{A}}}^{1}\left(
x\right) +II_{^{m_{A}}}^{1}\left( x\right) +II_{l_{A}}^{2}\left( x\right)
\right) +c.
\end{eqnarray*}

Hence,
\begin{eqnarray*} \\
&&\int_{\overline{I_{M}}}\left\vert \frac{\sigma _{n}a\left( x\right) }{%
V^{2}\left( n\right) }\right\vert ^{1/2}d\mu \left( x\right)
\\
&\leq &\frac{c}{V\left( n\right) }\left( \sum_{A=1}^{s}\int_{\overline{I_{M}}%
}\left\vert II_{l_{A}}^{1}\left( x\right) \right\vert ^{1/2}d\mu \left(
x\right) \right. \\
&&\left. +\int_{\overline{I_{M}}}\left\vert II_{m_{A}}^{1}\left( x\right)
\right\vert ^{1/2}d\mu \left( x\right) +\int_{\overline{I_{M}}}\left\vert
II_{l_{A}}^{2}\left( x\right) \right\vert ^{1/2}d\mu \left( x\right) \right)
+c.
\end{eqnarray*}

Since  $s\leq 4V\left( n\right) $ we obtain that Theorem \ref{th5.1.1} will be proved if we show that
\begin{equation} \label{5.1.11.1}
\int_{\overline{I_{M}}}\left\vert II_{\alpha _{A}}^{1}\left( x\right)
\right\vert ^{1/2}d\mu \left( x\right) \leq c<\infty ,\text{\ }\int_{%
\overline{I_{M}}}\left\vert II_{l_{A}}^{2}\left( x\right) \right\vert
^{1/2}d\mu \left( x\right) \leq c<\infty ,\text{\ }
\end{equation}%
where $\alpha _{A}=l_{A}$ or $\alpha _{A}=m_{A},$ $A=1,...,s$.

Let $t\in I_{M}$ and $x\in I_{l+1}\left( e_{k}+e_{l}\right) ,$ $0\leq
k<l<\alpha _{A}\leq M$ or $0\leq k<l\leq M\leq \alpha _{A}.$ Since $x+t\in
I_{l+1}\left( e_{k}+e_{l}\right),$ by applying Lemma \ref{lemma4} we can conclude that
\begin{equation} \label{5.1.10a}
K_{2^{\alpha _{A}}}\left( x+t\right) =0\text{ \ \ and \ \ }II_{\alpha
_{A}}^{1}\left( x\right) =0.
\end{equation}

Let $x\in I_{l+1}\left( e_{k}+e_{l}\right),$ $0\leq k<\alpha _{A}\leq l\leq M.$ Then $x+t\in I_{l+1}\left( e_{k}+e_{l}\right),$ where $t\in I_{M}$ and if we apply again Lemma \ref{lemma4} we get that
\begin{equation} \label{5.1.10b}
2^{\alpha _{A}}\left\vert K_{2^{\alpha _{A}}}\left( x+t\right) \right\vert\leq 2^{\alpha _{A}+k}\text{ \ \ and \ \ }II_{\alpha _{A}}^{1}\left(x\right) \leq 2^{\alpha_{A}+k}.
\end{equation}

Analogously to (\ref{5.1.10b}) for $0\leq \alpha _{A}\leq k<l\leq M$ we can prove that
\begin{equation} \label{5.1.10c}
2^{\alpha_{A}}\left\vert K_{2^{\alpha _{A}}}\left( x+t\right) \right\vert\leq 2^{2\alpha_{A}},\text{ \ }II_{\alpha _{A}}^{1}\left( x\right)\leq2^{2\alpha_{A}},\text{\ }t\in I_{M},\text{\ }x\in I_{l+1}\left(
e_{k}+e_{l}\right).
\end{equation}

Let $0\leq \alpha _{A}\leq M-1,$ where $A=1,...,s.$ According to (\ref{2.1.2}) and (\ref{5.1.10a}-\ref{5.1.10c}) we find that
\begin{eqnarray*}
&&\int_{\overline{I_{M}}}\left\vert II_{\alpha _{A}}^{1}\left( x\right)
\right\vert ^{1/2}d\mu \left( x\right) \\
&=&\overset{M-2}{\underset{k=0}{\sum }}\overset{M-1}{\underset{l=k+1}{\sum }}%
\int_{I_{l+1}\left( e_{k}+e_{l}\right) }\left\vert II_{\alpha
_{A}}^{1}\left( x\right) \right\vert ^{1/2}d\mu \left( x\right)\\
&+&\overset{M-1%
}{\underset{k=0}{\sum }}\int_{I_{M}\left( e_{k}\right) }\left\vert
II_{\alpha _{A}}^{1}\left( x\right) \right\vert ^{1/2}d\mu \left( x\right) \\
&\leq & c\overset{\alpha _{A}-1}{\underset{k=0}{\sum }}\overset{M-1}{\underset{%
l=\alpha _{A}+1}{\sum }}\int_{I_{l+1}\left( e_{k}+e_{l}\right) }2^{\left(
\alpha _{A}+k\right) /2}d\mu \left( x\right) \\
& +& c\overset{M-2}{\underset{%
k=\alpha _{A}}{\sum }}\overset{M-1}{\underset{l=k+1}{\sum }}%
\int_{I_{l+1}\left( e_{k}+e_{l}\right) }2^{\alpha _{A}}d\mu \left( x\right) \\
&+& c\overset{\alpha _{A}-1}{\underset{k=0}{\sum }}\int_{I_{M}\left(
e_{k}\right) }2^{\left( \alpha _{A}+k\right) /2}d\mu \left( x\right)
+c\overset{M-1}{\underset{k=\alpha _{A}}{\sum }}\int_{I_{M}\left( e_{k}\right)
}2^{\alpha _{A}}d\mu \left( x\right) \\
&\leq & c\overset{\alpha _{A}-1}{\underset{k=0}{\sum }}\overset{M-1}{\underset{%
l=\alpha _{A}+1}{\sum }}\frac{2^{\left( \alpha _{A}+k\right) /2}}{2^{l}}+c%
\overset{M-2}{\underset{k=\alpha _{A}}{\sum }}\overset{M-1}{\underset{l=k+1}{%
\sum }}\frac{2^{\alpha _{A}}}{2^{l}} \\
&+& c\overset{\alpha _{A}-1}{\underset{k=0}{%
\sum }}\frac{2^{\left( \alpha _{A}+k\right) /2}}{2^{M}}+c\overset{M-1}{%
\underset{k=\alpha _{A}}{\sum }}\frac{2^{\alpha _{A}}}{2^{M}}
\leq c<\infty.
\end{eqnarray*}

Let $\alpha _{A}\geq M.$ Analogously to $II_{\alpha _{A}}^{1}\left( x\right)$ we can prove (\ref{5.1.11.1}), for $A=1,...,s.$

Now, prove boundedness of $II_{l_{A}}^{2}$. Let $t\in I_{M}$ and $x\in I_{i}\backslash I_{i+1},$ $i\leq l_{A}-1.$ Since $x+t\in I_{i}\backslash
I_{i+1},$ if we apply first equality of Lemma \ref{lemma1} we get that
\begin{equation} \label{5.1.13a}
II_{l_{A}}^{2}\left( x\right) =0.
\end{equation}

Let $x\in I_{i}\backslash I_{i+1},$ $l_{A}\leq i\leq m_{A}.$ Since  $n\geq
2^{M}$ and $t\in I_{M},$ if we apply first equality of Lemma \ref{lemma1} we get that
\begin{equation} \label{5.1.13b}
II_{l_{A}}^{2}\left( x\right) \leq
2^{M}\int_{I_{M}}2^{l_{A}}\sum_{k=l_{A}}^{i}D_{2^{k}}\left( x+t\right) d\mu\left( t\right) \leq c2^{l_{A}+i}.
\end{equation}

Let $x\in I_{i}\backslash I_{i+1},$ $m_{A}<i\leq M-1.$ Then  $x+t\in I_{i}\backslash I_{i+1},$ for any $t\in I_{M}$ and by first equality of Lemma \ref{lemma1} we have that
\begin{equation} \label{5.1.13c}
II_{l_{A}}^{2}\left( x\right) \leq c2^{M}\int_{I_{M}}2^{l_{A}+m_{A}}\leq
c2^{l_{A}+m_{A}}.
\end{equation}

Let  $0\leq l_{A}\leq m_{A}\leq M.$ Then, in the view of (\ref{2.1.2})  and (\ref{5.1.13a}-\ref{5.1.13c}) we can conclude that
\begin{eqnarray*}
&&\int_{\overline{I_{M}}}\left\vert II_{l_{A}}^{2}\left( x\right) \right\vert
^{1/2}d\mu \left( x\right) \\
&=&\left(
\sum_{i=0}^{l_{A}-1}+\sum_{i=l_{A}}^{m_{A}}+\sum_{i=m_{A}+1}^{M-1}\right)
\int_{I_{i}\backslash I_{i+1}}\left\vert II_{l_{A}}^{2}\left( x\right)
\right\vert ^{1/2}d\mu \left( x\right) \\
&\leq& c\sum_{i=l_{A}}^{m_{A}}\int_{I_{i}\backslash I_{i+1}}2^{\left(
l_{A}+i\right) /2}d\mu \left( x\right) \\
&+& c\sum_{i=m_{A}+1}^{M-1}\int_{I_{i}\backslash I_{i+1}}2^{\left(
l_{A}+m_{A}\right) /2}d\mu \left( x\right) \\
&\leq & c\sum_{i=l_{A}}^{m_{A}}2^{\left( l_{A}+i\right) /2}\frac{1}{2^{i}} \\
&+& c\sum_{i=m_{A}+1}^{M-1}2^{\left(l_{A}+m_{A}\right) /2}\frac{1}{2^{i}}\leq c<\infty .
\end{eqnarray*}

Analogously, we can prove same estimations in the cases
$0\leq l_{A}\leq M<m_{A}$ and $M\leq l_{A}\leq m_{A}.$

Now, we prove part b) of Theorem \ref{th5.1.1}. According to (\ref{5.1.30}), there exists increasing sequence $\left\{ \alpha _{k}:\text{ }k\in\mathbb{N}_{+}\right\} \subset \left\{ n_{k}\text{ }:k\in\mathbb{N}_{+}\right\} $ of natural numbers such that
\begin{equation} \label{5.1.2aaa}
\sum_{k=1}^{\infty }\frac{\Phi ^{1/4}\left( \alpha _{k}\right)}{V^{1/2}\left(
\alpha _{k}\right) }\leq c<\infty.
\end{equation}

Let $f=\left(f_{n},n\in\mathbb{N}_{+}\right)$ be martingale form Example \ref{example2.2.1}, where
\begin{equation*}
\lambda _{k}:=\Phi ^{1/2}\left( \alpha _{k}\right) /V\left( \alpha
_{k}\right).
\end{equation*}

According to (\ref{5.1.2aaa}) we get that condition (\ref{3.3.2aa}) is fulfilled and it follows that $f=\left(f_{n},n\in\mathbb{N}_{+}\right)$.

By applying (\ref{3.3.10AA}) we get that
\begin{equation} \label{5.1.5aa}
\widehat{f}(j)
\end{equation}
\begin{equation*}
=\left\{
\begin{array}{ll}
2^{\left\vert \alpha _{k}\right\vert }\Phi ^{1/2}\left( \alpha _{k}\right)/V\left( \alpha _{k}\right) , & \text{\thinspace \thinspace }j\in \left\{2^{\left\vert \alpha _{k}\right\vert },...,2^{_{\left\vert \alpha_{k}\right\vert +1}}-1\right\} ,\text{ }k \in\mathbb{N}_{+} \\
0\,, & \text{\thinspace }j\notin \bigcup\limits_{k=0}^{\infty }\left\{2^{_{\left\vert \alpha _{k}\right\vert }},..., 2^{_{\left\vert \alpha_{k}\right\vert +1}}-1\right\}.
\end{array}
\right.
\end{equation*}

Let $2^{\left\vert \alpha _{k}\right\vert }<j<\alpha _{k}.$ If we apply (\ref{3.3.11AA}) we get that
\begin{equation} \label{5.1.sn}
S_{j}f=S_{2^{\left\vert \alpha_{k}\right\vert }}f+\frac{w_{2^{\left\vert \alpha_{k}\right\vert }}D_{j-2^{\left\vert \alpha _{k}\right\vert }}\Phi^{1/2}\left( \alpha _{k}\right) }{V\left( \alpha _{k}\right)}
\end{equation}

Hence,
\begin{eqnarray} \label{5.1.7aaa}
&&\frac{\sigma _{_{\alpha _{k}}}f}{\Phi \left( \alpha _{k}\right) } \\ \notag
&=&\frac{1}{\Phi \left( \alpha _{k}\right) \alpha_{k}}\sum_{j=1}^{2^{\left\vert \alpha_{k}\right\vert }}S_{j}f \\ \notag
&+&\frac{1}{\Phi \left( \alpha _{k}\right) \alpha _{k}}\sum_{j=2^{\left\vert \alpha _{k}\right\vert }+1}^{\alpha_{k}}S_{j}f
\\ \notag
&=&\frac{\sigma _{_{2^{\left\vert \alpha _{k}\right\vert }}}f}{\Phi \left(\alpha _{k}\right) \alpha _{k}} \\ \notag
&+&\frac{\left( \alpha _{k}-2^{\left\vert\alpha_{k}\right\vert }\right) S_{2^{\left\vert \alpha _{k}\right\vert }}f}{\Phi \left(\alpha _{k}\right) \alpha _{k}} \\ \notag
&+&\frac{w_{2^{\left\vert \alpha_{k}\right\vert }}2^{\left\vert \alpha_{k}\right\vert }\Phi^{1/2}\left(\alpha _{k}\right)}{\Phi \left( \alpha_{k}\right) V\left( \alpha_{k}\right) \alpha _{k}}\sum_{j=2^{_{\left\vert\alpha_{k}\right\vert }}+1}^{\alpha _{k}}D_{j-2^{\left\vert \alpha_{k}\right\vert }} \\ \notag
&=& III_{1}+III_{2}+III_{3}.
\end{eqnarray}%

For $ III_{3} $ we can conclude that
\begin{eqnarray} \label{5.1.9aaa}
&&\left\vert III_{3}\right\vert \\ \notag
&=&\frac{2^{\left\vert \alpha _{k}\right\vert }\Phi ^{1/2}\left( \alpha
_{k}\right) }{\Phi \left( \alpha _{k}\right) V\left( \alpha _{k}\right)
\alpha _{k}}\left\vert \sum_{j=1}^{\alpha _{k}-2^{_{\left\vert \alpha
_{k}\right\vert }}}D_{j}\right\vert \\ \notag
&=&\frac{2^{\left\vert \alpha
_{k}\right\vert }\Phi ^{1/2}\left( \alpha _{k}\right) }{\Phi \left( \alpha
_{k}\right) V\left( \alpha _{k}\right) \alpha _{k}}\left( \alpha
_{k}-2^{_{\left\vert \alpha _{k}\right\vert }}\right) \left\vert K_{\alpha
_{k}-2^{_{\left\vert \alpha _{k}\right\vert }}}\right\vert \\ \notag
&\geq& \frac{c\left( \alpha _{k}-2^{_{\left\vert \alpha _{k}\right\vert }}\right)\left\vert K_{\alpha _{k}-2^{_{\left\vert \alpha _{k}\right\vert
}}}\right\vert }{\Phi ^{1/2}\left( \alpha _{k}\right) V\left( \alpha
_{k}\right).}
\end{eqnarray}

Let
\begin{equation*}
\alpha_{k}=\sum_{i=1}^{r_{k}}\sum_{k=l_{i}^{k}}^{m_{i}^{k}}2^{k},
\end{equation*}
where
\begin{equation*}
m_{1}^{k}\geq l_{1}^{k}>l_{1}^{k}-2\geq m_{2}^{k}\geq
l_{2}^{k}>l_{2}^{k}-2\geq ...\geq m_{s}^{k}\geq l_{s}^{k}\geq 0.
\end{equation*}

Since (see theorems \ref{th4.1.1} and \ref{th5.1.1})
\begin{equation*}
\left\Vert III_{1}\right\Vert
_{1/2}\leq c, \left\Vert III_{2}\right\Vert _{1/2}\leq c,
\end{equation*}
and
\begin{equation*}
\mu \left\{E_{l_{i}^{k}}\right\} \geq 1/2^{l_{i}^{k}-1},
\end{equation*}
By combining (\ref{5.1.7aaa}), (\ref{5.1.9aaa}) and Lemma \ref{lemma7} we get that
\begin{eqnarray*}
&&\int_{G}\left\vert \sigma _{\alpha _{k}}f(x)/\Phi \left( \alpha _{k}\right)
\right\vert ^{1/2}d\mu \left( x\right) \\
&\geq & \left\Vert III_{3}\right\Vert
_{1/2}^{1/2}-\left\Vert III_{2}\right\Vert _{1/2}^{1/2}-\left\Vert
III_{1}\right\Vert _{1/2}^{1/2} \\
&\geq & c\text{ }\underset{i=2}{\overset{r_{k}-2}{\sum }}\int_{E_{l_{i}^{k}}}%
\left\vert 2^{2l_{i}^{k}}/\left( \Phi ^{1/2}\left( \alpha _{k}\right)
V\left( \alpha _{k}\right) \right) \right\vert ^{1/2}d\mu \left( x\right) -2c \\
&\geq & c\overset{r_{k}-2}{\underset{i=2}{\sum }}1/\left( V^{1/2}\left( \alpha
_{k}\right) \Phi ^{1/4}\left( \alpha _{k}\right) \right)-2c \\
&\geq&
cr_{k}/\left( V^{1/2}\left( \alpha _{k}\right) \Phi ^{1/4}\left( \alpha
_{k}\right) \right) \\
&\geq & cV^{1/2}\left( \alpha _{k}\right) /\Phi ^{1/4}\left( \alpha _{k}\right)
\rightarrow \infty ,\text{ as }k\rightarrow \infty .
\end{eqnarray*}

Theorem \ref{th5.1.1} is proved.
\QED

\begin{theorem} \label{th5.1.2}
a) Let $0<p<1/2,$ $f\in H_{p}(G).$ Then there exists an absolute constant $c_{p}$ depending only on $p$ such that
\begin{equation*}
\text{ }\left\Vert \sigma _{n}f\right\Vert _{H_{p}(G)}\leq c_{p}2^{d\left(
n\right) \left( 1/p-2\right) }\left\Vert f\right\Vert _{H_{p}(G)}.
\end{equation*}

b) Let $0<p<1/2$ and $\Phi \left( n\right):\mathbb{N}_{+}\rightarrow \lbrack
1,\infty ) $ be non-decreasing function such that
\begin{equation}  \label{5.1.31aaa}
\sup_{k\in \mathbb{N}_+}d\left( n_{k}\right) =\infty ,\text{ \ \ }\overline{\underset{k\rightarrow \infty }{\lim }}\frac{2^{d\left( n_{k}\right) \left(1/p-2\right) }}{\Phi \left( n_{k}\right) }=\infty .
\end{equation}%
Then there exist a martingale $f\in H_{p}(G),$ such that
\begin{equation*}
\underset{k\in\mathbb{N}_+}{\sup}\left\Vert\frac{\sigma_{n_{k}}f}{\Phi \left(n_{k}\right)}\right\Vert_{weak-L_{p}(G)}=\infty.
\end{equation*}
\end{theorem}
{\bf Proof}:
Let $n\in \mathbb{N}.$ Analogously to (\ref{5.1.12l}) it is sufficient to prove that
\begin{equation*}
\int\limits_{\overline{I}_{M}}\left( 2^{d\left( n\right) \left( 2-1/p\right)
}\left\vert \sigma _{n}\left( a\right) \right\vert \right) ^{p}d\mu \leq
c_{p}<\infty ,
\end{equation*}%
for every $ p $-atom $a$, where $I$ denotes support of the atom.

Analogously to Theorem \ref{th5.1.1} we may assume that $a$
is $ p $-atom with support $\ I=I_{M}$, $\mu \left( I_{M}\right)=2^{-M}$ and  $n>2^{M}.$
Since $\left\Vert a\right\Vert _{\infty }\leq
2^{M/p} $ we can conclude that
\begin{eqnarray*}
&& 2^{d\left( n\right) \left( 2-1/p\right) }\left\vert \sigma _{n}a\right\vert \\
&\leq & 2^{d\left( n\right) \left( 2-1/p\right)}\left\Vert a\right\Vert _{\infty }\int_{I_{M}}\left\vert K_{n}\left(
x+t\right) \right\vert d\mu \left( t\right) \\
&\leq & 2^{d\left( n\right) \left( 2-1/p\right) }2^{M/p}\int_{I_{M}}\left\vert
K_{n}\left( x+t\right) \right\vert d\mu \left( t\right).
\end{eqnarray*}

Let $x\in I_{l+1}\left( e_{k}+e_{l}\right) ,\,0\leq k,l\leq \left[ n\right]
\leq M.$ Then, by applying Lemma \ref{lemma4} we get that $K_{n}\left( x+t\right)=0, $ where $t\in I_{M}$ and hence,
\begin{equation} \label{5.1.12a}
2^{d\left( n\right) \left( 2-1/p\right) }\left\vert \sigma _{n}a \right\vert =0.
\end{equation}

Let $x\in I_{l+1}\left( e_{k}+e_{l}\right) ,\,\left[ n\right] \leq k,l\leq M$
or  $k\leq \left[ n\right] \leq l\leq M.$ Then Lemma \ref{lemma6} follows that
\begin{eqnarray} \label{5.1.12}
2^{d\left( n\right) \left( 2-1/p\right) }\left\vert \sigma _{n}a \right\vert &\leq &  2^{d\left( n\right) \left( 2-1/p\right)}2^{M\left( 1/p-2\right) +k+l}  \\ \notag
&\leq & c_{p}2^{\left[ n\right] \left( 1/p-2\right) +k+l}.
\end{eqnarray}

By combining (\ref{2.1.2}), (\ref{5.1.12a}) and (\ref{5.1.12}) we can conclude that
\begin{eqnarray*}
&&\int_{\overline{I_{M}}}\left\vert 2^{d\left( n\right) \left( 2-1/p\right)
}\sigma _{n}a\left( x\right) \right\vert ^{p}d\mu \left( x\right) \\
&\leq& \left( \overset{\left[ n\right] -2}{\underset{k=0}{\sum }}\overset{%
\left[ n\right] -1}{\underset{l=k+1}{\sum }}+\overset{\left[ n\right] -1}{%
\underset{k=0}{\sum }}\overset{M-1}{\underset{l=\left[ n\right] }{\sum }}+%
\overset{M-2}{\underset{k=\left[ n\right] }{\sum }}\overset{M-1}{\underset{%
l=k+1}{\sum }}\right) \int_{I_{l+1}\left( e_{k}+e_{l}\right) }\left\vert
2^{d\left( n\right) \left( 2-1/p\right) }\sigma _{n}a\left( x\right)
\right\vert ^{p}d\mu \left( x\right) \\
&+&\overset{M-1}{\underset{k=0}{\sum }}\int_{I_{M}\left( e_{k}\right)
}\left\vert 2^{d\left( n\right) \left( 2-1/p\right) }\sigma _{n}a\left(
x\right) \right\vert ^{p}d\mu \left( x\right) \\
&\leq & c_{p}\overset{M-2}{%
\underset{k=\left[ n\right] }{\sum }}\overset{M-1}{\underset{l=k+1}{\sum }}%
\frac{1}{2^{l}}2^{\left[ n\right] \left( 2p-1\right) }2^{p\left( k+l\right) } \\
&+& c_{p}\overset{\left[ n\right] }{\underset{k=0}{\sum }}\overset{M-1}{%
\underset{l=\left[ n\right] +1}{\sum }}\frac{1}{2^{l}}2^{\left[ n\right]
\left( 2p-1\right) }2^{p\left( k+l\right) } \\
&+&\frac{c_{p}2^{\left[ n\right]
\left( 2p-1\right) }}{2^{M}}\overset{\left[ n\right] }{\underset{k=0}{\sum }}%
2^{p\left( k+M\right) }<c_{p}<\infty .
\end{eqnarray*}

Now, we prove part b) of Theorem \ref{th5.1.2}. According to   (\ref{5.1.31aaa}) there exists an increasing sequence of natural numbers $\left\{ \alpha _{k}:\text{ }%
k\in \mathbb{N}_+\right\} \subset \left\{n_{k}:\text{ }k\in \mathbb{N}_+\right\},$ such that
$\alpha _{0}\geq 3$ and
\begin{equation} \label{5.1.121}
\sum_{\eta =0}^{\infty }u^{-p}\left( \alpha _{\eta }\right) <c_{p}<\infty, \text{ \ \ }
\text{\ }u\left( \alpha _{k}\right) =2^{d\left( \alpha _{k}\right) \left(
1/p-2\right) /2}/\Phi ^{1/2}\left( \alpha _{k}\right) .
\end{equation}

Let $f$ be martingale from Example \ref{example2.2.1}, where
\begin{equation*}
\lambda_{k}=u^{-1}\left( \alpha _{k}\right),
\end{equation*}

If we apply  (\ref{5.1.121}) we get that (\ref{3.3.2aa}) is fulfilled and it follows that $f\in H_{p}(G).$ According to (\ref{3.3.10AA}) we have that
\begin{equation} \label{5.1.6aa}
\widehat{f}(j)
=\left\{
\begin{array}{ll}
2^{\left\vert \alpha _{k}\right\vert \left( 1/p-1\right) }/u\left( \alpha
_{k}\right) , & j\in \left\{ 2^{\left\vert \alpha _{k}\right\vert
},...,2^{_{\left\vert \alpha _{k}\right\vert +1}}-1\right\} ,\text{ }%
k\in\mathbb{N}_+, \\
0\,, & j\notin \bigcup\limits_{k=0}^{\infty }\left\{ 2^{\left\vert \alpha
_{k}\right\vert },...,2^{_{\left\vert \alpha _{k}\right\vert +1}}-1\right\} .
\end{array}
\right.
\end{equation}

Let $2^{\left\vert \alpha _{k}\right\vert }<j<\alpha _{k}.$ Then, analogously to (\ref{5.1.sn}) and (\ref{5.1.7aaa}), if we apply (\ref{5.1.6aa}) we get that
\begin{eqnarray*}
&&\frac{\sigma _{_{\alpha _{k}}}f}{\Phi \left( \alpha _{k}\right)} \\ &=&\frac{%
\sigma _{_{2^{\left\vert \alpha _{k}\right\vert }}}f}{\Phi \left( \alpha
_{k}\right) \alpha _{k}}+\frac{\left( \alpha _{k}-2^{\left\vert \alpha
_{k}\right\vert }\right) S_{2^{\left\vert \alpha _{k}\right\vert }}f}{\Phi
\left( \alpha _{k}\right) \alpha _{k}} \\
&+&\frac{2^{\left\vert \alpha _{k}\right\vert \left( 1/p-1\right) }}{\Phi
\left( \alpha _{k}\right) u\left( \alpha _{k}\right) \alpha _{k}}%
\sum_{j=2^{_{\left\vert \alpha _{k}\right\vert }}}^{\alpha _{k}-1}\left(
D_{_{j}}-D_{2^{\left\vert \alpha _{k}\right\vert }}\right) \\
&=& IV_{1}+IV_{2}+IV_{3}.
\end{eqnarray*}

Let $\alpha _{k}\in \mathbb{N}$ and $E_{\left[ \alpha _{k}\right] }:=I_{_{%
\left[ \alpha _{k}\right] +1}}\left( e_{\left[ \alpha _{k}\right] -1}+e_{%
\left[ \alpha _{k}\right] }\right) .$ Since $\left[ \alpha
_{k}-2^{\left\vert \alpha _{k}\right\vert }\right] =\left[ \alpha _{k}\right],$
analogously to (\ref{5.1.9aaa}), if we apply Lemma \ref{lemma7} for $IV_{3}$ we have the following estimation
\begin{eqnarray*}
\left\vert IV_{3}\right\vert &=&\frac{2^{\left\vert \alpha _{k}\right\vert
\left( 1/p-1\right) }}{\Phi \left( \alpha _{k}\right) u\left( \alpha
_{k}\right) \alpha _{k}}\left( \alpha _{k}-2^{\left\vert \alpha
_{k}\right\vert }\right) \left\vert K_{\alpha _{k}-2^{\left\vert \alpha
_{k}\right\vert }}\right\vert \\
&=&\frac{2^{\left\vert \alpha _{k}\right\vert \left( 1/p-1\right) }}{\Phi
\left( \alpha _{k}\right) u\left( \alpha _{k}\right) \alpha _{k}}\left\vert
2^{\left[ \alpha _{k}\right] }K_{\left[ \alpha _{k}\right] }\right\vert \\
&\geq &
\frac{2^{\left\vert \alpha _{k}\right\vert \left( 1/p-2\right) }2^{2\left[
\alpha _{k}\right] -4}}{\Phi \left( \alpha _{k}\right) u\left( \alpha
_{k}\right) } \\
&\geq& 2^{\left\vert \alpha _{k}\right\vert \left( 1/p-2\right) /2}2^{2\left[
\alpha _{k}\right] -4}/\Phi ^{1/2}\left( \alpha _{k}\right) .
\end{eqnarray*}
Hence,
\begin{eqnarray*}
&&\left\Vert IV_{3}\right\Vert _{weak-L_{p}(G)}^{p} \\
&\geq & \left( \frac{2^{\left\vert \alpha _{k}\right\vert \left( 1/p-2\right)
/2}2^{2\left[ \alpha _{k}\right] -4}}{\Phi ^{1/2}\left( \alpha _{k}\right) }%
\right) ^{p}\mu \left\{ x\in G:\text{ }\left\vert IV_{3}\right\vert \geq
\frac{2^{\left\vert \alpha _{k}\right\vert \left( 1/p-2\right) /2}2^{2\left[
\alpha _{k}\right] -4}}{\Phi ^{1/2}\left( \alpha _{k}\right) }\right\} \\
&\geq & c_{p}\left( 2^{2\left[ \alpha _{k}\right] +\left\vert \alpha
_{k}\right\vert \left( 1/p-2\right) /2}/\Phi ^{1/2}\left( \alpha _{k}\right)
\right) ^{p}\mu (E_{\left[ \alpha _{k}\right] })\\
&\geq & c_{p}\left( 2^{\left( \left\vert \alpha _{k}\right\vert -\left[ \alpha
_{k}\right] \right) \left( 1/p-2\right) }/\Phi \left( \alpha _{k}\right)
\right) ^{p/2} \\
&=& c_{p}\left( 2^{d\left( \alpha _{k}\right) \left( 1/p-2\right)
}/\Phi \left( \alpha _{k}\right) \right) ^{p/2}\rightarrow \infty ,\text{\ as }k\rightarrow \infty .
\end{eqnarray*}

By combining Corollary \ref{cor4.1.2} and first part of Theorem \ref{th5.1.2} we find that
\begin{equation*}
\left\Vert IV_{1}\right\Vert _{weak-L_{p}(G)}\leq
c_{p}<\infty , \text{ \ \ } \left\Vert IV_{2}\right\Vert _{weak-L_{p}(G)}\leq
c_{p}<\infty.
\end{equation*}

On the other hand, for sufficiently large $ n $ we can conclude that
\begin{eqnarray*}
&&\left\Vert \sigma _{\alpha _{k}}f\right\Vert _{weak-L_{p}(G)}^{p} \\
&\geq&\left\Vert IV_{3}\right\Vert _{weak-L_{p}(G)}^{p}-\left\Vert
IV_{2}\right\Vert _{weak-L_{p}(G)}^{p}-\left\Vert IV_{1}\right\Vert
_{weak-L_{p}(G)}^{p} \\
&\geq & \frac{1}{2}\left\Vert IV_{3}\right\Vert _{weak-L_{p}(G)}^{p}\rightarrow \infty ,\text{ as }k\rightarrow \infty .
\end{eqnarray*}

Theorem \ref{th5.1.2} is proved.
\QED

The proofs of Corollaries \ref{cor5.1.1}-\ref{cor5.1.3} are similar to the proofs of Corollaries \ref{cor4.1.2}-\ref{cor4.1.4}. So, we leave out the details of proofs and just present these results:

\begin{corollary}\label{cor5.1.1}
Let $p>0$ and $f\in H_{p}(G)$. Then
\begin{equation*}
\left\Vert \sigma_{2^{k}}f-f\right\Vert _{H_p(G)}\rightarrow 0,\text{ as }k\rightarrow\infty.
\end{equation*}
\end{corollary}

\begin{corollary}\label{cor5.1.2}
Let $p>0$ and $f\in H_{p}(G)$. Then
\begin{equation*}
\left\Vert \sigma_{2^{k}+2^{k-1}}f-f\right\Vert _{H_p(G)}\rightarrow 0,\text{ as }k\rightarrow\infty.
\end{equation*}
\end{corollary}

\begin{corollary}\label{cor5.1.3}
Let $0<p<1/2$. Then there exists a martingale $f\in H_{p}(G)$, such that
\begin{equation*}
\left\Vert \sigma_{2^{k}+1}f-f\right\Vert _{weak-L_{p}(G)}\nrightarrow 0,\text{ as }k\rightarrow\infty.
\end{equation*}

On the other hand, for any $f\in H_{1/2}(G)$ the following is true:
\begin{equation*}
\left\Vert \sigma_{2^{k}+1}f-f\right\Vert _{H_{1/2}(G)}\rightarrow 0,\text{ as }k\rightarrow\infty.
\end{equation*}

\end{corollary}

\subsection{Modulus of continuity and convergence in norm of subsequences of  Fej\'er means with respect to the one-dimensional Walsh-Fourier series on the martingale Hardy spaces}

\text{ \qquad } In this section we apply Theorem \ref{th5.1.1} and Theorem \ref{th5.1.2} to find necessary and sufficient conditions for modulus of continuity of martingale $ f\in H_p $, for which subsequences of  Fej\'er means with respect to the one-dimensional Walsh-Fourier series converge in $ H_p $-norm.

First, we prove the following result:

\begin{theorem} \label{theorem5.2.1}
a) Let $f\in H_{1/2}(G), \ \ \sup_{k\in \mathbb{N}_+}V\left( n_{k}\right) =\infty$ and
\begin{equation} \label{5.2.cond}
\omega _{H_{p}(G)}\left( 1/2^{\left\vert n_{k}\right\vert },f\right) =o\left(
1/V^{2}\left( n_{k}\right) \right) ,\text{ \ as \ }k\rightarrow \infty .
\end{equation}%
Then
\begin{equation*}
\left\Vert \sigma_{n_{k}}f-f\right\Vert _{H_{1/2}(G)}\rightarrow 0,\text{ as }k\rightarrow\infty.
\end{equation*}

b) Let $\sup_{k\in \mathbb{N}_+}V\left( n_{k}\right) =\infty .$ Then there exists a martingale $f\in H_{1/2}(G),$ such that
\begin{equation} \label{5.2.cond2}
\omega _{H_{1/2}(G)}\left( 1/2^{\left\vert n_{k}\right\vert },f\right) =O\left(
1/V^{2}\left( n_{k}\right) \right) ,\text{ \ as \ }k\rightarrow \infty
\end{equation}%
and
\begin{equation} \label{5.2.kn3}
\left\Vert \sigma _{n_{k}}f-f\right\Vert _{H_{1/2}(G)}\nrightarrow 0, \text{ \ as \ }k\rightarrow \infty .
\end{equation}
\end{theorem}

{\bf Proof}:
Let $f\in H_{1/2}(G)$ and $2^{k}<n\leq 2^{k+1}.$ Then
\begin{eqnarray*} \\
&&
\left\Vert \sigma _{n}f-f\right\Vert _{H_{1/2}(G)}^{1/2}
\\
&&
\leq \left\Vert \sigma _{n}f-\sigma _{n}S_{2^{k}}f\right\Vert
_{H_{1/2}(G)}^{1/2}\\
&&
+\left\Vert \sigma _{n}S_{2^{k}}f-S_{2^{k}}f\right\Vert
_{H_{1/2}(G)}^{1/2}\\
&&
+\left\Vert S_{2^{k}}f-f\right\Vert _{H_{1/2}(G)}^{1/2}
\\
&&
=\left\Vert \sigma _{n}\left( S_{2^{k}}f-f\right) \right\Vert
_{H_{1/2}(G)}^{1/2}\\
&&
+\left\Vert S_{2^{k}}f-f\right\Vert _{H_{1/2}(G)}^{1/2}
\\
&&+\left\Vert
\sigma _{n}S_{2^{k}}f-S_{2^{k}}f\right\Vert _{H_{1/2}(G)}^{1/2}
\\
&&
\leq c\left( V\left( n\right) +1\right) \omega
_{H_{1/2}(G)}^{1/2}\left( 1/2^{k},f\right)\\
&&+\left\Vert \sigma
_{n}S_{2^{k}}f-S_{2^{k}}f\right\Vert _{H_{1/2}(G)}^{1/2}.\\
\end{eqnarray*}

It is evident that
\begin{eqnarray*}
&&\sigma _{n}S_{2^{k}}f-S_{2^{k}}f\\
&&
=\frac{2^{k}}{n}\left( S_{2^{k}}\sigma
_{2^{k}}f-S_{2^{k}}f\right) \\
&&
 =\frac{2^{k}}{n}S_{2^{k}}\left( \sigma
_{2^{k}}f-f\right).
\end{eqnarray*}

Let $p>0.$ By combining Corollaries \ref{cor4.1.2} and \ref{cor5.1.1} we can conclude that
\begin{eqnarray*}
&&\left\Vert \sigma _{n}S_{2^{k}}f-S_{2^{k}}f\right\Vert _{H_{1/2}(G)}^{1/2}
\\
&&
\leq \frac{2^{k/2}}{n^{1/2}}\left\Vert S_{2^{k}}\left( \sigma
_{2^{k}}f-f\right) \right\Vert _{H_{1/2}(G)}^{1/2}
\\
&&
\leq \left\Vert \sigma
_{2^{k}}f-f\right\Vert _{H_{1/2}(G)}^{1/2}\rightarrow 0,\text{\ as \ } k\rightarrow
\infty.
\end{eqnarray*}

Now, we prove part b) of Theorem \ref{theorem5.2.1}. Since $\sup_{k\in \mathbb{N}_+}V(\alpha _{k})=\infty,$ then there exists a martingale  $\{\alpha _{k}:k\in \mathbb{N}_+\}\subset \{n_{k}:k\in \mathbb{N}_+\}$ such that $V(\alpha _{k})\uparrow \infty,$ as $k\rightarrow \infty $ and
\begin{equation} \label{5.2.4.7}
V^{2}(\alpha _{k})\leq V(\alpha _{k+1}).
\end{equation}

Let $f$ be martingale from Example \ref{example2.2.1}, where
\begin{equation*}
\lambda_{k}=V^{-2}(\alpha _{k}),
\end{equation*}

If we apply (\ref{5.2.4.7}) we get that condition (\ref{3.3.2aa}) is fulfilled and it follows that $f\in H_{p}(G).$
By using (\ref{3.3.10AA}) we find that

\begin{equation} \label{5.2.4.13}
\widehat{f}(j)=\left\{
\begin{array}{ll}
2^{_{\left\vert \alpha _{k}\right\vert }}/V^{2}(\alpha _{k}), & \text{%
\thinspace }j\in \left\{ 2^{\left\vert \alpha _{k}\right\vert
},...,2^{_{\left\vert \alpha _{k}\right\vert +1}}-1\right\} ,\text{ }%
k\in\mathbb{N}_+ \\
0\,, & j\notin \bigcup\limits_{k=0}^{\infty }\left\{ 2^{_{\left\vert \alpha
_{k}\right\vert }},...,2^{_{\left\vert \alpha _{k}\right\vert +1}}-1\right\}
.\text{ }%
\end{array}%
\right.
\end{equation}

By combining (\ref{3.3.2aa0}) and (\ref{5.2.4.7}) we can conclude that
\begin{eqnarray} \label{4.2.4.12}
&& w_{H_{1/2}(G)}(1/2^n,f)
=\left\Vert f-S_{2^{n}}f\right\Vert_{H_{1/2}(G)} \\ \notag
&&\leq
\sum\limits_{i=n+1}^{\infty }1/V^{2}(\alpha _{i})= O\left( 1/V^{2}(\alpha
_{n})\right),\text{ \ as \ } n\rightarrow \infty.
\end{eqnarray}

Let $2^{_{\left\vert \alpha _{k}\right\vert }}<j\leq \alpha _{k}.$ By using (\ref{3.3.11AA}) we get that
\begin{equation*}
S_{j}f=S_{2^{_{\left\vert \alpha _{k}\right\vert }}}f+\frac{%
2^{_{\left\vert \alpha _{k}\right\vert }}w_{2^{_{\left\vert \alpha
_{k}\right\vert }}}D_{j-2^{_{\left\vert \alpha _{k}\right\vert }}}}{%
V^{2}(\alpha _{k})}.
\end{equation*}%
Hence,
\begin{eqnarray}  \label{5.2.nn}
&&\sigma _{\alpha _{k}}f-f \\ \notag
&=&\frac{2^{_{\left\vert \alpha _{k}\right\vert }}}{%
\alpha _{k}}\left( \sigma _{2^{_{\left\vert \alpha _{k}\right\vert
}}}f-f\right) \\ \notag
&+&\frac{\alpha _{k}-2^{_{\left\vert \alpha _{k}\right\vert }}}{\alpha _{k}}%
\left(S_{2^{_{\left\vert \alpha _{k}\right\vert }}}f-f\right) \\ \notag
&+&\frac{%
2^{_{\left\vert \alpha _{k}\right\vert }}w_{2^{_{\left\vert \alpha
_{k}\right\vert }}}\left( \alpha _{k}-2^{_{\left\vert \alpha _{k}\right\vert
}}\right) K_{\alpha _{k}-2^{_{\left\vert \alpha _{k}\right\vert }}}}{\alpha_{k}V^{2}(\alpha _{k})}.
\end{eqnarray}

According to (\ref{1.S2n000}),  (\ref{fe22222}) and (\ref{5.2.nn}) we have that
\begin{eqnarray} \label{5.2.a11}
&&\Vert \sigma _{\alpha _{k}}f-f\Vert _{1/2}^{1/2} \\ \notag
&\geq & \frac{c}{V(\alpha _{k})}%
\Vert \left( \alpha _{k}-2^{_{\left\vert \alpha _{k}\right\vert }}\right)
K_{\alpha _{k}-2^{_{\left\vert \alpha _{k}\right\vert }}}\Vert _{1/2}^{1/2} \\ \notag
&-&\left( \frac{2^{_{\left\vert \alpha _{k}\right\vert }}}{\alpha _{k}}\right)
^{1/2}\Vert \sigma _{2^{_{\left\vert \alpha _{k}\right\vert }}}f-f\Vert
_{1/2}^{1/2} \\ \notag
&-&\left( \frac{\alpha _{k}-2^{_{\left\vert \alpha _{k}\right\vert
}}}{\alpha _{k}}\right) ^{1/2}\Vert S_{2^{_{\left\vert \alpha
_{k}\right\vert }}}f-f\Vert _{1/2}^{1/2}.
\end{eqnarray}

Let
$$\alpha _{k}=\sum_{i=1}^{r_{k}}\sum_{k=l_{i}^{k}}^{m_{i}^{k}}2^{k},$$
where
$$m_{1}^{k}\geq l_{1}^{k}>l_{1}^{k}-2\geq m_{2}^{k}\geq
l_{2}^{k}>l_{2}^{k}-2>...>m_{s}^{k}\geq l_{s}^{k}\geq 0$$
and
$$E_{l_{i}^{k}}:=I_{_{l_{i}^{k}+1}}\left(e_{l_{i}^{k}-1}+e_{l_{i}^{k}}\right).$$

By using Lemma \ref{lemma7} we get that
\begin{eqnarray} \label{5.2.33}
&&\int_{G}\left\vert \left( \alpha _{k}-2^{_{\left\vert \alpha _{k}\right\vert
}}\right) K_{\alpha _{k}-2^{_{\left\vert \alpha _{k}\right\vert }}}\left(
x\right) \right\vert ^{1/2}d\mu   \\ \notag
&\geq &\frac{1}{16}\underset{i=2}{\overset{r_{k}-2}{\sum }}%
\int_{E_{l_{i}^{k}}}\left\vert \left( \alpha _{k}-2^{_{\left\vert \alpha
_{k}\right\vert }}\right) K_{\alpha _{k}-2^{_{\left\vert \alpha
_{k}\right\vert }}}\left( x\right) \right\vert ^{1/2}d\mu \left( x\right)
\\ \notag
&\geq & \frac{1}{16}\underset{i=2}{\overset{r_{k}-2}{\sum }}\frac{1}{%
2^{l_{i}^{k}}}2^{l_{i}^{k}}
\geq cr_{k}\geq cV(\alpha _{k}).
\end{eqnarray}

By combining estimations (\ref{5.2.a11}-\ref{5.2.33}), Corollaries \ref{cor4.1.2} and \ref{cor5.1.1} we get that (\ref{5.2.kn3}) holds true and Theorem \ref{theorem5.2.1} is proved.
\QED

\begin{theorem}\label{theorem5.2.2}
a) Let $0<p<1/2,$ $f\in H_{p}(G)$, $\ \sup_{k\in \mathbb{N}_+}d\left(
n_{k}\right) =\infty $ and
\begin{equation}
\omega _{H_{p}(G)}\left( 1/2^{\left\vert n_{k}\right\vert },f\right) =o\left(
1/2^{d\left( n_{k}\right) \left( 1/p-2\right) }\right) ,\text{ as \ }%
k\rightarrow \infty .  \label{5.2.cond3}
\end{equation}%
Then
\begin{equation} \label{5.2.fe2}
\left\Vert\sigma_{n_{k}}f-f\right\Vert_{H_{p}(G)}\rightarrow0,\text{\ as \ } k\rightarrow \infty.
\end{equation}

b) Let $\sup_{k\in \mathbb{N}_+}d\left( n_{k}\right) =\infty .$ Then there exists a martingale $f\in H_{p}(G)$ $\left( 0<p<1/2\right) ,$\ \ such that
\begin{eqnarray} \label{cond4}
&&\omega _{H_{p}(G)}\left( 1/2^{\left\vert n_{k}\right\vert },f\right)= O\left(
1/2^{d\left( n_{k}\right) \left( 1/p-2\right) }\right) ,\text{ \ as \ }%
k\rightarrow \infty
\end{eqnarray}%
and
\begin{equation}
\left\Vert\sigma_{n_{k}}f-f\right\Vert_{weak-L_{p}(G)}\nrightarrow
0,\text{\ as \ }k\rightarrow \infty.
\label{kn4}
\end{equation}
\end{theorem}

{\bf Proof}:
\textbf{\ }Let $0<p<1/2.$ Then under condition (\ref{5.2.cond3}) if we repeat steps of the proof of Theorem \ref{theorem5.2.1}, we immediately get that (\ref{5.2.fe2}) holds.

Let prove part b) of Theorem \ref{theorem5.2.2}. Since  $\sup_{k}d\left(
n_{k}\right) =\infty ,$ there exists $\{\alpha _{k}:k\in \mathbb{N}_+\}\subset\{n_{k}:k\in \mathbb{N}_+\}$ such that $\sup_{k\in \mathbb{N}_+}d\left( \alpha _{k}\right) =\infty $
and
\begin{equation} \label{5.2.4.18}
2^{2d\left( \alpha _{k}\right) \left( 1/p-2\right) }\leq 2^{d\left( \alpha
_{k+1}\right) \left( 1/p-2\right) }.
\end{equation}

Let $f$ be a martingale from Lemma \ref{example2.2.1}, where
\begin{equation*}
\lambda_{k}=2^{-\left( 1/p-2\right) d\left( \alpha _{i}\right)}.
\end{equation*}

If we use (\ref{5.2.4.18}) we conclude that condition (\ref{3.3.2aa}) is fulfilled and it follows that $f\in H_{p}(G).$

According to (\ref{3.3.10AA}) we get that

\begin{equation} \label{5.2.4.22}
\widehat{f}(j)=\left\{
\begin{array}{ll}
2^{\left( 1/p-2\right) \left[ \alpha _{k}\right] }, & \text{\thinspace
\thinspace }j\in \left\{ 2^{\left\vert \alpha _{k}\right\vert
},...,2^{_{\left\vert \alpha _{k}\right\vert +1}}-1\right\} ,\text{ }%
k\in \mathbb{N}_+ \\
0\,, & \text{\thinspace }j\notin \bigcup\limits_{n=0}^{\infty }\left\{
2^{_{\left\vert \alpha _{n}\right\vert }},...,2^{_{\left\vert \alpha
_{n}\right\vert +1}}-1\right\} .\text{ }%
\end{array}%
\right.
\end{equation}

By combining (\ref{3.3.2aa0}) and (\ref{5.2.4.18}) we have that
\begin{eqnarray} \label{4.21}
&&\omega _{H_{p}(G)}(1/2^{\left\vert \alpha _{k}\right\vert },f)  \\ \notag
&\leq & \sum\limits_{i=k}^{\infty }1/2^{d\left( \alpha _{i}\right) \left(
1/p-2\right) } \\ \notag
&=& O\left( 1/2^{d\left( \alpha _{k}\right) \left( 1/p-2\right)
}\right)\text{ \ as \ } k \rightarrow \infty.
\end{eqnarray}

Analogously to the proof of previous theorem, if we use also Corollaries \ref{cor4.1.2} and \ref{cor5.1.1}, for the sufficiently large $ k $ we can conclude that
\begin{eqnarray} \label{5.2.nnmm00}
&&\Vert \sigma _{\alpha _{k}}f-f\Vert _{weak-L_{p}(G)}^{p} \\ \notag
&\geq& 2^{\left(
1-2p\right) \left[ \alpha _{k}\right] }\Vert \left( \alpha
_{k}-2^{\left\vert \alpha _{k}\right\vert }\right) K_{\alpha
_{k}-2^{_{\left\vert \alpha _{k}\right\vert }}}\Vert _{weak-L_{p}(G)}^{p} \\  \notag
&-&\left( \frac{2^{\left\vert \alpha _{k}\right\vert }}{\alpha _{k}}\right)
^{p}\Vert \sigma _{2^{_{\left\vert \alpha _{k}\right\vert }}}f-f\Vert
_{weak-L_{p}(G)}^{p} \\ \notag
&-&\left( \frac{\alpha _{k}-2^{_{\left\vert \alpha
_{k}\right\vert }}}{\alpha _{k}}\right) ^{p}\Vert S_{2^{_{\left\vert \alpha
_{k}\right\vert }}}f-f\Vert _{weak-L_{p}(G)}^{p} \\ \notag
&\geq & 2^{\left(
1-2p\right) \left[ \alpha _{k}\right]-1 }\Vert \left( \alpha
_{k}-2^{\left\vert \alpha _{k}\right\vert }\right) K_{\alpha
_{k}-2^{_{\left\vert \alpha _{k}\right\vert }}}\Vert _{weak-L_{p}(G)}^{p}
\end{eqnarray}

Let $x\in E_{\left[ \alpha _{k}\right] }.$ Lemma \ref{lemma7} follows that
\begin{eqnarray*}
&&\mu \left( x\in G:\left( \alpha _{k}-2^{_{\left\vert \alpha _{k}\right\vert
}}\right) \left\vert K_{\alpha _{k}-2^{_{\left\vert \alpha _{k}\right\vert
}}}\right\vert \geq 2^{2\left[ \alpha _{k}\right] -4}\right) \\
&&\geq \mu \left(
E_{\left[ \alpha _{k}\right] }\right) \geq 1/2^{\left[ \alpha _{k}\right]
-4}
\end{eqnarray*}
and
\begin{eqnarray} \label{5.2.nnmm}
&& 2^{2p\left[ \alpha _{k}\right] -4}\mu \left( x\in G:\left( \alpha
_{k}-2^{_{\left\vert \alpha _{k}\right\vert }}\right) \left\vert K_{\alpha
_{k}-2^{_{\left\vert \alpha _{k}\right\vert }}}\right\vert \geq 2^{2\left[
\alpha _{k}\right] -4}\right) \\ \notag
&& \geq 2^{\left( 2p-1\right) \left[ \alpha _{k}%
\right] -4}.
\end{eqnarray}

Hence, by combining (\ref{1.S2n000}),  (\ref{fe22222}), (\ref{5.2.nnmm00}) and (\ref{5.2.nnmm}) we get that
\begin{equation*}
\left\Vert \sigma _{n_{k}}f-f\right\Vert _{weak-L_p(G)}\nrightarrow
0,\,\,\,\text{as}\,\,\,k\rightarrow \infty.
\end{equation*}%
The proof of Theorem \ref{theorem5.2.2} is complete.
\QED

By using Theorem \ref{theorem5.2.2} we easily get an important result which was proved in \cite{tep6}:

\begin{corollary}\label{corollary5.2.2ss2}
a) Let  $f\in H_{1/2}(G)$ and
\begin{equation*}
\omega _{H_{1/2}(G)}\left( 1/2^{k},f\right) =o\left(\frac{1}{k^2}\right) ,\text{ as \ }%
k\rightarrow \infty .
\end{equation*}%
Then
\begin{equation*}
\left\Vert \sigma _{k}f-f\right\Vert _{H_{1/2}(G)}\rightarrow
0,\,\,\,\text{ \ as \ }k\rightarrow \infty.
\end{equation*}

b) There exists a martingale $f\in H_{1/2}(G),$ for which
\begin{equation*}
\omega _{H_{1/2}(G)}\left( 1/2^{k},f\right) =O\left(\frac{1}{k^2}\right),\text{ \ as \ }
k\rightarrow \infty
\end{equation*}%
and
\begin{equation*}
\left\Vert \sigma_{k}f-f\right\Vert_{1/2}\nrightarrow0,\,\,\,\text{ \ as \ }k\rightarrow\infty .
\end{equation*}
\end{corollary}

\begin{corollary}\label{corollary5.2.2ss}
a) Let $0<p<1/2,$ $f\in H_{p}(G)$ and

\begin{equation*}
\omega _{H_{p}(G)}\left( 1/2^{k},f\right) =o\left(
1/2^{k( 1/p-2)}\right) ,\text{ as \ }%
k\rightarrow \infty .
\end{equation*}%
Then
\begin{equation*}
\left\Vert \sigma _{k}f-f\right\Vert _{H_{p}(G)}\rightarrow
0,\,\,\,\text{as\thinspace \thinspace \thinspace }k\rightarrow \infty .
\end{equation*}

b) Then there exists a martingale $f\in H_{p}(G)$ $\left( 0<p<1/2\right) ,$\ \ for which

\begin{equation*}
\omega _{H_{p}(G)}\left( 1/2^{k},f\right) =O\left(
1/2^{k\left( 1/p-2\right) }\right) ,\text{ \ as \ }%
k\rightarrow \infty
\end{equation*}%
and
\begin{equation*}
\left\Vert \sigma _{k}f-f\right\Vert _{weak-L_{p}(G)}\nrightarrow
0,\,\,\,\text{ \ as \ }k\rightarrow \infty .
\end{equation*}
\end{corollary}

\subsection{Strong convergence of  Fej\'er means with respect to the one-dimensional Walsh-Fourier series on the martingale Hardy spaces}

\text{ \qquad } In this section we consider strong convergence results of Fej\'er means with respect to the one-dimensional Walsh-Fourier series in the martingale Hardy spaces, when $ 0<p\leq 1/2 $ (for details see \cite{tep5}).

The following is true:

\begin{theorem} \label{theorem5.3.1}
a) Let $0<p\leq 1/2$ and $f\in H_{p}(G)$. Then there exists a constant $c_{p},$ depending only on $p$, such that

\begin{equation*}
\frac{1}{\log ^{\left[ 1/2+p\right] }n}\overset{n}{\underset{m=1}{\sum }}%
\frac{\left\Vert \sigma _{m}f\right\Vert _{H_{p}(G)}^{p}}{m^{2-2p}}\leq
c_{p}\left\Vert f\right\Vert _{H_{p}(G)}^{p}.
\end{equation*}%
b) Let $0<p<1/2,$ $\Phi :\mathbb{N}_{+}\rightarrow\lbrack 1,\infty )$ be non-decreasing function, such that $\Phi \left( n\right) \uparrow \infty $ and
\begin{equation*}
\overline{\underset{k\rightarrow \infty }{\lim }}\frac{k^{
2-2p }}{\Phi \left({k}\right) }=\infty .
\end{equation*}

Then there exists a martingale $f\in H_{p}(G),$ such that
\begin{equation*}
\underset{m=1}{\overset{\infty }{\sum }}\frac{\left\Vert \sigma
_{m}f\right\Vert _{weak-L_{p}(G)}^{p}}{\Phi \left( m\right) }=\infty .
\end{equation*}
\end{theorem}

{\bf Proof}:
Suppose that
\begin{equation*}
\frac{1}{\log ^{\left[ 1/2+p\right] }n}\overset{n}{\underset{m=1}{\sum }}%
\frac{\left\Vert \sigma _{m}f\right\Vert _{p}^{p}}{m^{2-2p}}\leq
c_{p}\left\Vert f\right\Vert _{H_{p}(G)}^{p}.
\end{equation*}

By combining (\ref{1.S2n}), (\ref{sigmamax}) and Lemma \ref{lemma3.2.10} we can conclude that
\begin{eqnarray} \label{5.3.5.3}
&&\frac{1}{\log ^{\left[ 1/2+p\right] }n}\overset{n}{\underset{m=1}{\sum }}%
\frac{\left\Vert \sigma _{m}f\right\Vert _{H_{p}(G)}^{p}}{m^{2-2p}} \\
&& \leq\frac{1}{\log ^{\left[ 1/2+p\right] }n}\overset{n}{\underset{m=1}{\sum }}%
\frac{\left\Vert \sigma _{m}f\right\Vert _{p}^{p}}{m^{2-2p}}+\|\widetilde{\sigma}_{\#}^{\ast }f\|_{H_p(G)}+\|\widetilde{S}_{\#}^{\ast }f\|_{H_p(G)}  \notag \\
&&  \leq c_{p}\left\Vert f\right\Vert _{H_{p}(G)}^{p}. \notag
\end{eqnarray}

According to Lemma \ref{lemma3.2.5} and (\ref{5.3.5.3}) Theorem \ref{theorem5.3.1} will be proved if we show that

\begin{equation*}
\frac{1}{\log ^{\left[ 1/2+p\right] }n}\overset{n}{\underset{m=1}{\sum }}%
\frac{\left\Vert \sigma _{m}a\right\Vert _{p}^{p}}{m^{2-2p}}\leq c<\infty ,%
\text{ \ \ \ }m=2,3,...
\end{equation*}%
for any  $ p $-atom $a$. We may assume that $a$ is $ p $-atom, with support  $ I$, $\mu \left( I\right) =2^{-M}$ and $I=I_{M}.$ It is evident that $\sigma _{n}\left( a\right) =0,$ when $n\leq 2^{M}.$ Therefore, we may assume that $n>2^{M}.$

Let $x\in I_{M}.$ Since $\sigma _{n}$ is bounded from $L_{\infty }(G)$ to $L_{\infty }(G)$ (The boundedness follows fact that Fej\'er kernels are uniformly bounded in the space $L_{1}(G)$, which is proved in Lemma \ref{lemma3}) and $\left\Vert
a\right\Vert _{\infty }\leq 2^{M/p}$ we can conclude that
\begin{equation*}
\int_{I_{M}}\left\vert \sigma _{m}a\left( x\right) \right\vert ^{p}d\mu
\left( x\right) \leq \left\Vert \sigma _{m}a\right\Vert _{\infty
}^{p}/2^{M}
\end{equation*}
\begin{equation*}
\leq \left\Vert a\right\Vert _{\infty
}^{p}/2^{M}\leq c<\infty ,\text{ }0<p\leq 1/2.
\end{equation*}

Let $0<p\leq 1/2.$ Then
\begin{equation*}
\frac{1}{\log ^{\left[ 1/2+p\right] }n}\overset{n}{\underset{m=1}{\sum }}%
\frac{\int_{I_{M}}\left\vert \sigma _{m}a\left( x\right) \right\vert
^{p}d\mu \left( x\right) }{m^{2-2p}}
\end{equation*}
\begin{equation*}
\leq \frac{c}{\log ^{\left[ 1/2+p\right]
}n}\overset{n}{\underset{m=1}{\sum }}\frac{1}{m^{2-2p}}\leq c<\infty .
\end{equation*}

It is evident that
\begin{equation*}
\left\vert \sigma _{m}a\left( x\right) \right\vert \leq
\int_{I_{M}}\left\vert a\left( t\right) \right\vert \left\vert K_{m}\left(
x+t\right) \right\vert d\mu \left( t\right)
\end{equation*}%
\begin{equation*}
\leq
2^{M/p}\int_{I_{M}}\left\vert K_{m}\left( x+t\right) \right\vert d\mu \left(
t\right) .
\end{equation*}

Lemma \ref{lemma4} follows that
\begin{equation} \label{5.3.12}
\left\vert \sigma _{m}a\left( x\right) \right\vert \leq \frac{%
c2^{k+l}2^{M\left( 1/p-1\right) }}{m},\text{ \ \ }x\in I_{l+1}\left(
e_{k}+e_{l}\right) ,\,0\leq k<l<M
\end{equation}%
and
\begin{equation} \label{5.3.12a}
\left\vert \sigma _{m}a\left( x\right) \right\vert \leq c2^{M\left(
1/p-1\right) }2^{k},\text{ \ }x\in I_{M}\left( e_{k}\right) ,\,0\leq k<M.
\end{equation}

If we use identity (\ref{2.1.2}) and (\ref{5.3.12}-\ref{5.3.12a}) we get that
\begin{eqnarray} \label{7.2}
&&\int_{\overline{I_{M}}}\left\vert \sigma _{m}a\left( x\right) \right\vert
^{p}d\mu \left( x\right)  \\
&=&\overset{M-2}{\underset{k=0}{\sum }}\overset{M-1}{\underset{l=k+1}{\sum }}%
\int_{I_{l+1}\left( e_{k}+e_{l}\right) }\left\vert \sigma _{m}a\left(
x\right) \right\vert ^{p}d\mu \left( x\right)
\notag \\
& &
+\overset{M-1}{\underset{k=0}{%
\sum }}\int_{I_{M}\left( e_{k}\right) }\left\vert \sigma _{m}a\left(
x\right) \right\vert ^{p}d\mu \left( x\right)  \notag \\
&\leq &c\overset{M-2}{\underset{k=0}{\sum }}\overset{M-1}{\underset{l=k+1}{%
\sum }}\frac{1}{2^{l}}\frac{2^{p\left( k+l\right) }2^{M\left( 1-p\right) }}{%
m^{p}}+c\overset{M-1}{\underset{k=0}{\sum }}\frac{1}{2^{M}}2^{M\left(
1-p\right) }2^{pk}  \notag \\
&\leq &\frac{c2^{M\left( 1-p\right) }}{m^{p}}\overset{M-2}{\underset{k=0}{%
\sum }}\overset{M-1}{\underset{l=k+1}{\sum }}\frac{2^{p\left( k+l\right) }}{%
2^{l}}+c\overset{M-1}{\underset{k=0}{\sum }}\frac{2^{pk}}{2^{pM}}  \notag \\
&\leq &\frac{c2^{M\left( 1-p\right) }M^{\left[ 1/2+p\right] }}{m^{p}}+c.
\notag
\end{eqnarray}

Hence,
\begin{eqnarray*}
&&\frac{1}{\log ^{\left[ 1/2+p\right] }n}\overset{n}{\underset{m=2^{M}+1}{%
\sum }}\frac{\int_{\overline{I_{M}}}\left\vert \sigma _{m}a\left( x\right)
\right\vert ^{p}d\mu \left( x\right) }{m^{2-2p}} \\
&\leq &\frac{1}{\log ^{\left[ 1/2+p\right] }n}\left( \overset{n}{\underset{%
m=2^{M}+1}{\sum }}\frac{c2^{M\left( 1-p\right) }M^{\left[ 1/2+p\right] }}{%
m^{2-p}}+\overset{n}{\underset{m=2^{M}+1}{\sum }}\frac{c}{m^{2-2p}}\right)
<c<\infty .
\end{eqnarray*}

The proof of part a) of theorem \ref{theorem5.3.1} is complete.

\bigskip Now, we prove part b) of Theorem \ref{theorem5.3.1}.  Let $\Phi \left( n\right) $ non-decreasing function satisfying the condition
\begin{equation} \label{5.3.12j}
\underset{k\rightarrow \infty }{\lim }\frac{2^{\left( \left\vert
n_{k}\right\vert +1\right) \left( 2-2p\right) }}{\Phi \left( 2^{\left\vert
n_{k}\right\vert +1}\right) }=\infty .
\end{equation}

According to (\ref{5.3.12j}), there exists an increasing sequence  $\left\{ \alpha _{k}:%
\text{ }k\in\mathbb{N}_+\right\} \subset \left\{ n_{k}:\text{ }k\in\mathbb{N}_+\right\} $ such that
\begin{equation} \label{5.3.122}
\left\vert \alpha _{k}\right\vert \geq 2,\text{ \ \ \ where \ }k\in\mathbb{N}_+
\end{equation}%
and
\begin{eqnarray} \label{5.3.121}
&&\sum_{\eta =0}^{\infty }\frac{\Phi ^{1/2}\left( 2^{\left\vert \alpha _{\eta
}\right\vert +1}\right) }{2^{\left\vert \alpha _{\eta }\right\vert \left(
1-p\right) }} \\ \notag
&=& 2^{1-p}\sum_{\eta =0}^{\infty }\frac{\Phi ^{1/2}\left(
2^{\left\vert \alpha _{\eta }\right\vert +1}\right) }{2^{\left( \left\vert
\alpha _{\eta }\right\vert +1\right) \left( 1-p\right) }}<c<\infty .
\end{eqnarray}

Let $f=\left(f_{n},\text{ }n\in\mathbb{N}_+\right) \in H_{p}(G)$ be a martingale from the Example \ref{example2.2.1}, where
\begin{equation*}
\lambda _{k}=\frac{\Phi ^{1/2p}\left( 2^{\left\vert \alpha _{k}\right\vert
+1}\right) }{2^{\left( \left\vert \alpha _{k}\right\vert \right) \left(
1/p-1\right) }}
\end{equation*}

By combining (\ref{3.3.2aa}) and (\ref{5.3.121}) we get that $f\in H_{p}(G).$
According to (\ref{3.3.10AA}) we have that
\begin{eqnarray} \label{5.3.6aa}
&&\widehat{f}(j)=\left\{
\begin{array}{l}
\Phi ^{1/2p}\left(2^{\left\vert \alpha _{k}\right\vert +1}\right) 
\text{ \ if \ }j\in \left\{ 2^{\left\vert \alpha
_{k}\right\vert },...,2^{\left\vert \alpha _{k}\right\vert +1}-1\right\} ,%
\text{ }k\in\mathbb{N}_+, \\
0\text{ },\text{  \ if \ } j\notin
\bigcup\limits_{k=0}^{\infty }\left\{ 2^{\left\vert \alpha _{k}\right\vert
},...,2^{\left\vert \alpha _{k}\right\vert +1}-1\right\} .\text{ }
\end{array}
\right. 
\end{eqnarray}

Let $2^{\left\vert \alpha _{k}\right\vert }<n<2^{\left\vert \alpha
_{k}\right\vert +1}.$ Then

\begin{equation} \label{5.3.7aa}
\sigma _{_{n}}f=\frac{1}{n}\sum_{j=1}^{2^{\left\vert \alpha _{k}\right\vert
}}S_{j}f+\frac{1}{n}\sum_{j=2^{\left\vert \alpha _{k}\right\vert
}+1}^{n}S_{j}f=III+IV.
\end{equation}

It is evident that
\begin{equation} \label{5.3.7aaa}
S_{j}f=0,\,\ \text{ if \thinspace \thinspace }0\leq j\leq 2^{\left\vert \alpha
_{1}\right\vert }
\end{equation}

Let $2^{\left\vert \alpha _{s}\right\vert }<j\leq 2^{\left\vert
\alpha _{s}\right\vert +1},$ where $s=1,2,...,k.$ If we apply (\ref{3.3.11AA}) we get that
\begin{eqnarray} \label{5.3.8aa}
S_{j}f &=&\sum_{\eta =0}^{s-1}\Phi ^{1/2p}\left( 2^{\left\vert \alpha _{\eta}\right\vert +1}\right) \left( D_{2^{\left\vert \alpha _{\eta }\right\vert+1}}-D_{2^{\left\vert \alpha _{\eta }\right\vert }}\right)  \\
&&+\Phi ^{1/2p}\left( 2^{\left\vert \alpha _{s}\right\vert +1}\right) w_{2^{\left\vert \alpha _{s}\right\vert }}D_{j-2^{\left\vert \alpha _{s}\right\vert }}.  \notag
\end{eqnarray}

Let $2^{\left\vert \alpha _{s}\right\vert +1}\leq j\leq 2^{\left\vert \alpha_{s+1}\right\vert },$ $s=0,1,...k-1.$ Then if we use (\ref{3.3.12AA}) we can conclude that
\begin{equation} \label{5.3.10aaaa}
S_{j}f=\sum_{\eta =0}^{s}\Phi ^{1/2p}\left( 2^{\left\vert \alpha _{\eta
}\right\vert +1}\right) \left( D_{2^{\left\vert \alpha _{\eta }\right\vert
+1}}-D_{2^{\left\vert \alpha _{\eta }\right\vert }}\right).
\end{equation}

Let $x\in I_{2}\left( e_{0}+e_{1}\right).$ Since (see Lemmas \ref{lemma1} and \ref{lemma4})
\begin{equation} \label{5.3.40}
D_{2^{n}}\left( x\right) =K_{2^{n}}\left( x\right) =0,\text{ where }n\geq 2
\end{equation}
by combining (\ref{5.3.122}) and (\ref{5.3.7aaa}-\ref{5.3.40}) we get that
\begin{equation}  \label{5.3.9aaa}
III=\frac{1}{n}\sum_{\eta =0}^{k-1}\Phi ^{1/2p}\left( 2^{\left\vert \alpha
_{\eta }\right\vert +1}\right) \sum_{v=2^{\left\vert \alpha _{\eta
}\right\vert }+1}^{2^{\left\vert \alpha _{\eta }\right\vert +1}}D_{v}\left(
x\right)
\end{equation}%
\begin{equation*}
=\frac{1}{n}\sum_{\eta =0}^{k-1}\Phi ^{1/2p}\left( 2^{\left\vert \alpha
_{\eta }\right\vert +1}\right) \left( 2^{\left\vert \alpha _{\eta
}\right\vert +1}K_{2^{\left\vert \alpha _{\eta }\right\vert +1}}\left(
x\right) -2^{\left\vert \alpha _{\eta }\right\vert }K_{2^{\left\vert \alpha
_{\eta }\right\vert }}\left( x\right) \right) =0.
\end{equation*}

If we use (\ref{5.3.8aa}) when $s=k$  for $IV$ we can write that

\begin{eqnarray} \label{5.3.9aa}
IV &=&\frac{n-2^{\left\vert \alpha _{k}\right\vert }}{n}\sum_{\eta
=0}^{k-1}\Phi ^{1/2p}\left( 2^{\left\vert \alpha _{\eta }\right\vert
+1}\right) \left( D_{2^{\left\vert \alpha _{\eta }\right\vert
+1}}-D_{2^{\left\vert \alpha _{\eta }\right\vert }}\right)   \\ \notag
&+&\frac{\Phi ^{1/2p}\left( 2^{\left\vert \alpha _{k}\right\vert +1}\right)
}{n}\sum_{j=2^{_{\left\vert \alpha _{k}\right\vert }}+1}^{n}w_{2^{\left\vert \alpha _{k}\right\vert }}
D_{j-2^{\left\vert \alpha _{k}\right\vert }} \\
&=& IV_{1}+IV_{2}.
\notag
\end{eqnarray}

By combining (\ref{5.3.122}) and (\ref{5.3.40}) we can conclude that

\begin{equation}
IV_{1}=0,\text{ \ where \ }x\in I_{2}\left( e_{0}+e_{1}\right) .  \label{8aaaa}
\end{equation}

Let $\alpha _{k}\in \mathbb{A}_{0,2},$ $2^{\left\vert \alpha _{k}\right\vert
}<n<2^{\left\vert \alpha _{k}\right\vert +1}$ and $x\in I_{2}\left(
e_{0}+e_{1}\right) $. Since $n-2^{_{\left\vert \alpha _{k}\right\vert }}\in
\mathbb{A}_{0,2},$ Lemmas \ref{lemma0} and \ref{lemma3} and (\ref{5.3.40}) follows that
\begin{eqnarray}  \label{5.3.31}
&&\left\vert IV_{2}\right\vert
=\frac{\Phi ^{1/2p}\left( 2^{\left\vert \alpha _{k}\right\vert +1}\right)
}{n}\left\vert \sum_{j=1}^{n-2^{\left\vert _{\alpha _{k}}\right\vert
}}D_{_{j}}\left( x\right) \right\vert  \notag\\
&=&\frac{\Phi ^{1/2p}\left( 2^{\left\vert \alpha _{k}\right\vert +1}\right)
}{n}\left\vert \left( n-2^{_{\left\vert \alpha _{k}\right\vert }}\right)
K_{n-2^{_{\left\vert \alpha _{k}\right\vert }}}\left( x\right) \right\vert \notag\\
&\geq& \frac{\Phi ^{1/2p}\left( 2^{\left\vert \alpha _{k}\right\vert
+1}\right) }{2^{\left\vert \alpha _{k}\right\vert +1}}. \notag
\end{eqnarray}

Let $0<p<1/2$ and $n\in \mathbb{A}_{0,2}.$ By combining (\ref{5.3.7aa}-\ref{5.3.31}) we get that

\begin{eqnarray} \label{10aaa}
&&\left\Vert \sigma _{n}f\right\Vert _{weak-L_{p}(G)}^{p}   \\
&\geq &\frac{c_{p}\Phi ^{1/2}\left( 2^{\left\vert \alpha _{k}\right\vert
+1}\right) }{2^{p\left( \left\vert \alpha _{k}\right\vert +1\right) }}\mu
\left\{ x\in I_{2}\left( e_{0}+e_{1}\right) :\text{ }\left\vert \sigma
_{n}f\right\vert \geq \frac{c_{p}\Phi ^{1/2p}\left( 2^{\left\vert \alpha
_{k}\right\vert +1}\right) }{2^{\left\vert \alpha _{k}\right\vert +1}}%
\right\}   \notag \\
&\geq &\frac{c_{p}\Phi ^{1/2}\left( 2^{\left\vert \alpha _{k}\right\vert
+1}\right) }{2^{p\left( \left\vert \alpha _{k}\right\vert +1\right) }}\mu
\left\{ I_{2}\left( e_{0}+e_{1}\right) \right\} \notag\\
&\geq & \frac{c_{p}\Phi
^{1/2}\left( 2^{\left\vert \alpha _{k}\right\vert +1}\right) }{2^{p\left(
\left\vert \alpha _{k}\right\vert +1\right) }}.  \notag
\end{eqnarray}

Hence,
\begin{eqnarray*}
&&\underset{n=1}{\overset{\infty }{\sum }}\frac{\left\Vert \sigma
_{n}f\right\Vert _{weak-L_{p}(G)}^{p}}{\Phi \left( n\right) }\\
&\geq & \underset{\left\{ n\in \mathbb{A}_{0,2}:\text{ }2^{\left\vert \alpha _{k}\right\vert}<n<2^{\left\vert \alpha _{k}\right\vert +1}\right\} }{\sum }\frac{\left\Vert \sigma _{n}f\right\Vert _{weak-L_{p}(G)}^{p}}{\Phi \left(n\right) } \\
&\geq &\frac{1}{\Phi ^{1/2}\left( 2^{\left\vert \alpha _{k}\right\vert
+1}\right) }\underset{\left\{ n\in \mathbb{A}_{0,2}:\text{ }2^{\left\vert
\alpha _{k}\right\vert }<n<2^{\left\vert \alpha _{k}\right\vert +1}\right\} }%
{\sum }\frac{1}{2^{p\left( \left\vert \alpha _{k}\right\vert +1\right) }} \\
&\geq &\frac{c_{p}2^{\left( 1-p\right) \left( \left\vert \alpha
_{k}\right\vert +1\right) }}{\Phi ^{1/2}\left( 2^{\left\vert \alpha
_{k}\right\vert +1}\right) }\rightarrow \infty ,\text{ \qquad as \ \ }%
k\rightarrow \infty .
\end{eqnarray*}

The proof of Theorem \ref{theorem5.3.1} is complete.
\QED

\begin{theorem} \label{theorem5.3.2}
Let $f\in H_{1/2}(G).$ Then
\begin{equation*}
\underset{n\in \mathbf{%
\mathbb{N}
}_{+}}{\sup }\underset{\left\Vert f\right\Vert _{H_{p}}\leq 1}{\sup }\frac{1%
}{n}\underset{m=1}{\overset{n}{\sum }}\left\Vert \sigma _{m}f\right\Vert
_{1/2}^{1/2}=\infty .
\end{equation*}
\end{theorem}

{\bf Proof}:
Let $0<p\leq 1$ and
\begin{equation*}
f_{k}({x}):=2^{k}\left( D_{2^{k+1}}(x)-D_{2^k}(x)\right)
\end{equation*}%
Since
\begin{equation*}
\text{supp}(f_{k})=I_{k},\text{ \ \ }\int_{I_{k}}a_{k}d\mu =0
\end{equation*}
and
\begin{equation*}
\left\Vert f_{k}\right\Vert _{\infty }\leq  2^{2k}=(\text{supp}f_{k})^{-2},
\end{equation*}%
we conclude that $f_{k}$ is $ 1/2 $-atom, for every $k\in \mathbb{N}$

Moreover, if we use orthogonality of Walsh functions we get that 
\begin{equation*}
S_{2^{n}}\left(f_{k},x)\right)
\end{equation*}
\begin{equation*}
=\left\{
\begin{array}{ll}
0, & n=0,...,k, \\
\left( D_{2^{k+1}}(x)-D_{2^k}(x)\right), & n\geq k+1,
\end{array}%
\right.
\end{equation*}%
and
\begin{eqnarray*}
&&\sup\limits_{n\in \mathbb{N}}\left\vert S_{2^{n}}\left(f_{k},x\right)\right\vert \\
&=&\left\vert \left( D_{2^{k+1}}(x)-D_{2^k}(x)\right) \right\vert,
\end{eqnarray*}
where ${x}\in G$.

By combining first equality of Lemma \ref{lemma0} and Lemma \ref{lemma3.2.3} we obtain that
\begin{eqnarray*}
&&\left\Vert a_{k}\right\Vert _{H_{p}(G)} \\
&=&{2^{k}}
\left\Vert \sup\limits_{n\in \mathbb{N}}\left\vert S_{2^{n}}\left( D_{2^{k+1}}(x)-D_{2^k}(x)\right)\right\vert\right\Vert _{1/2} \\
&=&{2^{k}}
\left\Vert\left( D_{2^{k+1}}(x)-D_{2^k}(x)\right)\right\Vert _{1/2} \\
&=&{2^{k}}
\left\Vert D_{2^k}(x)\right\Vert _{1/2} \\
&\leq & {2^{k}}\cdot 2^{-k}\leq 1.
\end{eqnarray*}

It is easy to easy to show that
\begin{equation} \label{5.3.14a}
\widehat{f}_{m}\left( i\right) =\left\{
\begin{array}{l}
\text{ }2^{m},\text{ if }i=2^{m},...,2^{m+1}-1, \\
\text{ }0,\text{otherwise}%
\end{array}%
\right.
\end{equation}
and
\begin{equation}  \label{5.2.14}
S_{i}f_{m}=\left\{
\begin{array}{l}
2^{m}\left( D_{i}-D_{2^{m}}\right) ,\text{ \ if }i=2^{m}+1,...,2^{m+1}-1, \\
\text{ }f_{m},\text{ \ if }i\geq 2^{m+1}, \\
0,\text{ \ otherwise.}
\end{array}%
\right.
\end{equation}

Let $0<n<2^{m}.$ By using first equality of Lemma \ref{lemma0} we have that

\begin{eqnarray} \label{5.3.16b}
&&\left\vert \sigma _{n+2^{m}}f_{m}\right\vert \\ &=&\frac{1}{n+2^{m}}\left\vert
\overset{n+2^{m}}{\underset{j=2^{m}+1}{\sum }}S_{j}f_{m}\right\vert  \notag\\
&=&\frac{1}{n+2^{m}}\left\vert 2^{m}\overset{n+2^{m}}{\underset{j=2^{m}+1}{%
\sum }}\left( D_{j}-D_{2^{m}}\right) \right\vert  \notag \\
&=&\frac{1}{n+2^{m}}\left\vert 2^{m}\overset{n}{\underset{j=1}{\sum }}\left(
D_{j+2^{m}}-D_{2^{m}}\right) \right\vert  \notag \\
&=&\frac{1}{n+2^{m}}\left\vert 2^{m}\overset{n}{\underset{j=1}{\sum }}%
D_{j}\right\vert=\frac{2^{m}}{n+2^{m}}n\left\vert K_{n}\right\vert .  \notag
\end{eqnarray}

Let
$$n=\sum_{i=1}^{s}\sum_{k=l_{i}}^{m_{i}}2^{k},$$
where
$$0\leq l_{1}\leq m_{1}\leq l_{2}-2<l_{2}\leq m_{2}\leq ...\leq l_{s}-2<l_{s}\leq m_{s}.$$
By applying Lemma \ref{lemma7} and (\ref{5.3.16b}) we find that
\begin{equation*}
\left\vert\sigma_{n+2^{m}}f_{m}\left( x\right)\right\vert\geq c2^{2l_{i}}, \text{ \ \ \ where \ \ \ }x\in I_{l_{i}+1}\left(e_{l_{i}-1}+e_{l_{i}}\right).
\end{equation*}
Hence,
\begin{eqnarray*}
&&\int_{G}\left\vert \sigma _{n+2^{m}}f_{m}(x)\right\vert ^{1/2}d\mu \left(
x\right) \\
&\geq &\text{ }\underset{i=0}{\overset{s}{\sum }}\int_{I_{l_{i}+1}\left(
e_{l_{i}-1}+e_{l_{i}}\right) }\left\vert \sigma
_{n+2^{m}}f_{m}(x)\right\vert ^{1/2}d\mu \left( x\right) \\
&\geq &c\overset{s}{\underset{i=0}{\sum }}\frac{1}{2^{l_{i}}}2^{l_{i}}\geq cs\geq cV\left( n\right) .
\end{eqnarray*}

According to the second estimation of Lemma \ref{lemma2} we can conclude that
\begin{eqnarray*}
&&\underset{n\in \mathbf{\mathbb{N}}_{+}}{\sup }\underset{\left\Vert f\right\Vert _{H_{p}}\leq 1}{\sup}\frac{1}{n}\underset{k=1}{\overset{n}{\sum }}\left\Vert \sigma _{k}f\right\Vert_{1/2}^{1/2} \\
&\geq &\frac{1}{2^{m+1}}\overset{2^{m+1}-1}{\underset{k=2^{m}+1}{\sum }}
\left\Vert \sigma _{k}f_{m}\right\Vert _{1/2}^{1/2} \\
&\geq &\frac{c}{2^{m+1}}\underset{k=2^{m}+1}{\overset{2^{m+1}-1}{\sum }}
V\left( k-2^{m}\right) \\
&\geq& \frac{c}{2^{m+1}}\underset{k=1}{\overset{2^{m}-1}{\sum }}V\left( k\right)\geq c\log m\rightarrow \infty ,\text{ \ as \ }m\rightarrow \infty .
\end{eqnarray*}

The proof is complete.
\QED

\newpage

\section{Convergence and summability of partial sums with respect to the two-dimensional Walsh-Fourier series on the martingale Hardy spaces}

\subsection{Basic notations}

\text{ \qquad }  Let denote by $\overrightarrow{x}$ the two-dimensional vector $\overrightarrow{x}:=\left( x^{1},x^{2}\right)$ and by $G^{2}$ the direct product of two Walsh groups. Let $I_{n}^2:=I_{n}\left( 0\right)\times I_{n}\left( 0\right) $ for any $n\in \mathbb{N}$ and $\overline{I_{n}^2}:=G^2\backslash I_{n}^2$.

The norms (or quasi-norms) of the spaces of $L_{p}(G^{2})$ space is defined by

\begin{equation*}
\left\Vert f\right\Vert _{p}:=\left( \int_{G^{2}}  \left\vert f(\overrightarrow{x})\right\vert
^{p}d\mu (\overrightarrow{x})\right) ^{1/p}\qquad \left( 0<p<\infty \right).
\end{equation*}

The space $weak-L_{p}(G^{2}) $ consists of all functions $f$ for which

\begin{equation*}
\left\Vert f\right\Vert _{weak-L_{p}(G^2)}:=\underset{\lambda >0}{\sup}
\lambda \mu \left (\overrightarrow{x}\in{G^{2}}:  |f|>\lambda \right)^{1/p}<+\infty.
\end{equation*}

Two-dimensional Walsh system is defined by

\begin{equation*}
w_{n_{1},n_{2}}\left( x^{1},x^{2}\right):=w_{n_{1}}\left( x^{1}\right)w_{n_{2}}\left( x^{2}\right).
\end{equation*}

The two-dimensional Walsh system is orthonormal and complete  $L_{2}\left( {G}^2\right) $ (see
\cite{sws}).

For $f\in L_{1}\left( G^{2}\right) $ the following number
\begin{eqnarray*}
\widehat{f}\left(n_{1},n_{2}\right) :=\int\limits_{G^{2}}f({\overrightarrow{x}})w_{n_{1},n_{2}}({\overrightarrow{x}})d\mu({\overrightarrow{x}})
\end{eqnarray*}
is called $\left(n_{1},n_{2}\right) $-th Fourier coeficients of function  $f$.

$\left({n_{1},n_{2}}\right) $-th rectangular partial sum $S_{n_{1},n_{2}}$ of function $f$ is defined by:
\begin{equation*}
S_{n_{1},n_{2}}(f;{\overrightarrow{x}}):=\sum\limits_{i_{1}=0}^{n_{1}-1}\sum\limits_{i_{2}=0}^{n_{2}-1}
\widehat{f}\left( i_{1},i_{2}\right) w_{i_{1},i_{2}}\left( {\overrightarrow{x}}\right) .
\end{equation*}

$n$-th Marcinkiewicz (Marcinkiewicz-Fej\'er) means of the two-dimensional Walsh-Fourier series of function $f$ is defined by

\begin{equation*}
\mathcal{M}_{n}(f;{\overrightarrow{x}}):=\frac{1}{n}\sum
\limits_{k=0}^{n}S_{k,...,k}(f;{\overrightarrow{x}}).
\end{equation*}

Dirichlet and Marcinkiewicz kernels of the two-dimensional Walsh-Fourier series are defined by

\begin{equation*}
D_{n_{1},n_{2}}\left({\overrightarrow{x}}\right)=D_{n_{1}}\left( x^{1}\right)D_{n_{2}}\left( x^{2}\right).
\end{equation*}
and
\begin{equation*}
K_{n}({\overrightarrow{x}}):=
\frac{1}{n}\sum\limits_{k=0}^{n}D_{k,k}({\overrightarrow{x}}).
\end{equation*}

For the partial sums of the two-dimensional Walsh-Fourier let as define 
\begin{equation*}
S_{M}^{\left( 1\right) }f\left( x^1,x^2\right) :=\int\limits_{G}f\left(
s,x^2\right) D_{M}\left( x^1+s\right) d\mu \left( s\right)
\end{equation*}
and
\begin{equation*}
S_{N}^{\left( 2\right) }f\left( x^1,x^2\right) :=\int\limits_{G}f\left(
x^1,t\right) D_{N}\left(x^2+t\right) d\mu \left( t\right) .
\end{equation*}

For the partial sums of the two-dimensional Walsh-Fourier let define the following maximal operators $S_{\#}^{\ast}$ and $\widetilde{S}_{\#}^{\ast}$  by

\begin{equation*}
S_{\#}^{\ast}\left( f;x^{1},x^{2}\right): =\sup_{n\in \mathbb{N_+}} \left\vert S_{n,n}\left( f;x^{1},x^{2}\right) \right\vert
\end{equation*}
and
\begin{equation*}
\widetilde{S}_{\#}^{\ast}\left( f;x^{1},x^{2}\right)
:=\sup_{n\in \mathbb{N}}\left\vert S_{2^n,2^n}\left( f;x_{1},x_{2}\right)
\right\vert .
\end{equation*}

We define the  maximal operator  $\mathcal{M}^{\ast }$  and restricted maximal operator  $\widetilde{\mathcal{M}}_{\#}^{\ast}$ of Marcinkiewicz means by

\begin{equation*}
\mathcal{M}^{\ast }\left( f;x^{1},x^{2}\right) :=\sup_{n\in \mathbb{N_+}} \left\vert \mathcal{M}%
_n\left( f;x_{1},x_{2}\right) \right\vert
\end{equation*}
and
\begin{equation*}
\widetilde{\mathcal{M}}_{\#}^{\ast}\left( f;x^{1},x^{2}\right)
:=\sup_{n\in \mathbb{N}}\left\vert\mathcal{M}_{2^n}\left( f;x_{1},x_{2}\right)
\right\vert .
\end{equation*}

For the partial sums of the two-dimensional Walsh-Fourier let define the following weighted maximal operators
\begin{equation*}
\widetilde{\mathcal{M}}^{\ast }(f;x_{1},x_{2}):=\sup_{n\in \mathbb{N_+}}\frac{\left\vert
\mathcal{M}_{n}(f;x_{1},x_{2})\right\vert }{\log ^{3/2}\left( n+1\right) }
\end{equation*}
and
\begin{equation*}
\widetilde{\mathcal{M}}^{\ast ,p}\left( f;x_{1},x_{2}\right) =\sup_{n\in \mathbb{N_+}}\left\vert\frac{\mathcal{M}_{n}(f;x_{1},x_{2})}{(n+1)^{2/p-3}}\right\vert.\quad
\end{equation*}

The $\sigma $-algebra generated by the two-dimensional cubes
\begin{equation*}
I_{n}^{2}\left(\overrightarrow{x}\right):=I_{n}\left( x^{1}\right) \times I_{n}\left( x^{2}\right),
\end{equation*} is defined by $\digamma _{n}\left( n\in \mathbb{N}\right)$.

The conditional expectation operator with respect to
$\digamma _{n}\left( n\in \mathbb{N}\right) $
is denoted by $E_{n}$ and in our concrete case we have the following explicit expression for it:
\begin{eqnarray*}
E_{n}f(\overrightarrow{x})
=S_{2^{n},2^{n}}f\left( \overrightarrow{x}\right) \\ =\sum_{k_{1}=0}^{2^{n}-1}\sum_{k_{2}=0}^{2^{n}-1}
\widehat{f}\left( k_{1},k_{2}\right) w_{ k_{1},k_{2}}(\overrightarrow{x}) \\
=\frac{1}{\left| I_{n}^{2}\left( \overrightarrow{x}\right) \right| }\int_{I_{n}^{2}\left( \overrightarrow{x}\right)
}f(\overrightarrow{x})d\mu (\overrightarrow{x}),
\end{eqnarray*}
where $\left| I_{n}^{2}\left( \overrightarrow{x}\right) \right| =2^{-2n}$  denotes measure of cube $I_{n}^{2}\left(\overrightarrow{x}\right)$.

Sequence $f=\left(f_{n},\text{ }n\in \mathbb{N}\right) $ of functions $f_{n}\in L_{1}\left( G^2\right) $   is called dyadic martingales if (for details see \cite{sws})

$\left( i\right) $ $f_{n}$ is measurable with respect to $\sigma-$algebra $\digamma _{n}$, for any $n\in
\mathbb{N}$,

$\left( ii\right) $ $E_{n}f_{m}=f_{n}$ for every $n\leq m$.

The maximal function of martingale $f$ is defined by

\begin{equation*}
f^{\ast }=\sup_{n\in \mathbb{N}}\left\vert f_{n}\right\vert .
\end{equation*}

If $f\in L_{1}\left( G^2\right) ,$ then it is weell-known that the maximal operator is defined by

\begin{equation*}
f^{\ast }\left( \overrightarrow{x}\right) =\sup\limits_{n\in \mathbb{N}}\frac{1}{\mu \left(
I_{n}^{2}\left( \overrightarrow{x}\right) \right) }\left\vert \int\limits_{I_{n}^{2}\left( \overrightarrow{x}\right)
}f\left( \overrightarrow{u}\right) d\mu \left( \overrightarrow{u}\right) \right\vert.
\end{equation*}

For $0<p<\infty $ the one-parameter martingale Hardy space  $H_{p}\left( G^2\right)$ consist of all martingales for which
\begin{equation*}
\left\Vert f\right\Vert _{H_{p}\left( G^2\right)}:=\left\Vert f^{\ast }\right\Vert
_{p}<\infty .
\end{equation*}

Next, we define $p$-atoms, which are very important to characterize martingale Hardy spaces.
 
A function  $a$ is called a $p$-atom, if there exists an interval  $I^2$, such that  

\begin{equation*}
\int_{I^2}ad\mu =0, \\  \qquad \left\| a\right\| _{\infty }\leq \mu \left( I^2\right) ^{-1/p}, \\ \qquad \text{supp}\left( a\right) \subset I^2.
\end{equation*}

It is easy to check that for every martingale  $f=(f_{n},n\in \mathbb{N}) $ and for every $({k_{1},k_{2})}\in {\mathbb{N}}^{2}$ the limit
\begin{equation*}
\widehat{f}\left({k_{1},k_{2}}\right):=\lim_{{{n_{1},n_{2}}}\rightarrow \infty }\int_{{G}^{2}}f_{{{n_{1},n_{2}}}}\left(
{\overrightarrow{x}}\right) w_{{k_{1},k_{2}}}\left( {\overrightarrow{x}}\right) d\mu \left({\overrightarrow{x}}\right)
\end{equation*}
exists and it is called $\left({k_{1},k_{2}}\right)$-th Walsh-Fourier coefficients of $f$.

If $f\in L_{1}\left( G^2\right)$ and $(E_{n}f:n\in \mathbb{N})$ is regular martingale, then
\begin{eqnarray*}
\widehat{f}\left({k_{1},k_{2}}\right)=\int_{{G^2}}f\left( { {\overrightarrow{x}}}\right) w_{{k_{1},k_{2}}}({\overrightarrow{x}})
d\mu \left({{\overrightarrow{x}}}\right) 
=\widehat{f}\left({k_{1},k_{2}}\right) ,\text{ }{k_{1},k_{2}}\in \mathbb{N}.
\end{eqnarray*}

For the two-dimensional case modulus of continuity in $H_{p}(G^2)$ spaces can be defined as 
\begin{equation*}
\omega _{H_{p}\left( G^2\right)}\left( \frac{1}{2^{n}},f\right) :=\left\Vert
f-S_{2^{n},2^{n}}f\right\Vert _{H_{p}\left( G^2\right)}.
\end{equation*}

It is necessary to describe how can be understood difference $f-S_{2^{n},2^{n}}f$ where $f$ is martingale and $S_{2^{n},2^{n}}f$ is function. The following is true:

\begin{remark} \label{remark4.1.1}
Let $0<p\leq 1.$ Since
\begin{equation*}
S_{2^{n},2^{n}}f=f^{\left( n\right) }\in L_1(G^2),\text{ \ where }f=\left( f^{\left( n\right)}:n\in \mathbb{N}\right) \in H_{p}(G^2)
\end{equation*}%
and
\begin{eqnarray*}
&&\left( S_{2^{k},2^{k}}f^{\left( n\right) }:k \in \mathbb{N}\right) \\ &=&\left( S_{2^{k},2^{k}}S_{2^{n},2^{n}},k \in \mathbb{N}\right) \\
&=&\left( S_{2^{0},2^{0}}f,\ldots ,S_{2^{n-1},2^{n-1}}f,S_{2^{n},2^{n}}f,S_{2^{n},2^{n}}f,\ldots \right) \\
&=&\left( f^{\left( 0\right) },\ldots ,f^{\left( n-1\right) },f^{\left(
n\right) },f^{\left( n\right) },\ldots \right)
\end{eqnarray*}
we obtain that $  f-S_{2^{n},2^{n}}f $ is martingale for which
\begin{equation*}
\left( f-S_{2^{n},2^{n}}f\right) ^{\left( k\right) }=\left\{
\begin{array}{ll}
0, & k=0,.\ldots ,n, \\
f^{\left( k\right) }-f^{\left( n\right) }, & k\geq n+1,%
\end{array}%
\right.
\end{equation*}
and norm
$$ \left\Vert f-S_{2^{n},2^{n}}f\right\Vert _{H_{p}\left( G^2\right)}  $$   can be understood as $ H_p(G^2) $ norm of martingale
$$  f-S_{2^{n},2^{n}}f=(\left( f-S_{2^{n},2^{n}}f\right) ^{\left( k\right) }, \ \  k\in \mathbb{N}) $$
\end{remark}

\subsection{Auxillary lemmas}

\text{ \qquad }  In the following lemmas we investigate estimations of Marcinkiewicz means of the two-dimensional Walsh-Fourier series (see Lemma \ref{lemma8}-Lemma \ref{lemma11}).

Glukhov \cite{GLU} proved that the following is true:

\begin{lemma} \label{lemma8}
There exists an absolute constant $c$, such that
$$
\sup_{n \in \mathbb{N}} \int_{G^2} |K_n (x^1,x^2)| d\mu (x^1,x^2)\leq c.
$$
\end{lemma}

The following lemma is proved in \cite{gog1}:

\begin{lemma} \label{lemma9}
\label{lemma-gog1} Let $n\geq 2^{N}$, $(x^{1},x^{2})\in (I_{l^{1}}\backslash
I_{l^{1}+1})\times (I_{m^{2}}\backslash I_{m^{2}+1})$ and $0\leq l^{1}\leq
m^{2}<N$. Then
\begin{eqnarray*}
&&\int_{I_{N}\times I_{N}}|K_{n}(x^{1}+t^{1},x^{2}+t^{2})|d\mu (t^{1},t^{2})
\\
&\leq &\frac{c}{n2^{2N}}\left\{
2^{l^{1}-m^{2}}%
\sum_{r^{1}=l^{1}+1}^{m^{2}+1}2^{r^{1}}D_{2^{m^{2}+1}}(x^{1}+e_{l^{1}}+e_{r^{1}})\sum_{s=m^{2}+1}^{N}D_{2^{s}}(x^{2}+e_{m^{2}}+x_{m^{2}+1,s-1}^{1})\right.
\\
&&\left.
+2^{l^{1}+m^{2}}\sum_{s=l^{1}}^{m^{2}}%
\sum_{r^{1}=l^{1}+1}^{s}D_{2^{s}}(x^{1}+e_{l^{1}}+e_{r^{1}})\right\},
\end{eqnarray*}%
where  $$x_{i,j}:=\sum_{s=i}^{j}x_{s}e_{s}, \text{ \ \ } (x_{i,i-1}=0).$$
\end{lemma}

For our further investigation we need the following lemma (for details see \cite{gog1}):

\begin{lemma} \label{lemma10}
\label{lemma-gog2} Let $(x^{1},x^{2})\in I_{N}\times (I_{m^{2}}\backslash
I_{m^{2}+1})$ and $0\leq m^{2}<N.$ Then
\begin{equation*}
\int_{I_{N}\times I_{N}}|K_{n}(x^{1}+t^{1},x^{2}+t^{2})|d\mu
(t^{1},t^{2})
\end{equation*}

\begin{equation*}
\leq c\frac{2^{m^{2}}}{n2^{N}}%
\sum_{s=m^{2}}^{N-1}D_{2^{s}}(x^{2}+e_{m^{2}}),\text{ \ when \ }n>2^N.
\end{equation*}
\end{lemma}

We also need the following lemma proved by Goginava \cite{gog4}:

\begin{lemma} \label{lemma11}
Let
\begin{equation*}
x^{1}\in I_{4A}\left(
0,...,0,x_{4m}^{1}=1,0,...,0,x_{4l}^{1}=1,x_{4l+1}^{1},...,x_{4A-1}^{1}%
\right)
\end{equation*}%
and
\begin{equation*}
x^{2}\in I_{4A}\left(
0,...,0,x_{4l}^{2}=1,x_{4l+1}^{1}...,x_{4q-1}^{1},1-x_{4q}^{1},x_{4q+1}^{2},...,x_{4A-1}^{2}\right) .
\end{equation*}
Then
\begin{equation*}
n_{A-1}\left\vert K_{n_{A-1}}\left( x^{1},x^{2}\right) \right\vert \geq
2^{4q+4l+4m-3},
\end{equation*}%
where
$$n_{A}=2^{4A}+2^{4A-4}+...+2^{4}+2^{0}.$$
\end{lemma}

Hardy martingale spaces $H_{p}\left( G^2\right) $ have atomic decomposition for $0<p\leq 1$. The following is true (for details see \cite{S}, \cite{We1} and \cite{We3}):

\begin{lemma} \label{lemma3.2.4dim2}
A martingale $f=\left(f_{n},\text{ }n\in \mathbb{N}%
\right) $ belongs to $H_{p}(G^2)\left( 0<p\leq 1\right) $ if and only if there exist a sequence of p-atoms
 $\left( a_{k},\text{ }k\in \mathbb{N}\right) $ and sequence of real numbers 
 $\left( \mu _{k},\text{ }k\in \mathbb{N}\right) $ such that for every $n\in \mathbb{N}$

\begin{equation} \label{3.2.1}
\qquad \sum_{k=0}^{\infty }\mu _{k}S_{2^{n},2^{n}}a_{k}=f_{n}
\end{equation}%
and
\begin{equation*}
\qquad \sum_{k=0}^{\infty }\left\vert \mu _{k}\right\vert ^{p}<\infty.
\end{equation*}%
Moreover,
\begin{equation*}
\left\Vert f\right\Vert _{H_{p}(G^2)}\backsim \inf \left(
\sum_{k=0}^{\infty }\left\vert \mu _{k}\right\vert ^{p}\right) ^{1/p},
\end{equation*}
where infimum is taken over all decomposition of $ f $ of the form  (\ref{3.2.1}).

\end{lemma}

\begin{lemma} \label{lemma3.2.5dim2}
Let $0<p\leq 1$ and $T$ be $ \sigma  $-sublinear operator, such that 
\begin{equation*}
\int\limits_{G^2}\left\vert Ta(\overrightarrow{x})\right\vert ^{p}d\mu(\overrightarrow{x}) \leq
c_{p}<\infty,
\end{equation*}
for any  $p$-atom  $a$. Then there exists an absolute constant $ c_p, $ such that
\begin{equation} \label{3.2.5000dim2}
\left\Vert Tf\right\Vert _{p}\leq c_{p}\left\Vert f\right\Vert
_{H_{p}(G)}.
\end{equation}

In addition, if $T$  is bounded from $L_{\infty}(G^2)$ to $L_{\infty}(G^2)$, then we only have to prove that
\begin{equation*}
\int\limits_{\overset{-}{I^2}}\left\vert Ta\left(
{\overrightarrow{x}}\right)\right\vert ^{p}d\mu\left(
{\overrightarrow{x}}\right) \leq
c_{p}<\infty,
\end{equation*}
for every $p$-atom  $a$, where  $I^2$ denotes support of  $a$ and $\overline{I^2}:=G^2\backslash I^2$.
\end{lemma}

In special cases there exists simpler ways how to calculate $ H_p(G^2) $-norm
of martingale $ f\in H_p(G^2) $ (For details see e.g. \cite{S}, \cite{We1} and \cite{We3}):

\begin{lemma}  \label{lemma3.2.3dim2}
Let $g\in L_{1}\left( G^2\right)$ and $f:=(E_{n}g:n\in \mathbb{N})$ be regular martingale. Then  $H_{p}\left( G^2\right)\text{ }(0<p\leq1)$ norm is calculated by
\begin{equation*}
\left\Vert f\right\Vert _{H_{p}(G^2)}=\left\Vert \sup\limits_{n\in \mathbb{N}}|S_{2^n,2^n}g|\right\Vert _{p}.
\end{equation*}
\end{lemma}

Proofs of Lemma \ref{lemma3.2.8dim2} and Lemma \ref{lemma3.2.9} are proved in \cite{tep16}, \cite{tep18},  \cite{tep15}.

\begin{lemma} \label{lemma3.2.8dim2}
Let $0<p\leq 1$, $2^{k}\leq n<2^{k+1}$ and $S_{n,n}f$ be  $(n,n)$-th partial sum, where $f\in H_{p}(G^2)$. Then for any fixed $n\in \mathbb{N},$ we have the following estimation:
\begin{eqnarray*}
\left\Vert S_{n,n}f\right\Vert _{H_{p}(G^2)} 
&\leq& \left\Vert \sup_{0\leq l\leq k}\left\vert S_{2^{l},2^{l}}f\right\vert \right\Vert _{p}+\left\Vert S_{n,n}f\right\Vert _{p} \\
&&\leq \left\Vert \widetilde{S}_{\#}^{\ast }f\right\Vert _{p}+\left\Vert
S_{n,n}f\right\Vert _{p}.
\end{eqnarray*}
\end{lemma}

{\bf Proof}:
Let consider the following martingale
\begin{equation*}
f_{\#}:=\left( S_{2^{k},2^{k}}S_{n,n}f,\text{ }k\in\mathbb{N}_+\right)
=\left(S_{2^{0},2^{0}},S_{2^{k},2^{k}}f,...,\ \ \ S_{n,n}f,...,S_{n,n}f,...\right).
\end{equation*}%
By using Lemma \ref{lemma3.2.3dim2} we immediately get that
\begin{eqnarray*}
\left\Vert S_{n,n}f\right\Vert _{H_{p}(G^2)}^{p} 
&\leq & \left\Vert \sup_{0\leq l\leq k}\left\vert S_{2^{l},2^{l}}f\right\vert \right\Vert _{p}^{p}+\left\Vert
S_{n,n}f\right\Vert _{p}^{p} \\
&&\leq \left\Vert \widetilde{S}_{\#}^{\ast }f\right\Vert
_{p}^{p}+\left\Vert S_{n,n}f\right\Vert _{p}^{p}.
\end{eqnarray*}

The proof is complete.
\QED

\begin{lemma} \label{lemma3.2.9}
Let $0<p\leq1$, $2^{k}\leq n<2^{k+1}$ and $\mathcal{M}_{n}f$ be $n$-th Marcinkiewicz means, where $f\in H_{p}(G^2)$. Then for any fixed $n\in \mathbb{N},$ we get that
\begin{eqnarray*}
\left\Vert \mathcal{M}_{n}f\right\Vert _{H_{p}(G^2)}^{p}
&\leq& \left\Vert \sup_{0\leq l \leq k}\left\vert \mathcal{M}_{2^{l}}f\right\vert \right\Vert _{p}^{p}+\left\Vert \sup_{0\leq l\leq k}\left\vert S_{2^{l}}f\right\vert \right\Vert _{p}^{p}+\left\Vert \mathcal{M}
_{n}f\right\Vert _{p}^{p} \\
&\leq &\left\Vert \widetilde{\mathcal{M}}_{\#}^{\ast
}f\right\Vert _{p}^{p}+\left\Vert \widetilde{S}_{\#}^{\ast
}f\right\Vert _{p}^{p}+\left\Vert \mathcal{M}_{n}f\right\Vert _{p}^{p}.
\end{eqnarray*}
\end{lemma}

{\bf Proof}:
Let consider the following martingale
\begin{equation*}
f_{\#}=\left( S_{2^{k}}\mathcal{M}_{n}f,\text{ }k \in \mathbb{N}\right)
\end{equation*}%
\begin{equation*}
=\left( \frac{2^{0}\mathcal{M}_{2^{0}}}{n}+\frac{(n-2^{0})S_{2^{0}}f}{n},...,\frac{2^{k}\mathcal{M} _{2^{k}}f}{n}+\frac{(n-2^{k})S_{2^{k}}f}{n},
\mathcal{M}_{n}f,...,\mathcal{M}_{n}f,...\right).
\end{equation*}
According to Lemma \ref{lemma3.2.3dim2} we immediately get
\begin{eqnarray*}
\left\Vert \mathcal{M}_{n}f\right\Vert _{H_{p}(G^2)}^{p}
&\leq &\left\Vert \sup_{0\leq l\leq k}\left\vert \mathcal{M}_{2^{l}}f\right\vert \right\Vert _{p}^{p}+\left\Vert \sup_{0\leq l\leq k}\left\vert S_{2^{l},2^{l}}f\right\vert \right\Vert _{p}^{p}+\left\Vert
S_{n,n}f\right\Vert _{p}^{p} \\
&\leq & \left\Vert \widetilde{\mathcal{M}}_{\#}^{\ast }f\right\Vert
_{p}^{p}+\left\Vert \widetilde{S}_{\#}^{\ast }f\right\Vert
_{p}^{p}+\left\Vert \mathcal{M}_{n}f\right\Vert _{p}^{p}.
\end{eqnarray*}

The proof is complete.
\QED

\subsection{Strong convergence of partial sums with respect to the two-dimensional Walsh-Fourier series on the martingale Hardy spaces}

\text{ \qquad }In this section we investigate strong convergence of partial sums with respect to the two-dimensional Walsh-Fourier series on the martingale Hardy spaces when $ 0<p\leq 1 $ (see \cite{tep15}).
The following theorem is true:

\begin{theorem} \label{theorem6.1.1}
a) Let $0<p<1$ and $f\in H_{p}(G^2) $. Then there exists an absolute constant $ c_p, $ depending only on $ p, $ such that
\begin{equation*}
\sum\limits_{n=1}^{\infty }\frac{\left\| S_{n,n}f\right\| _{H_p(G^2)}^{p}}{n^{3-2p}}%
\leq c_{p}\left\| f\right\| _{H_{p}(G^2)}^{p}.
\end{equation*}

b) Let $0<p<1$ and $\Phi :\mathbb{N}\rightarrow [1,$ $\infty )$ be non-decreasing function, satisfying the condition
\begin{equation} \label{6.1.0}
\underset{n\rightarrow
\infty }{\lim }\Phi \left( n\right) =+\infty .
\end{equation}

Then there exists a martingale $f\in H_{p}(G^2)$, such that

\begin{equation*}
\underset{n=1}{\overset{\infty }{\sum }}\frac{\left\| S_{n,n}f\right\|
_{weak-L_{p}(G^2)}^{p}\Phi \left( n\right) }{n^{3-2p}}=\infty .
\end{equation*}

\end{theorem}

{\bf Proof}:
Suppose that
\begin{equation*}
\sum\limits_{n=1}^{\infty }\frac{\left\| S_{n,n}f\right\| _{p}^{p}}{n^{3-2p}}%
\leq c_{p}\left\| f\right\| _{H_{p}(G^2)}^{p}.
\end{equation*}

By combining Lemma \ref{lemma3.2.8} and inequality (\ref{1.Snnmax}) we can conclude that
\begin{eqnarray*}
{\sum_{n=1}^{\infty}}\frac{\left\Vert S_{n,n}f\right\Vert _{H_{p}(G^2)}^{p}}{n^{3-2p}} \leq{\sum_{n=1}^{\infty}}\frac{\left\Vert S_{n,n}f\right\Vert _{p}^{p}}{n^{3-2p}}+\left\Vert \widetilde{S}_{\#}^{\ast}f\right\Vert _{p}^{p} 
\leq  \left\Vert f\right\Vert _{H_{p}(G^2)}^{p}.
\end{eqnarray*}

According to Lemma \ref{lemma3.2.5dim2} we only have to prove that
\begin{equation} \label{6.1.main}
\sum\limits_{n=1}^{\infty }\frac{\left\| S_{n,n}a\right\| _{p}^{p}}{n^{3-2p}}%
\leq c_{p}<\infty ,
\end{equation}
for avery $p$-atom $a$.

Let $a$ be  $p$-atom with support $I_{N}\left( z^{1}\right) \times I_{N}\left( z^{2}\right) $, where $\mu \left(
I_{N}\right) =\mu \left( I_{N}\right) =2^{-N}$. Without loss the generality we may assume that $z^{1}=z^{2}=0.$

Let $\left(x^1,x^2\right)\in\overline{I}_{N}\times \overline{I}_{N}$. Then
\begin{equation*}
D_{2^{i}}\left(x^1+t^1\right) 1_{I_{N}}\left(t^1\right) =0, \text{ \ when \ } i\geq N
\end{equation*}
and
\begin{equation*}
D_{2^{i}}\left(x^2+t^2\right) 1_{I_{N}}\left(t^2\right) =0, \text{ \ when \ } i\geq N.
\end{equation*}
If we apply $w_{2^{j}}\left(x^i+t^i\right)=w_{2^{j}}\left(x^i\right),$ where $t^i\in I_{N}, \ \ i=1 \vee 2$ and $j<N$ according to both equality of Lemma \ref{lemma1} we get that
\begin{eqnarray*}
&&S_{n,n}a\left(x^1,x^2\right) \\
&=&\int\limits_{G\times G}a\left(t^1,t^2\right) D_{n}\left(x^1+t^1\right)
D_{n}\left(x^2+t^2\right) d\mu \left(t^1,t^2\right) \\
&=&\int\limits_{I_{N}\times I_{N}}a\left(t^1,t^2\right) D_{n}\left( x^1+t^1\right)
D_{n}\left(x^2+t^2\right) d\mu \left(t^1,t^2\right) \\
&=&\int\limits_{I_{N}\times I_{N}}a\left(t^1,t^2\right) w_{n}\left(
x^1+t^1+x^2+t^2\right) \sum\limits_{i=0}^{N-1}n_{i}w_{2^{i}}\left( x^1+t^1\right)
D_{2^{i}}\left( x^1+t^1\right) \\
&&\times \sum\limits_{j=0}^{N-1}n_{j}w_{2^{j}}\left(x^2+t^2\right)
D_{2^{j}}\left(x^2+t^2\right) d\mu \left(t^1,t^2\right) \\
&=& w_{n}\left( x^1\right) \sum\limits_{i=0}^{N-1}n_{i}w_{2^{i}}\left( x^2\right)
D_{2^{i}}\left( x^1\right) w_{n}\left( x^2\right)
\sum\limits_{j=0}^{N-1}n_{j}w_{2^{j}}\left( x^2\right) D_{2^{j}}\left( x^2\right)\\
&&\times \int\limits_{I_{N}\times I_{N}}a\left(t^1,t^2\right) w_{n}\left(
t^1+t^2\right) d\mu \left(t^1,t^2\right) \\
&=&w_{n}\left(x^1+x^2\right) \sum\limits_{i=0}^{N-1}n_{i}w_{2^{i}}\left(
x^1\right) D_{2^{i}}\left( x^1\right)
\sum\limits_{j=0}^{N-1}n_{j}w_{2^{j}}\left( x^2\right) D_{2^{j}}\left( x^2\right)
\\
&&\times \int\limits_{I_{N}}\left( \int\limits_{I_{N}}a\left( t^2+\tau
,t^2\right) d\mu \left( t^2\right) \right) w_{n}\left( \tau \right) d\mu \left(
\tau \right) \\
&=&w_{n}\left(x^1+x^2\right) \sum\limits_{i=0}^{N-1}n_{i}w_{2^{i}}\left(
x^1\right) D_{2^{i}}\left( x^1\right)
\sum\limits_{j=0}^{N-1}n_{j}w_{2^{j}}\left( x^2\right) D_{2^{j}}\left(
x^2\right) \int\limits_{I_{N}}\Phi \left( \tau \right) w_{n}\left( \tau
\right) d\mu \left( \tau \right) \\
&=&w_{n}\left( x^1+x^2\right) \sum\limits_{i=0}^{N-1}n_{i}w_{2^{i}}\left(
x^1\right) D_{2^{i}}\left( x^1\right)
\sum\limits_{j=0}^{N-1}n_{j}w_{2^{j}}\left(x^2\right) D_{2^{j}}\left(x^2\right) \widehat{\Phi }\left( n\right) ,
\end{eqnarray*}
where
\begin{equation*}
\Phi \left( \tau \right) =\int\limits_{I_{N}}a\left( t^i+\tau ,t^i\right) d\mu\left(t^i\right) \text{ \ and \ } i=1\vee 2.
\end{equation*}

Let $x\in I_{s}\backslash I_{s+1},$ where $ i=1\vee 2. $ By using again Lemma \ref{lemma1} can conclude that
\begin{equation*}
\overset{N-1}{\underset{i=0}{\sum }}D_{2^{i}}\left( x\right) \leq c2^{s}.
\end{equation*}

By using (\ref{2.1.2}) we have that
\begin{eqnarray} \label{6.1.dir3}
&&\int_{\overline{I}_{N}}\left( \overset{N-1}{\underset{i=0}{\sum }}%
D_{2^{i}}\left( x\right) \right) ^{p}d\mu \left( x\right)   \\
&\leq &c_{p}\underset{s=0}{\overset{N-1}{\sum }}\int_{I_{s}\backslash
I_{s+1}}2^{ps}d\mu \left( x\right)  \notag \\
&\leq &c_{p}\underset{s=0}{\overset{\infty }{\sum }}2^{\left( p-1\right) s}
\notag \\
&<&c_{p}<\infty ,\text{ \qquad }0<p<1.  \notag
\end{eqnarray}

From (\ref{6.1.dir3}) we get that
\begin{eqnarray*}
&&\sum\limits_{n=1}^{\infty }\frac{1}{n^{3-2p}}\int\limits_{\overline{I}%
_{N}\times \overline{I}_{N}}\left| S_{n,n}a\left(x^1,x^2\right) \right|
^{p}d\mu \left(x^1,x^2\right) \\
&\leq &\sum\limits_{n=1}^{\infty }\frac{\left| \widehat{\Phi }\left(
n\right) \right| ^{p}}{n^{3-2p}}\int\limits_{\overline{I}_{N}}\left(
\sum\limits_{i=0}^{N-1}D_{2^{i}}\left( x^1\right) \right) ^{p}d\mu \left(
x^1\right) \int\limits_{\overline{I}_{N}}\left(
\sum\limits_{i=0}^{N-1}D_{2^{i}}\left( x^2\right) \right) ^{p}d\mu \left(
x^2\right)  \\
&\leq & c_{p}\sum\limits_{n=1}^{\infty }\frac{\left| \widehat{\Phi }\left(
n\right) \right| ^{p}}{n^{3-2p}}.
\end{eqnarray*}

Let $n<2^{N}$. Since $w_{n}\left( \tau \right) =1,$ where $\tau \in I_{N}$ we obtain that
\begin{eqnarray*}
\widehat{\Phi }\left( n\right)&=&\int\limits_{I_{N}}\Phi \left( \tau
\right) w_{n}\left( \tau \right) d\mu \left( \tau \right) \\
&=&\int\limits_{I_{N}}\left( \int\limits_{I_{N}}a\left(t^2+\tau ,t^2\right)d\mu \left(t^2\right) \right) w_{n}\left( \tau \right) d\mu \left( \tau\right) \\
&=&\int\limits_{I_{N}\times I_{N}}a\left(t^1,t^2\right) d\mu \left(t^1,t^2\right)=0.
\end{eqnarray*}

So, we may assume that $n\geq 2^{N}.$ If we apply Holder's inequality we get that
\begin{eqnarray} \label{6.1.1p}
&&\sum\limits_{n=1}^{\infty }\frac{\left| \widehat{\Phi }\left( n\right)
\right| ^{p}}{n^{3-2p}}   \\
&\leq &\left( \sum\limits_{n=2^{N}}^{\infty }\left| \widehat{\Phi }\left(
n\right) \right| ^{2}\right) ^{p/2}\left( \sum\limits_{n=2^{N}}^{\infty }%
\frac{1}{n^{\left( 3-2p\right) \cdot \left( 2/\left( 2-p\right) \right) }}%
\right) ^{\left( 2-p\right) /2}  \notag \\
&\leq &\left( \frac{1}{2^{N\left( 2\left( 3-2p\right) /\left( 2-p\right)
-1\right) }}\right) ^{\left( 2-p\right) /2}\left( \int\limits_{G}\left| \Phi
\left( \tau \right) \right| ^{2}d\mu \left( \tau \right) \right) ^{p/2}
\notag \\
&\leq &\frac{c_{p}}{2^{N\left( 4-3p\right) /2}}\left(
\int\limits_{I_{N}}\left| \int\limits_{I_{N}}a\left( t^2+\tau ,t^2\right) d\mu
\left( t\right) \right| ^{2}d\mu \left( \tau \right) \right) ^{p/2}  \notag
\\
&\leq &\frac{c_{p}}{2^{N\left( 4-3p\right) /2}}\left\| a\right\| _{\infty
}^{p}\frac{1}{2^{Np/2}}\frac{1}{2^{Np}}  \notag \\
&\leq &\frac{c_{p}}{2^{N\left( 4-3p\right) /2}}2^{2N}\frac{1}{2^{3pN/2}}%
<c_{p}<\infty .  \notag
\end{eqnarray}

Let $\left(x^1,x^2\right) \in \overline{I}_{N}\times I_{N}$. Then
\begin{eqnarray*}
&&S_{n,n}a\left(x^1,x^2\right) \\
&=&w_{n}\left( x^1\right) \sum\limits_{j=0}^{N-1}n_{j}w_{2^{j}}\left( x^1\right)
D_{2^{j}}\left( x^1\right) \\
&&\times \int\limits_{G\times G}a\left(t^1,t^2\right) w_{n}\left(t^1\right)
D_{n}\left(x^2+t^2\right) d\mu \left(t^1,t^2\right) \\
&=&w_{n}\left( x^1\right) \sum\limits_{j=0}^{N-1}n_{j}w_{2^{j}}\left( x^1\right)D_{2^{j}}\left(x^1\right) \int\limits_{G}S_{n}^{\left( 2\right) }a\left(t^1,x^2\right) w_{n}\left(t^1\right) d\mu \left(t^1\right) \\
&=&w_{n}\left( x^1\right) \sum\limits_{j=0}^{N-1}n_{j}w_{2^{j}}\left( x^1\right)D_{2^{j}}\left( x^1\right) \widehat{S}_{n}^{\left( 2\right) }a\left(n,x^2\right) .
\end{eqnarray*}

By applying (\ref{6.1.dir3}) we get that
\begin{eqnarray*}
&&\sum\limits_{n=1}^{\infty }\frac{1}{n^{3-2p}}\int\limits_{\overline{I}%
_{N}\times I_{N}}\left| S_{n,n}a\left( x^1,x^2\right) \right| ^{p}d\mu \left(
x^1,x^2\right) \\
&\leq &\sum\limits_{n=1}^{\infty }\frac{1}{n^{3-2p}}\int\limits_{\overline{I}%
_{N}\times I_{N}}\left( \sum\limits_{j=0}^{N-1}D_{2^{j}}\left( x^1\right)
\left| \widehat{S}_{n}^{\left( 2\right) }a\left( n,x^2\right) \right| \right)^{p}d\mu \left( x^1,x^2\right) \\
&\leq &\sum\limits_{n=1}^{\infty }\frac{1}{n^{3-2p}}\int_{\overline{I}%
_{N}}\left( \overset{N-1}{\underset{i=0}{\sum }}D_{2^{i}}\left( x^1\right)
\right) ^{p}d\mu \left( x^1\right) \cdot \int\limits_{I_{N}}\left| \widehat{S}%
_{n}^{\left( 2\right) }a\left( n,x^2\right) \right| ^{p}d\mu \left(x^2\right) \\
&\leq &\sum\limits_{n=1}^{\infty }\frac{1}{n^{3-2p}}\int\limits_{I_{N}}%
\left| \widehat{S}_{n}^{\left( 2\right) }a\left( n,x^2\right) \right| ^{p}d\mu
\left(x^2\right) .
\end{eqnarray*}

Let $n<2^{N}$. By applyng definition of $p$-atom we get that
\begin{eqnarray*}
\widehat{S}_{n}^{\left( 2\right) }a\left( n,x^2\right)
&=&\int\limits_{G}\left( \int\limits_{G}a\left(t^1,t^2\right) D_{n}\left(
x^2+t^2\right) d\mu \left( t^2\right) \right) w_{n}\left( t^1\right) d\mu \left(t^1\right) \\
&=&D_{n}\left( x^2\right) \int\limits_{I_{N}\times I_{N}}a\left( t^1,t^2\right)
d\mu \left(t^1,t^2\right) =0.
\end{eqnarray*}

Therefore, we can suppose that $n\geq 2^{N}$. If follows that
\begin{eqnarray*}
&&\sum\limits_{n=1}^{\infty }\frac{1}{n^{3-2p}}\int\limits_{\overline{I}%
_{N}\times I_{N}}\left\vert S_{n,n}a\left( x^1,x^2\right) \right\vert ^{p}d\mu\left( x^1,x^2\right) \\
&\leq &\sum\limits_{n=2^{N}}^{\infty }\frac{1}{n^{3-2p}}\int\limits_{I_{N}}%
\left\vert \widehat{S}_{n}^{\left( 2\right) }a\left( n,x^2\right) \right\vert
^{p}d\mu \left(x^2\right)
\end{eqnarray*}%

Since
\begin{equation*}
\left\Vert S_{n}^{\left( 2\right) }a\left( n,x^2\right) \right\Vert _{2}\leq
c\left\Vert a\right\Vert _{2}
\end{equation*}%
if we use Holder's inequality we can conclude that
\begin{eqnarray*}
&&\int\limits_{I_{N}}\left\vert \widehat{S}_{n}^{\left(2\right)}a\left(
n,x^2\right)\right\vert ^{p}d\mu \left(x^2\right) \\
&\leq &\frac{c_{p}}{2^{N\left(1-p\right) }}\left(\int\limits_{I_{N}} \left\vert \widehat{S}_{n}^{\left(2\right)}a\left(n,x^2\right) \right\vert d\mu \left(x^2\right) \right) ^{p} \\
&=&\frac{c_{p}}{2^{N\left( 1-p\right)}}\left(\int\limits_{I_{N}} \left\vert\int\limits_{I_{N}}S_{n}^{\left(2\right)}a\left(t^1,x^2\right) w_{n}\left(t^1\right)d\mu \left(t^1\right)\right\vert d\mu \left(x^2\right)\right)^{p}\\
&=&\frac{c_{p}}{2^{N\left( 1-p\right) }}\left( \int\limits_{I_{N}}\left\vert
\int\limits_{I_{N}}\left( \int\limits_{I_{N}}a\left(t^1,t^2\right) D_{n}\left(x^2+t^2\right) d\mu \left( t^2\right) \right) w_{n}\left( t^1\right) d\mu \left(
t^1\right) \right\vert d\mu \left( x^2\right) \right) ^{p} \\
&\leq &\frac{c_{p}}{2^{N\left( 1-p\right) }}\left( \int\limits_{I_{N}}\left(
\int\limits_{I_{N}}\left\vert \int\limits_{I_{N}}a\left( t^1,t^2\right)
D_{n}\left( x^2+t^2\right) d\mu \left( t^2\right) \right\vert d\mu \left( x^2\right)
\right) d\mu \left( t^1\right) \right) ^{p} \\
&\leq &\frac{c_{p}}{2^{N\left( 1-p\right) }}\left( \frac{1}{2^{N/2}}%
\int\limits_{I_{N}}\left( \int\limits_{I_{N}}\left\vert
\int\limits_{I_{N}}a\left( t^1,t^2\right) D_{n}\left( x^2+t^2\right) d\mu \left(t^2\right) \right\vert ^{2}d\mu \left(x^2\right) \right) ^{1/2}d\mu \left(t^1\right) \right) ^{p} \\
&\leq &\frac{c_{p}}{2^{N\left( 1-p\right) }}\left( \frac{1}{2^{N/2}}%
\int\limits_{I_{N}}\left( \int\limits_{I_{N}}\left\vert a\left(t^1,t^2\right)
\right\vert ^{2}d\mu \left( t^2\right) \right) ^{1/2}d\mu \left(t^1\right)
\right) ^{p} \\
&\leq &\frac{c_{p}}{2^{N\left( 1-p\right) }}\left( \frac{\left\Vert
a\right\Vert _{\infty }}{2^{N/2}}\frac{1}{2^{N}}\frac{1}{2^{N/2}}\right)
^{p}\leq \frac{c_{p}}{2^{N\left( 1-p\right) }}\left( \frac{2^{2N/p}}{2^{2N}}%
\right) ^{p}\leq c_{p}2^{N\left( 1-p\right) }.
\end{eqnarray*}%

Hence,
\begin{eqnarray} \label{6.1.2p}
&&\sum\limits_{n=1}^{\infty }\frac{1}{n^{3-2p}}\int\limits_{\overline{I}%
_{N}\times I_{N}}{\left\vert S_{n,n}a\left(x^1,x^2\right) \right\vert}^p d\mu
\left(x^1,x^2\right)   \\ \notag
&\leq & c_{p}\sum\limits_{n=2^{N}}^{\infty }\frac{1}{n^{3-2p}}2^{N\left(
1-p\right) } \\ \notag
&\leq &\frac{c_{p}}{2^{N\left( 1-p\right) }}\leq c_{p}<\infty .
\notag
\end{eqnarray}

Analogously, we can prove that
\begin{equation} \label{6.3p}
\sum\limits_{n=1}^{\infty }\frac{1}{n^{3-2p}}\int\limits_{I_{N}\times
\overline{I}_{N}}\left| S_{n,n}a\left(x^1,x^2\right) \right|^{p}d\mu \left(x^1,x^2\right) \leq c_{p}<\infty .
\end{equation}

Let $\left(x^1,x^2\right) \in I_{N}\times I_{N}$. Then by the definition of $p$-atom we get that
\begin{eqnarray*}
&&\int\limits_{I_{N}\times I_{N}}\left| S_{n,n}a\left(x^1,x^2\right) \right|
^{p}d\mu \left(x^1,x^2\right) \\
&\leq &\frac{1}{2^{N\left( 2-p\right) }}\left( \int\limits_{I_{N}\times
I_{N}}\left| S_{n,n}a\left(x^1,x^2\right) \right| ^{2}d\mu \left( x^1,x^2\right)\right) ^{p/2} \\
&\leq &\frac{1}{2^{N\left( 2-p\right) }}\left( \int\limits_{I_{N}\times
I_{N}}\left| a\left( x^1,x^2\right) \right| ^{2}d\mu \left(x^1,x^2\right) \right)^{p/2} \\
&\leq &\frac{\left\| a\right\| _{\infty }^{p}}{2^{N\left( 2-p\right) }}\frac{1}{2^{Np}} \\
&\leq & c_{p}\frac{1}{2^{N\left( 2-p\right) }}2^{2N}\frac{1}{2^{Np}}%
\leq c_{p}<\infty .
\end{eqnarray*}
It follows that
\begin{eqnarray} \label{6.1.4p}
&&\sum\limits_{n=1}^{\infty }\frac{1}{n^{3-2p}}\int\limits_{I_{N}\times
I_{N}}\left| S_{n,n}a\left( x^1,x^2\right) \right| d\mu \left( x^1,x^2\right) \\
&\leq &c_{p}\sum\limits_{n=1}^{\infty }\frac{1}{n^{3-2p}}\leq c_{p}<\infty .
\notag
\end{eqnarray}

By combining (\ref{6.1.main}-\ref{6.1.4p}) we get that Theorem \ref{theorem6.1.1} is proved.

Let $0<p<1$ and $\Phi \left( n\right) $ satisfies condition \ref{6.1.0}. Then there exists increasing sequence of natural numbers $\left\{ \alpha _{k}:\text{ }k\in \mathbb{N}_+\right\} $ such that

\begin{equation*}
\alpha _{0}\geq 2
\end{equation*}
and
\begin{equation} \label{6.1.2}
\sum_{k=0}^{\infty }{\Phi ^{-p/4}\left( 2^{2\alpha _{k}}\right) }%
<\infty .
\end{equation}

Let $f=\left(f_{n},\text{ }n\in\mathbb{N}_+\right)$ be a martingale 
\begin{equation*}
f_{n}(x^{1},x^{2})={\sum_{k=1}^n}\lambda _{k}{a_{k}(x^{1},x^{2})}
\end{equation*}
where
\begin{equation*}
a_{k}(x^{1},x^{2})=2^{2\alpha _{k}(1/p-1)}\left(
D_{2^{2\alpha_{k}+1}}(x^{1})-D_{2\alpha_{k}}(x^{1})\right) \left(
D_{2^{2\alpha_{k}+1}}(x^{2})-D_{2^{2\alpha_{k}}}(x^{2})\right)
\end{equation*}
and
\begin{equation*}
\lambda _{k}={\Phi ^{-1/4}\left( 2^{2\alpha _{k}}\right)}.
\end{equation*}

By applying (\ref{6.1.2}) and Lemma \ref{lemma3.2.4dim2} we obtain that $f\in H_{p}(G^2).$

It is evident that

\begin{equation} \label{6.1.5}
\widehat{f}(i,j)=\left\{
\begin{array}{l}
\frac{2^{2\alpha _{k}\left( 2/p-2\right) }}{\Phi ^{1/4}\left( 2^{2\alpha
_{k}}\right) },\,\, \\
\text{ თუ }\left( i,j\right) \in \left\{
2^{2\alpha _{k}},...,\text{ ~}2^{2\alpha _{k}+1}-1\right\}^2, k\in\mathbb{N}_+,
\\
0,\text{ \thinspace } \\
\text{\thinspace \thinspace თუ \thinspace }\left(i,j\right) \notin \bigcup\limits_{k=1}^{\infty }\left\{ 2^{2\alpha _{k}},...,%
\text{ ~}2^{2\alpha _{k}+1}-1\right\}^2.%
\end{array}%
\right.
\end{equation}

Let  $2^{2\alpha _{k}}<n<2^{2\alpha _{k}+1}$. By combining (\ref{6.1.5}) and first equality of Lemma \ref{lemma0} we get that

\begin{eqnarray} \label{6.1.13d}
&&S_{n,n}f\left( x^1,x^2\right) \\
&=&\sum_{i=0}^{2^{\alpha _{k-1}+1}-1}\sum_{j=0}^{2^{\alpha _{k-1}+1}-1}%
\widehat{f}(i,j)w_{i}\left( x^1\right) w_{j}\left( x^2\right)  \notag \\
&+&\sum_{i=2^{\alpha _{k}}}^{n-1}\sum_{j=2^{\alpha _{k}}}^{n-1}\widehat{f}%
(i,j)w_{i}\left( x^1\right) w_{j}\left( x^2\right)  \notag \\
&=&\sum_{\eta =0}^{k-1}\sum_{i=2^{\alpha _{\eta }}}^{2^{\alpha _{\eta
}+1}-1}\sum_{j=2^{\alpha _{\eta }}}^{2^{\alpha _{\eta }+1}-1}\widehat{f}%
(i,j)w_{i}\left( x^1\right) w_{j}\left( x^2\right)  \notag \\
&+&\sum_{i=2^{\alpha _{k}}}^{n-1}\sum_{j=2^{\alpha _{k}}}^{n-1}\widehat{f}%
(i,j)w_{i}\left( x^1\right) w_{j}\left( x^2\right)  \notag \\
&=&\sum_{\eta =0}^{k-1}\sum_{i=2^{\alpha _{\eta }}}^{2^{\alpha _{\eta
}+1}-1}\sum_{j=2^{\alpha _{\eta }}}^{2^{\alpha _{\eta }+1}-1}\frac{2^{\alpha
_{\eta }\left( 2/p-2\right) }}{\Phi ^{1/4}\left( 2^{\alpha _{\eta }}\right) }%
w_{i}\left( x^1\right) w_{j}\left( x^2\right)  \notag \\
&+&\sum_{i=2^{\alpha _{k}}}^{n-1}\sum_{j=2^{\alpha _{k}}}^{n-1}\frac{%
2^{\alpha _{k}\left( 2/p-2\right) }}{\Phi ^{1/4}\left( 2^{\alpha
_{k}}\right) }w_{i}\left( x^1\right) w_{j}\left( x^2\right)  \notag \\
&=&\sum_{\eta =0}^{k-1}\frac{2^{\alpha _{\eta }\left( 2/p-2\right) }}{\Phi
^{1/4}\left( 2^{\alpha _{\eta }}\right) }\left( D_{2^{\alpha _{\eta
}+1}}\left( x^1\right) -D_{2^{\alpha _{\eta }}}\left( x^1\right) \right) \left(
D_{2^{\alpha _{\eta }+1}}\left( x^2\right) -D_{2^{\alpha _{\eta }}}\left(
x^2\right) \right)  \notag \\
&+&\frac{2^{\alpha _{k}\left( 2/p-2\right) }}{\Phi ^{1/4}\left( 2^{\alpha
_{k}}\right) }\left( D_{_{n}}\left( x^1\right) -D_{2^{\alpha _{k}}}\left(
x^1\right) \right) \left( D_{n}\left( x^2\right) -D_{2^{\alpha _{k}}}\left(
x^2\right) \right)  \notag \\
&=&I+II.  \notag
\end{eqnarray}

Let $\left( x^1,x^2\right) \in \left( G\backslash I_{1}\right)^2$ and $n\ $ is odd number. Since $n-2^{2\alpha _{k}}\
$ is also odd, according to both equality of Lemma \ref{lemma1} we get that

\begin{eqnarray} \label{6.1.13a}
\left| II\right|
&=&\frac{2^{2\alpha _{k}\left( 2/p-2\right) }}{\Phi
^{1/4}\left( 2^{2\alpha _{k}}\right) }\left| w_{2^{2\alpha _{k}}}\left(
x^1\right) D_{n-2^{2\alpha _{k}}}\left( x^1\right) w_{2^{2\alpha _{k}}}\left(
x^2\right) D_{n-2^{2\alpha _{k}}}\left( x^2\right) \right|  \\
&=&\frac{2^{2\alpha _{k}\left( 2/p-2\right) }}{\Phi ^{1/4}\left( 2^{2\alpha
_{k}}\right) }\left| w_{2^{2\alpha _{k}}}\left( x^1\right) w_{n-2^{2\alpha
_{k}}}\left( x^1\right) D_{_{1}}\left( x^1\right) w_{2^{2\alpha _{k}}}\left(x^2\right) w_{n-2^{2\alpha _{k}}}\left(x^2\right) D_{_{1}}\left(x^2\right) \right|
\notag \\
&=&\frac{2^{2\alpha _{k}\left( 2/p-2\right) }}{\Phi ^{1/4}\left( 2^{2\alpha
_{k}}\right) }.  \notag
\end{eqnarray}

If we apply again second equality of Lemma \ref{lemma1} according to condition  $\alpha _{n}\geq 2$ $\left( n\in \mathbb{N}\right)$ for $I$ we can conclude that

\begin{equation} \label{6.1.13b}
I=\sum_{\eta =0}^{k-1}\frac{2^{2\alpha _{k}\left( 2/p-2\right) }}{\Phi
^{1/4}\left( 2^{2\alpha _{\eta }}\right) }\left( D_{2^{2\alpha _{\eta
}+1}}\left( x^1\right) -D_{2^{2\alpha _{\eta }}}\left( x^1\right) \right) \left(
D_{2^{\alpha _{\eta }+1}}\left( x^2\right) -D_{2^{\alpha _{\eta }}}\left(
x^2\right) \right) =0.
\end{equation}
Hence,
\begin{eqnarray} \label{6.1.13}
&&\left\| S_{n,n}f\left( x^1,x^2\right) \right\| _{weak-L_{p}(G^2)}   \\
&\geq &\frac{2^{2\alpha _{k}\left( 2/p-2\right) }}{2\Phi ^{1/4}\left(
2^{2\alpha _{k}}\right) }\left( \mu \left\{ \left( x^1,x^2\right) \in \left(G\backslash I_{1}\right)^2:\left|
S_{n,n}f\left( x^1,x^2\right) \right| \geq \frac{2^{2\alpha _{k}\left(
2/p-2\right) }}{2\Phi ^{1/4}\left( 2^{2\alpha _{k}}\right) }\right\} \right)
^{1/p}  \notag \\
&\geq &\frac{2^{2\alpha _{k}\left( 2/p-2\right) }}{2\Phi ^{1/4}\left(
2^{2\alpha _{k}}\right) }\left|( G\backslash I_{1}\right)^2| \geq \frac{c_{p}2^{2\alpha _{k}\left(
2/p-2\right) }}{\Phi ^{1/4}\left( 2^{2\alpha _{k}}\right) }.  \notag
\end{eqnarray}

By using (\ref{6.1.13}) we get that

\begin{eqnarray} \label{6.1.14}
&&\underset{n=1}{\overset{2^{2\alpha _{k}+1}-1}{\sum }}\frac{\left\|
S_{n,n}f\right\| _{weak-L_{p}(G^2)}^{p}\Phi \left( n\right) }{n^{3-2p}} \\
&\geq &\underset{n=2^{2\alpha _{k}}+1}{\overset{2^{2\alpha _{k}+1}-1}{\sum }}%
\frac{\left\| S_{n,n}f\right\| _{weak-L_{p}(G^2)}^{p}\Phi \left( n\right) }{%
n^{3-2p}}  \notag \\
&\geq &c_{p}\Phi \left( 2^{2\alpha _{k}}\right) \underset{n=2^{\alpha
_{k}-1}+1}{\overset{2^{\alpha _{k}}-1}{\sum }}\frac{\left\|
S_{2n+1,2n+1}f\right\| _{weak-L_{p}(G^2)}^{p}}{\left( 2n+1\right)^{3-2p}}  \notag
\\
&\geq &c_{p}\Phi \left( 2^{2\alpha _{k}}\right) \frac{2^{2\alpha _{k}\left(
1-p\right) }}{\Phi ^{1/4}\left( 2^{2\alpha _{k}}\right) }\underset{%
n=2^{2\alpha _{k}-1}+1}{\overset{2^{2\alpha _{k}}-1}{\sum }}\frac{1}{\left(
2n+1\right) ^{3-2p}}  \notag \\
&\geq &c_{p}\Phi ^{3/4}\left( 2^{2\alpha _{k}}\right) \rightarrow \infty ,%
\text{ \quad as \quad }k\rightarrow \infty .  \notag
\end{eqnarray}

Theorem is proved.
\QED

\subsection{Strong convergence of Marcinkiewicz means with respect to the two-dimensional Walsh-Fourier series on the martingale Hardy spaces}

\text{ \qquad } In this section we consider strong convergence of Marcinkiewicz means with respect to the two-dimensional Walsh-Fourier series in the martingale Hardy spaces for $p=2/3$ (for details see Nagy and Tephnadze \cite{tep16}).

\begin{theorem} \label{th7.1.1}
Let $f\in H_{2/3}\left( G^{2}\right)$. Then there exists an absolute constant $c,$ such that
\begin{equation*}
\frac{1}{\log n}\sum_{m=1}^{n}\frac{\left\Vert \mathcal{M}_{m}f\right\Vert
_{H_{2/3}(G^2)}^{2/3}}{m}\leq c\left\Vert f\right\Vert _{H_{2/3}(G^2)}^{2/3}.
\end{equation*}%
\end{theorem}

{\bf Proof}:
Suppose that
\begin{equation} \label{7.1.5.4}
\frac{1}{\log n}{\sum_{m=1}^{n}}\frac{\left\Vert \mathcal{M}%
_{m}f\right\Vert _{2/3}^{2/3}}{m}\leq c\left\Vert f\right\Vert
_{H_{2/3}(G^2)}^{2/3}.
\end{equation}

By combining Lemma \ref{lemma3.2.9} and inequalities (\ref{1.Snnmax}), (\ref{1.Mnnmax}), (\ref{7.1.5.4}) we get that
\begin{eqnarray}  \label{6.1.eq-1}
&&\frac{1}{\log n}{\sum_{m=1}^{n}}\frac{\left\Vert \mathcal{M}%
_{m}f\right\Vert _{H_{2/3}(G^2)}^{2/3}}{m} \\ \notag
&\leq&\frac{1}{\log n}{\sum_{m=1}^{n}}\frac{\left\Vert \mathcal{M}%
_{m}f\right\Vert _{2/3}^{2/3}}{m} \\ \notag
&+&\left\Vert \widetilde{\mathcal{M}}_{\#}^{\ast}f\right\Vert _{2/3}^{2/3}+\left\Vert \widetilde{S}_{\#}^{\ast}f\right\Vert _{2/3}^{2/3}\\  \notag
&\leq& \left\Vert f\right\Vert _{H_{2/3}(G^2)}^{2/3}.
\end{eqnarray}

Since $\mathcal{M}_{n}$ is (see Lemma \ref{lemma8}) bounded from $L_{\infty }(G^2)$ to $L_{\infty }(G^2)$, if we use Lemma \ref{lemma3.2.5dim2} we only have to prove that
\begin{equation*}
\frac{1}{\log n}\sum_{m=1}^{n}\frac{\left\Vert \mathcal{M}_{m}a\right\Vert
_{2/3}^{2/3}}{m}<c<\infty,
\end{equation*}
for every $2/3$-atom $a$.

Let $a$ be $2/3$-atom with support $I^{2}$, where $\mu
(I^{2})=2^{-2N}$. Without loss the generality we may assume that $I^{2}:=I_{N}^2$. It is easy to show that ${\mathcal{M}}_{n}a=0$
for $n\leq 2^{N}.$ So, we may assume that $n>2^{N}$.

We can write that
\begin{eqnarray*}
&&\frac{1}{\log n}\sum_{m=1}^{n} \frac{\left\Vert \mathcal{M}%
_{m}a\right\Vert _{2/3}^{2/3}}{m} \\
& \leq& \frac{1}{\log n}\sum_{m=2^{N}}^{n}%
\frac{\left\Vert {\mathcal{M}_{m}a}\right\Vert _{2/3}^{2/3}}{m} \\
& \leq & \frac{1}{\log n}\sum_{m=2^{N}}^{n}\int_{I_{N}\times I_{N}}\frac{%
\left\vert {\mathcal{M}_{m}a}\right\vert ^{2/3}}{m}d\mu \\
&+&\frac{1}{\log n}%
\sum_{m=2^{N}}^{n}\int_{I_{N}\times \overline{I_{N}}}\frac{\left\vert {%
\mathcal{M}_{m}a}\right\vert ^{2/3}}{m}d\mu \\
& +&\frac{1}{\log n}\sum_{m=2^{N}}^{n}\int_{\overline{I_{N}}\times {I_{N}}}%
\frac{\left\vert {\mathcal{M}_{m}a}\right\vert ^{2/3}}{m}d\mu \\
&+&\frac{1}{%
\log n}\sum_{m=2^{N}}^{n}\int_{\overline{I_{N}}\times \overline{I_{N}}}\frac{%
\left\vert {\mathcal{M}_{m}a}\right\vert ^{2/3}}{m}d\mu \\
& =:&I_{1}+I_{2}+I_{3}+I_{4}.
\end{eqnarray*}%
By applying Lemma \ref{lemma8} we have that
\begin{eqnarray*}
I_{1} &\leq &\frac{1}{\log n}\sum_{m=2^{N}}^{\infty }\int_{I_{N}\times I_{N}}\frac{\left\vert \mathcal{M}_{m}a\right\vert ^{2/3}}{m}d\mu \\
&\leq &\frac{1}{\log n}\sum_{m=2^{N}}^{\infty }\frac{1}{m}\left\Vert
a\right\Vert _{\infty }^{2/3}/2^{2N} \\
&\leq &\frac{1}{\log n}\sum_{m=2^{N}}^{n}%
\frac{1}{m}<c<\infty.
\end{eqnarray*}%

Now, we estimate $I_{2}$. Set
\begin{equation*}
J_{t}:=I_{t}\backslash I_{t+1},\text{ \ \ } (t\in {\mathbb{N}}).
\end{equation*}

We introduce $\overline{I_{N}}$ and $J_{m^{2}}$ as the following disjoint union:
\begin{equation} \label{7.2.gio2}
\overline{I_{N}}=\bigcup_{m^{2}=0}^{N-1}J_{m^{2}},\quad
J_{m^{2}}=\bigcup_{q^{2}=m^{2}+1}^{N}I_{N}^{m^{2},q^{2}},
\end{equation}%
where
\begin{equation*}
I_{N}^{m^{2},q^{2}}:=%
\begin{cases}
I_{q^{2}+1}(0,...,0,x_{m^{2}}=1,0,...,0,x_{q^{2}}=1), & \text{ where }\text{ \ \ }
m^{2}<q^{2}<N, \\
I_{N}(0,...,0,x_{m^{2}}=1,0,...,0), & \text{ where }\text{ \ \ }q^{2}=N.%
\end{cases}%
\end{equation*}

Let $(x^{1},x^{2})\in I_{N}\times I_{N}^{m^{2},q^{2}}.$ According to Lemma \ref
{lemma10} we can conclude that
\begin{eqnarray*}
&&{|{\mathcal{M}}_{n}a(x^{1},x^{2})|}  \\
&\leq &\Vert a\Vert _{\infty
}\int_{I_{N}\times I_{N}}|K_{n}(x^{1}+t^{1},x^{2}+t^{2})|d\mu (t^{1},t^{2})
\\
&\leq &c2^{3N}\frac{2^{m^{2}}}{n2^{N}}%
\sum_{s=m^{2}}^{q^{2}}D_{2^{s}}(x^{2}+e_{m^{2}}) \\
&\leq &\frac{c2^{2N+m^{2}}}{n}\sum_{s=m^{2}}^{q^{2}}2^{s} \\
&\leq & \frac{
c2^{2N+m^{2}+q^{2}}}{n}.
\end{eqnarray*}%
Hence,
\begin{eqnarray*}
I_{2} &\leq &\frac{c2^{4N/3}}{\log n}\sum_{m=2^{N}}^{n}\sum_{m^{2}=0}^{N-1}%
\sum_{q^{2}=m^{2}+1}^{N}\int_{I_{N}\times I_{N}^{m^{2},q^{2}}}\frac{%
\left\vert {\mathcal{M}_{m}a}\right\vert ^{2/3}}{m}d\mu \\
&\leq &\frac{c2^{4N/3}}{\log n}\sum_{m=2^{N}}^{n}\sum_{m^{2}=0}^{N-1}%
\sum_{q^{2}=m^{2}+1}^{N}\int_{I_{N}\times I_{N}^{m^{2},q^{2}}}\frac{%
2^{2(m^{2}+q^{2})/3}}{m^{5/3}}d\mu \\
&\leq &\frac{c2^{4N/3}}{\log n}\sum_{m=2^{N}}^{n}\sum_{m^{2}=0}^{N-1}%
\sum_{q^{2}=m^{2}+1}^{N}\frac{2^{2(m^{2}+q^{2})/3}}{m^{5/3}}2^{-N-q^{2}} \\
&\leq &\frac{c2^{N/3}}{\log n}\sum_{m=2^{N}}^{\infty }\frac{1}{m^{5/3}}%
\sum_{m^{2}=0}^{N-1}2^{2m^{2}/3}\sum_{q^{2}=m^{2}+1}^{N}2^{-q^{2}/3} \\
&\leq &\frac{c2^{N/3}}{\log n}\sum_{m=2^{N}}^{\infty }\frac{2^{N/3}}{m^{5/3}} \\
&\leq &\frac{c2^{2N/3}}{\log n}\sum_{m=2^{N}}^{\infty }\frac{1}{m^{5m/3}} \\
&\leq&
\frac{c}{N}.
\end{eqnarray*}

Analogously, we can prove that $I_{3}\leq c<\infty $.

Next we prove boundedness of $I_{4}$. If we apply (\ref{7.2.gio2}) we get that
\begin{eqnarray*}
I_{4} &\leq &\frac{1}{\log n}\sum_{m=2^{N}}^{n}\sum_{l^{1}=0}^{N-1}%
\sum_{m^{2}=0}^{l^{1}-1}\int_{J_{l^{1}}\times J_{m^{2}}}\frac{\left\vert {%
\mathcal{M}_{m}a}\right\vert ^{2/3}}{m}d\mu \\
&+&\frac{1}{\log n}\sum_{m=2^{N}}^{n}\sum_{l^{1}=0}^{N-1}%
\sum_{m^{2}=l^{1}}^{N-1}\int_{J_{l^{1}}\times J_{m^{2}}}\frac{\left\vert {\mathcal{M}_{m}a}\right\vert ^{2/3}}{m}d\mu  \\
&=:& I_{4,1}+I_{4,2}.
\end{eqnarray*}%

Let consider $I_{4,2}$ (Analogously we can estimate $I_{4,1}$). For $(x^{1},x^{2})\in J_{l^{1}}\times J_{m^{2}}$ if we apply Lemma \ref{lemma9} we obtain that
\begin{eqnarray*}
&&{|{\mathcal{M}}_{n}a(x^{1},x^{2})|} \\
&&\leq \Vert a\Vert _{\infty
}\int_{I_{N}\times I_{N}}|K_{n}(x^{1}+t^{1},x^{2}+t^{2})|d\mu (t^{1},t^{2})
\\
&\leq &\frac{2^{N+l^{1}-m^{2}}}{n}%
\sum_{r^{1}=l^{1}+1}^{m^{2}+1}2^{r^{1}}D_{2^{m^{2}+1}}(x^{1}+e_{l^{1}}+e_{r^{1}})\sum_{s=m^{2}+1}^{N}D_{2^{s}}(x^{2}+e_{m^{2}}+x_{m^{2}+1,s-1}^{1})
\\
&&+\frac{2^{N+l^{1}+m^{2}}}{n}\sum_{s=l^{1}}^{m^{2}}%
\sum_{r^{1}=l^{1}+1}^{s}D_{2^{s}}(x^{1}+e_{l^{1}}+e_{r^{1}}).
\end{eqnarray*}%
It is evident that
\begin{eqnarray*}
&&\int_{J_{l^{1}}\times
J_{m^{2}}}D_{2^{m^{2}+1}}^{2/3}(x^{1}+e_{l^{1}}+e_{r^{1}})D_{2^{s}}^{2/3}(x^{2}+e_{m^{2}}+x_{m^{2}+1,s-1}^{1})d\mu (x^{1},x^{2}) \\
&&\leq c2^{2s/3-m^{2}/3-l^{1}}\leq c2^{-\left( m^{2}+s\right) /3}
\end{eqnarray*}%
and
\begin{eqnarray*}
&&\int_{J_{l^{1}}\times
J_{m^{2}}}D_{2^{s}}^{2/3}(x^{1}+e_{l^{1}}+e_{r^{1}})d\mu (x^{1},x^{2})\leq
c2^{2s/3-m^{2}-l^{1}}\leq c2^{-m^{2}-s/3}.
\end{eqnarray*}%
Hence, 
\begin{eqnarray*}
&&\int_{J_{l^{1}}\times J_{m^{2}}}\left\vert {\mathcal{M}_{m}a}\right\vert
^{2/3}d\mu  \\
&&\leq \frac{c2^{2\left( N+l^{1}-m^{2}\right) /3}}{m^{2/3}}%
\sum_{r^{1}=l^{1}+1}^{m^{2}+1}\sum_{s=m^{2}+1}^{N}2^{2r^{1}/3}\times \\
&&\times \int_{J_{l^{1}}\times
J_{m^{2}}}D_{2^{m^{2}+1}}^{2/3}(x^{1}+e_{l^{1}}+e_{r^{1}})D_{2^{s}}^{2/3}(x^{2}+e_{m^{2}}+x_{m^{2}+1,s-1}^{1})d\mu\left( x^{1},x^{2}\right)
\\
&&+\frac{c2^{2\left( N+l^{1}+m^{2}\right) /3}}{m^{2/3}}%
\sum_{s=l^{1}}^{m^{2}}\sum_{r^{1}=l^{1}+1}^{s}\int_{J_{l^{1}}\times
J_{m^{2}}}D_{2^{s}}^{2/3}(x^{1}+e_{l^{1}}+e_{r^{1}})d\mu (x^{1},x^{2}) \\
&\leq & \frac{c2^{2(N+l^{1}-m^{2})/3}}{m^{2/3}}%
\sum_{r^{1}=l^{1}+1}^{m^{2}+1}2^{2r^{1}/3}\sum_{s=m^{2}+1}^{N}2^{-\left(
m^{2}+s\right) /3} +\frac{c2^{2(N+l^{1}+m^{2})/3}}{m^{2/3}}%
\sum_{s=l^{1}}^{m^{2}}\sum_{r^{1}=l^{1}+1}^{s}2^{-m^{2}-s/3} \\
& \leq& \frac{c2^{2(N+l^{1}-m^{2})/3}}{m^{2/3}}2^{-2m^{2}/3}%
\sum_{r^{1}=l^{1}+1}^{m^{2}+1}2^{2r^{1}/3}+\frac{c2^{2\left(
N+l^{1}+m^{2}\right) /3}}{m^{2/3}}\sum_{s=l_{1}}^{m^{2}}\left(
s-l_{1}-1\right) 2^{-m^{2}-s/3} \\
&\leq & \frac{c2^{2\left( N+l^{1}-m^{2}\right) /3}}{m^{2/3}}+\frac{%
c2^{\left( 2N+l^{1}-m^{2}\right) /3}}{m^{2/3}}.
\end{eqnarray*}
and
\begin{equation*}
\begin{split}
I_{4,2}& \leq \frac{c}{\log n}\sum_{m=2^{N}}^{n}\frac{1}{m}%
\sum_{l^{1}=0}^{N-1}\sum_{m^{2}=l^{1}}^{N-1}\frac{2^{2(N+l^{1}-m^{2})/3}+2^{%
\left( 2N-m^{2}+l^{1}\right) /3}}{m^{2/3}} \\
& \leq \frac{c2^{2N/3}}{\log n}\sum_{m=2^{N}}^{n}\frac{1}{m^{5/3}}%
\sum_{l^{1}=0}^{N-1}1\leq c.
\end{split}%
\end{equation*}
\QED

\begin{corollary} \label{cor7.1.1}
Let $f\in H_{2/3}\left( G^{2}\right) .$ Then
\begin{equation*}
\lim_{n\rightarrow \infty }\frac{1}{\log n}\sum_{m=1}^{n}\frac{\left\Vert
\mathcal{M}_{m}f-f\right\Vert _{H_{2/3}\left( G^{2}\right)}^{2/3}}{m}=0
\end{equation*}%
and
\begin{equation*}
\lim_{n\rightarrow \infty }\frac{1}{\log n}\sum_{m=1}^{n}\frac{\left\Vert
\mathcal{M}_{m}f\right\Vert _{H_{2/3}\left(G^{2}\right)}^{2/3}}{m} =\left\Vert f\right\Vert_{H_{2/3}\left( G^{2}\right)}^{2/3}.
\end{equation*}
\end{corollary}

\subsection{Modulus of continuity and convergence in norm of Marcinkiewicz means with respect to the two-dimensional Walsh-Fourier series on the martingale Hardy spaces}

\text{ \qquad }In this section we investigate necessary and sufficient conditions for modulus of continuity, which provide convergence in norm of Marcinkiewicz means with respect to the two-dimensional Walsh-Fourier series in $ H_{2/3} $-norm (For details see \cite{tep17}).

\begin{theorem} \label{theorem7.2.1}
 a) Let $f\in H_{2/3}(G^2)$ and
\begin{equation}  \label{7.2.10A}
\omega_{H_{2/3}(G^2)} \left( \frac{1}{2^{k}},f\right)=o\left( \frac{1}{k^{3/2}}
\right) ,\text{\quad as \quad }k\rightarrow \infty.
\end{equation}
Then
\begin{equation*}
\left\Vert \mathcal{M}_{n}f -f\right\Vert
_{H_{2/3}(G^2)}\rightarrow 0, \text{ \quad as \quad }n\rightarrow \infty .
\end{equation*}

b) There exists a martingale $f\in H_{2/3}(G^2),$ such that
\begin{equation*}
\omega _{H_{2/3}(G^2)} \left( \frac{1}{2^{2^{k}}},f\right)=O\left( \frac{1}{%
2^{3k/2}}\right) , \text{ \quad as \quad }k\rightarrow \infty
\end{equation*}
and
\begin{equation*}
\left\Vert \mathcal{M}_{n}f-f\right\Vert_{2/3}\nrightarrow 0%
\text{ \quad as \quad }n\rightarrow \infty.
\end{equation*}
\end{theorem}

{\bf Proof}:
In \cite{nagy} it was proved that (see inequality (\ref{1.c22})) the following inequality is true:
\begin{equation} \label{7.2.5.3.0}
\left\Vert \mathcal{M}_{n}f\right\Vert _{2/3}\leq c\log^{2/3} \left( n+1\right) \left\Vert f\right\Vert_{H_{2/3}(G^2)}.
\end{equation}

If we apply inequality (\ref{1.Mnnmax}) and Lemma \ref{lemma3.2.9} according to (\ref{7.2.5.3.0}) we get the following estimation
\begin{eqnarray} \label{7.2.5.3}
&&\left\Vert \mathcal{M}_{n}f\right\Vert _{H_{2/3}(G^2)}^{2/3}\\\notag
&&\leq \left\Vert \mathcal{M}_{n}f\right\Vert _{2/3}^{2/3}+\left\Vert \widetilde{\mathcal{M}}_{\#}^{\ast}f\right\Vert _{2/3}^{2/3}+\left\Vert \widetilde{S}_{\#}^{\ast}f\right\Vert _{2/3}^{2/3}\\\notag
&&\leq c\log \left( n+1\right) \left\Vert f\right\Vert _{H_{2/3}(G^2)}^{2/3}+c\left\Vert f\right\Vert _{H_{2/3}(G^2)}^{2/3}\\\notag
&&\leq c\log\left(n+1\right)\left\Vert f\right\Vert _{H_{2/3}(G^2)}^{2/3}.  \notag
\end{eqnarray}

Let $2^{N}<n\leq 2^{N+1}.$ If we use (\ref{7.2.5.3}) by simple calculations we have that
\begin{eqnarray*}
&&\left\Vert \mathcal{M}_{n}f-f\right\Vert _{H_{2/3}(G^2)}^{2/3} \\
&&\leq \left\Vert
\mathcal{M}_{n}f-\mathcal{M}_{n}S_{2^{N},2^{N}}f\right\Vert _{H_{2/3}(G^2)}^{2/3}\\
&&+\left\Vert \mathcal{M}_{n}S_{2^{N},2^{N}}f-S_{2^{N},2^{N}}f\right\Vert
_{H_{2/3}(G^2)}^{2/3}\\
&&+\left\Vert S_{2^{N},2^{N}}f-f\right\Vert _{H_{2/3}(G^2)}^{2/3} \\
&&=\left\Vert \mathcal{M}_{n}\left( S_{2^{N},2^{N}}f-f\right) \right\Vert
_{H_{2/3}}^{2/3} \\
&&+\left\Vert \mathcal{M}_{n}S_{2^{N},2^{N}}f-S_{2^{N},2^{N}}f\right\Vert
_{H_{2/3}(G^2)}^{2/3}\\
&&+\left\Vert S_{2^{N},2^{N}}f-f\right\Vert _{H_{2/3}(G^2)}^{2/3} \\
&&\leq c\left( \log \left( n+1\right) +1\right) \omega_{H_{2/3}(G^2)} ^{2/3}\left( \frac{1}{2^{N}},f\right) \\
&&+\left\Vert \mathcal{M}_{n}S_{2^{N},2^{N}}f-S_{2^{N},2^{N}}f\right\Vert
_{H_{2/3}(G^2)}^{2/3}.
\end{eqnarray*}

Let $2^{N}<n\leq 2^{N+1}.$ Then it is evident that
\begin{eqnarray*}
&&\mathcal{M}_{n}S_{2^{N},2^{N}}f-S_{2^{N},2^{N}}f \\
&=&\frac{1}{n}
\sum_{k=0}^{2^{N}}S_{k,k}S_{2^{N},2^{N}}f \\
&+&\frac{1}{n}
\sum_{k=2^{N}+1}^{n}S_{k,k}S_{2^{N},2^{N}}f-S_{2^{N},2^{N}}f \\
& =&\frac{1}{n}\sum_{k=0}^{2^{N}}S_{k,k}f \\
&+&\frac{\left( n-2^{N}\right)
S_{2^{N},2^{N}}f}{n}-S_{2^{N},2^{N}}f \\
&=&\frac{2^{N}}{n}\left( \mathcal{M}_{2^{N}}f-S_{2^{N},2^{N}}f\right) \\
&=&\frac{2^{N}}{n}\left( S_{2^{N},2^{N}}\mathcal{M}
_{2^{N}}f-S_{2^{N},2^{N}}f\right) \\
&=&\frac{2^{N}}{n}S_{2^{N},2^{N}}\left( \mathcal{M}_{2^{N}}f-f\right) .
\end{eqnarray*}
By combining (\ref{1.Snn0}) and (\ref{1.Mnn}) we get that
\begin{eqnarray} \label{7.2.10}
&&\left\Vert \mathcal{M}_{n}S_{2^{N},2^{N}}f-S_{2^{N},2^{N}}f\right\Vert
_{H_{2/3}(G^2)}^{2/3} \\
&&\leq \left( \frac{2^{N}}{n}\right) ^{2/3}\hspace{-10pt}%
\left\Vert S_{2^{N},2^{N}}\left( \mathcal{M}_{2^{N}}f-f\right) \right\Vert
_{H_{2/3}(G^2)}^{2/3}  \notag   \\
&&\leq
\left\Vert S_{2^{N},2^{N}}\left( \mathcal{M}_{2^{N}}f-f\right) \right\Vert
_{H_{2/3}(G^2)}^{2/3}  \notag   \\
&&\leq \left\Vert \mathcal{M}_{2^{N}}f-f\right\Vert
_{H_{2/3}(G^2)}^{2/3}\rightarrow 0,\text{ where }k\rightarrow \infty,  \notag
\end{eqnarray}

Hence, we immediately get that if
\begin{equation*}
\omega_{H_{2/3}(G^2)}\left( \frac{1}{2^{n}},f\right)=o\left( \frac{1}{n^{3/2}}%
\right) ,\text{ as }n\rightarrow \infty ,
\end{equation*}%
then
\begin{equation*}
\left\Vert \mathcal{M}_{n}f-f\right\Vert _{H_{2/3}(G^2)}\rightarrow 0,\text{ as
}n\rightarrow \infty .
\end{equation*}

Now, prove part b) of Theorem \ref{theorem7.2.1}. Let
\begin{equation*}
a_{i}(x^{1},x^{2})=2^{2^{i}}\left(
D_{2^{2^{i}+1}}(x^{1})-D_{2^{2^{i}}}(x^{1})\right) \left(
D_{2^{2^{i}+1}}(x^{2})-D_{2^{2^{i}}}(x^{2})\right)
\end{equation*}%
and
\begin{equation*}
f_{n}(x^{1},x^{2})={\sum_{i=1}^n}\frac{a_{i}(x^{1},x^{2})}{%
2^{3i/2}}.
\end{equation*}

Since
\begin{equation*}
{\sum_{i=1}^{\infty}}\left(\frac{1}{2^{3i/2}}\right)^{2/3}<c<\infty
\end{equation*}
\begin{eqnarray*}
&& S_{2^{n},2^{n}}a_{k}(x^{1},x^{2})=
\begin{cases}
a_{k}(x^{1},x^{2}), & \text{ \ if \ }2^{k}\leq n, \\
0, & \text{ \ if \ }2^{k}>n
\end{cases}%
\end{eqnarray*}%
and
\begin{eqnarray*}
&&\text{supp }a_{k} =I_{2^{k}}^{2}, \\
&&\int_{I_{2^{k}}^{2}}a_{k}d\mu = 0, \\
&&\left\Vert a_{k}\right\Vert _{\infty } \leq \mu (\text{supp }a_{k})^{-3/2},
\end{eqnarray*}%
by using Lemma \ref{lemma3.2.4dim2}, we conclude that $f\in H_{2/3}.$

On the other hand, if we apply Remark \ref{remark4.1.1} we immediately get that
\begin{eqnarray*}
&& f-S_{2^{n},2^{n}}f \\
&=&\left( f^{\left( 1\right) }-S_{2^{n},2^{n}}f^{\left( 1\right)
},...,f^{\left( n\right) }-S_{2^{n},2^{n}}f^{\left( n\right) },...,f^{\left(
n+k\right) }-S_{2^{n},2^{n}}f^{\left( n+k\right) },...\right) \\
&=&\left( 0,...,0,f^{\left( n+1\right) }-f^{\left( n\right) },...,f^{\left(
n+k\right) }-f^{\left( n\right) },...\right) \\
&=&\left( 0,...,0,...,\sum_{i=\log n+1}^{\log n+k} \frac{a_{i}(x)}{%
2^{3i/2}},...\right) ,\text{ \ }k\in \mathbb{N}_{+},
\end{eqnarray*}%
Hence,
\begin{eqnarray*}
&&\omega_{H_{2/3}(G^2)}\left( \frac{1}{2^{n}},f\right)\leq\sum_{i=\left[ \log n\right] }^\infty\frac{1}{2^{3i/2}}=O\left( \frac{1}{n^{3/2}}%
\right) \text{ \quad as \quad }k \rightarrow \infty .
\end{eqnarray*}%
where $\left[ \log n\right] $ denotes integer part of  $\log n$.

Set
\begin{eqnarray*}
n_{2^{A-2}}&=& 2^{4\cdot 2^{A-2}}+2^{4\cdot
2^{A-2}-4}+...+2^{4}+2^{0} \\
&=& 2^{2^{A}}+2^{2^{A}-4}+...+2^{4}+2^{0}.
\end{eqnarray*}
If we use Lemma \ref{lemma11} we get that
\begin{eqnarray} \label{7.2.nn}
&&\mathcal{M}_{n_{2^{k-2}}}f-f \\ \notag
&=&\frac{2^{2^{k}}\mathcal{M}_{2^{2^{k}}}f}{%
n_{2^{k-2}}}+\frac{1}{n_{2^{k-2}}}\sum_{j=2^{2^{k}}+1}^{n_{2^{k-2}}}S_{j,j}f \\ \notag
&-&\frac{2^{2^{k}}f}{n_{2^{k-2}}}-\frac{n_{2^{k-2}-1}f}{n_{2^{k-2}}}.  \\ \notag
\end{eqnarray}%

It is evident that

\begin{equation} \label{7.2.35}
\widehat{f}(i,j)=%
\begin{cases}
\frac{2^{2^{k}}}{2^{3k/2}}, & \text{ თუ }\left( i,j\right) \in \left\{
2^{2^{k}},...,2^{2^{k}+1}-1\right\} ^{2},\ k\in\mathbb{N} \\
0, & \text{ თუ }\left( i,j\right) \notin \bigcup\limits_{k=0}^{\infty
}\left\{ 2^{2^{k}},...,2^{2^{k}+1}-1\right\} ^{2}.%
\end{cases}
\end{equation}%

Let $2^{2^{k}}<j\leq n_{2^{k-1}}.$ Since $%
w_{v+2^{2^{k}}}=w_{2^{2^{k}}}w_{v}, $ when $v<2^{2^{k}}$, if we apply (\ref{7.2.35}) and first equality of Lemma \ref{lemma0} we obtain that
\begin{eqnarray*}
&& S_{j,j} f\left( x^{1},x^{2}\right) \\
&=& S_{2^{2^{k}},2^{2^{k}}}f\left(
x^{1},x^{2}\right) \\ &+&\sum_{v=2^{2^{k}}}^{j-1}\sum_{s=2^{2^{k}}}^{j-1}\widehat{
f}(v,s)w_{v,s}\left( x^{1},x^{2}\right) \\
&=& S_{2^{2^{k}},2^{2^{k}}}f\left( x^{1},x^{2}\right) \\
&+&\frac{2^{2^{k}}}{
2^{3k/2}}\sum_{v=0}^{j-2^{2^{k}}-1}\sum_{s=0}^{j-2^{2^{k}}-1}w_{v+2^{2^{k}}}
\left( x^{1}\right) w_{s+2^{2^{k}}}\left( x^{2}\right) \\
&=& S_{2^{2^{k}},2^{2^{k}}}f\left( x^{1},x^{2}\right) \\
&+&\frac{
2^{2^{k}}w_{2^{2^{k}}}\left( x^{1}\right) w_{2^{2^{k}}}\left( x^{2}\right) }{
2^{3k/2}}\sum_{v=0}^{j-2^{2^{k}}-1}\sum_{s=0}^{j-2^{2^{k}}-1}w_{v}\left(
x^{1}\right) w_{s}\left( x^{2}\right) \\
&=& S_{2^{2^{k}},2^{2^{k}}}f\left( x^{1},x^{2}\right) \\
&+&\frac{
2^{2^{k}}w_{2^{2^{k}}}\left( x^{1}\right) w_{2^{2^{k}}}\left( x^{2}\right)
D_{_{j-2^{2^{k}},j-2^{2^{k}}}}\left( x^{1},x^{2}\right) }{2^{3k/2}}.
\end{eqnarray*}
Hence,
\begin{eqnarray*}
&&\frac{1}{n_{2^{k-2}}}\sum_{j=2^{2^{k}}+1}^{n_{2^{k-2}}}S_{j,j}f\left(
x^{1},x^{2}\right) \\
&=&\frac{n_{2^{k-2}-1}S_{2^{2^{k}},2^{2^{k}}}f\left( x^{1},x^{2}\right) }{%
n_{2^{k-2}}}\\
&&+\frac{2^{2^{k}}w_{2^{2^{k}}}\left( x^{1}\right)
w_{2^{2^{k}}}\left( x^{2}\right) }{n_{2^{k-2}}2^{3k/2}}%
\sum_{j=1}^{n_{2^{k-2}-1}}D_{_{j,j}}\left( x^{1},x^{2}\right) \\
&=&\frac{n_{2^{k-2}-1}S_{2^{2^{k}},2^{2^{k}}}f\left( x^{1},x^{2}\right) }{n_{2^{k-2}}}\\
&&+\frac{2^{2^{k}}w_{2^{2^{k}}}\left( x^{1}\right)
w_{2^{2^{k}}}\left( x^{2}\right) n_{2^{k-2}-1}K_{n_{2^{k-2}-1}}\left(
x^{1},x^{2}\right) }{n_{2^{k-1}}2^{3k/2}}.
\end{eqnarray*}

By applying (\ref{7.2.nn}) we have that

\begin{eqnarray} \label{7.2.a11}
&&\Vert \mathcal{M}_{n_{2^{k-2}}}f-f\Vert _{2/3}^{2/3} \\
&&\geq\frac{c}{2^{k}}%
\Vert {n_{2^{k-2}-1}}K_{n_{2^{k-2}-1}}\Vert _{2/3}^{2/3}   \notag\\
&&-\left( \frac{2^{2^{k}}}{n_{2^{k-2}}}\right) ^{2/3}\Vert \mathcal{M}%
_{2^{2^{k}}}f-f\Vert _{2/3}^{2/3}  \notag \\
&&-\left( \frac{n_{2^{k-2}-1}}{n_{2^{k-2}}}\right) ^{2/3}\Vert
S_{2^{2^{k}},2^{2^{k}}}f-f\Vert _{2/3}^{2/3} \notag\\
&&\geq\frac{c}{2^{k}}
\Vert {n_{2^{k-2}-1}}K_{n_{2^{k-2}-1}}\Vert _{2/3}^{2/3}   \notag\\
&&-\Vert \mathcal{M}_{2^{2^{k}}}f-f\Vert _{2/3}^{2/3} \notag\\
&&-\Vert S_{2^{2^{k}},2^{2^{k}}}f-f\Vert _{2/3}^{2/3}.  \notag
\end{eqnarray}

Set
\begin{eqnarray*}
&& x^{1}\in I_{2^{k-2}}^{m,l} \\
&=:& I_{2^{k-2}}\left(
0,...,0,x_{4m}^{1}=1,0,...,0,x_{4l}^{1}=1,x_{4l+1}^{1},...,x_{2^{k-2}-1}^{1}%
\right)
\end{eqnarray*}%
and
\begin{eqnarray*}
&&x^{2}\in J_{2^{k-2}}^{l,q} \\
&=:& I_{2^{k-2}}\left(
0,...,0,x_{4l}^{2}=1,x_{4l+1}^{1}...,x_{4q-1}^{1},1-x_{4q}^{1},x_{4q+1}^{2},...,x_{2^{k-2}-1}^{2}\right) .
\end{eqnarray*}

According to Lemma \ref{lemma11} we get that
\begin{equation*}
n_{_{2^{k-2}-1}}\left\vert K_{n_{2^{k-2}-1}}\left( x^{1},x^{2}\right)
\right\vert \geq 2^{4q+4l+4m-3}.
\end{equation*}
Hence,
\begin{eqnarray*}
&&\int_{G}(n_{_{2^{k-2}-1}}|K_{n_{2^{k-2}-1}}(x^1,x^2)|)^{2/3}d\mu (x^1,x^2)\\
&&\geq c\sum_{m=1}^{2^{k-2}-3}\sum_{l=m+1}^{2^{k-2}-2}\sum_{q=l+1}^{2^{k-2}-1}%
\sum_{x_{4l+1}^{1}=0}^{1}...\sum_{x_{_{2^{k-2}-1}}^{1}=0}^{1}%
\sum_{x_{4q+1}^{2}=0}^{1}...\sum_{x_{_{2^{k-2}-1}}^{2}=0}^{1} \\
&&\int_{I_{2^{k-2}}^{m,l}\times J_{2^{k-2}}^{l,q}}\hspace{0cm}({n_{2^{k-2}-1}%
}|K_{n_{2^{k-2}-1}}(x^1,x^2)|)^{2/3}d\mu(x^1,x^2) \\
&& \geq
c\sum_{m=1}^{2^{k-2}-3}\sum_{l=m+1}^{2^{k-2}-2}\sum_{q=l+1}^{2^{k-2}-1}%
\sum_{x_{4l+1}^{1}=0}^{1}...\sum_{x_{_{2^{k-2}-1}}^{1}=0}^{1}%
\sum_{x_{4q+1}^{2}=0}^{1}...\sum_{x_{_{2^{k-2}-1}}^{2}=0}^{1} \\
&&\mu \left( I_{2^{k-2}}^{m,l}\times J_{2^{k-2}}^{l,q}\right) 2^{\left(
8q+8l+8m\right) /3} \\
\end{eqnarray*}
\begin{eqnarray*}
&& \geq c\sum_{m=1}^{2^{k-2}-3}\sum_{l=m+1}^{2^{k-2}-2}\sum_{q=l+1}^{2^{k-2}-1}2^{%
\left( 8q+8l+8m\right) /3}2^{2^{k-2}-4l}2^{2^{k-2}-4q}\left( \frac{1}{%
2^{2^{k-2}}}\right) ^{2} \\
&& \geq c\sum_{m=1}^{2^{k-2}-3}2^{8m/3}\sum_{l=m+1}^{2^{k-2}-2}2^{-4l/3}%
\sum_{q=l+1}^{2^{k-2}-1}2^{-4q/3}  \\
&&\geq  c\sum_{m=1}^{2^{k-2}-3}1\geq c2^{k}.
\end{eqnarray*}

By combining (\ref{1.Snn}), (\ref{1.Mnn}) and (\ref{7.2.a11}) we can conclude that
\begin{equation*}
{\limsup_{k\rightarrow \infty } }\Vert \mathcal{M}
_{n_{2^{k-2}}}f-f\Vert _{2/3}\geq c>0.
\end{equation*}
Theorem \ref{theorem7.2.1} is proved.
\QED

\end{document}